\documentclass[12pt]{article}
\usepackage{amsthm,amsfonts,amssymb,epsfig,graphics,amsmath,amsbsy,color,enumerate,tikz}

\theoremstyle{plain}
\newtheorem{theorem}{Theorem}
\newtheorem{lemma}{Lemma}
\newtheorem{proposition}{Proposition}
\newtheorem{remark}{Remark}
\newtheorem{definition}{Definition}

\newtheorem{claim}{Claim}

\newcommand{\mas}{\operatorname{Mas}}

\newcommand{\dom}{\operatorname{dom}}
\newcommand{\ran}{\operatorname{ran}}
\newcommand{\colspan}{\operatorname{colspan}}
\newcommand{\AC}{\operatorname{AC}}
\newcommand{\ess}{\operatorname{ess}}
\newcommand{\pt}{\operatorname{pt}}
\newcommand{\ind}{\operatorname{ind}}

\numberwithin{equation}{section}
\numberwithin{lemma}{section}
\numberwithin{theorem}{section}
\numberwithin{remark}{section}
\numberwithin{claim}{section}
\numberwithin{corollary}{section}
\numberwithin{proposition}{section}
\numberwithin{definition}{section}
\numberwithin{condition}{section}

\DeclareFontFamily{OT1}{pzc}{}
\DeclareFontShape{OT1}{pzc}{m}{it}{<-> s * [1.100] pzcmi7t}{}
\DeclareMathAlphabet{\mathpzc}{OT1}{pzc}{m}{it}

\title{Renormalized Oscillation Theory for Linear 
Hamiltonian Systems on $[0, 1]$ via the Maslov Index}

\author{Peter Howard and Alim Sukhtayev}

\begin{document}

\maketitle

\begin{abstract} 
Working with a general class of regular linear Hamiltonian systems, 
we show that renormalized oscillation 
results can be obtained in a natural way through consideration
of the Maslov index associated with appropriately 
chosen paths of Lagrangian subspaces of $\mathbb{C}^{2n}$.   
We verify that our applicability class includes Dirac 
and Sturm-Liouville systems, as well as a system arising 
from differential-algebraic equations for which the spectral
parameter appears nonlinearly.
\end{abstract}

\section{Introduction}\label{introduction}

For values $\lambda$ in some interval $I \subset \mathbb{R}$, we consider linear Hamiltonian systems  
\begin{equation} \label{hammy}
J y' = \mathbb{B} (x; \lambda) y; \quad x \in (0,1), \quad y (x; \lambda) \in \mathbb{C}^{2n},
\quad n \in \{1, 2, \dots \},
\end{equation}
where $J$ denotes the standard symplectic matrix 
\begin{equation*}
J 
=
\begin{pmatrix}
0_n & - I_n \\
I_n & 0_n
\end{pmatrix},
\end{equation*}
and an important feature of the analysis is that $\mathbb{B} (x; \lambda)$ is 
allowed to depend nonlinearly on the spectral parameter $\lambda$.
We assume throughout that for each $\lambda \in I$, $\mathbb{B} (x; \lambda)$
is a measurable matrix-valued function of $x$, self-adjoint for 
a.e. $x \in (0, 1)$, for which there exists a dominating function
$b_0 \in L^1 ((0,1), \mathbb{R})$ so that for each $\lambda \in I$,
$\|\mathbb{B} (x; \lambda)\| \le b_0 (x)$ for a.e. $x \in (0,1)$. 
Moreover, we assume $\mathbb{B}$ is differentiable in $\lambda$, 
and that there exists a dominating function
$b_1 \in L^1 ((0,1), \mathbb{R})$ so that for each $\lambda \in I$,
$\|\mathbb{B}_{\lambda} (x; \lambda)\| \le b_1 (x)$ for a.e. $x \in (0,1)$. 
(Here, $\|\cdot\|$ denotes any matrix norm.) For convenient reference, we refer to these 
basic assumptions as Assumptions {\bf (A)}. 

We consider two types of self-adjoint boundary conditions, 
{\it separated} and {\it generalized}.

\medskip
\noindent
{\bf (BC1).} We express separated self-adjoint boundary conditions as 
\begin{equation*}
\alpha y(0; \lambda) = 0; \quad \beta y(1; \lambda) = 0,
\end{equation*}
where we assume
\begin{equation*}
\begin{aligned}
\alpha &\in \mathbb{C}^{n \times 2n}, 
\quad \operatorname{rank} \alpha = n, \quad \alpha J \alpha^* = 0; \\
\beta &\in \mathbb{C}^{n \times 2n}, 
\quad \operatorname{rank} \beta = n, \quad \beta J \beta^* = 0.
\end{aligned}
\end{equation*}

\medskip
\noindent
{\bf (BC2).} We express general self-adjoint boundary conditions as
\begin{equation*}
\Theta
\begin{pmatrix}
y(0; \lambda) \\ y(1; \lambda) 
\end{pmatrix} = 0; \quad \Theta  \in \mathbb{C}^{2n \times 4n},
\quad \operatorname{rank} \Theta = 2n, \quad \Theta \mathcal{J}_{4n} \Theta^* = 0,
\end{equation*}
where 
\begin{equation*}
\mathcal{J}_{4n} := 
\begin{pmatrix}
-J & 0 \\
0 & J
\end{pmatrix}.
\end{equation*}
Here, and throughout, we use the superscript $*$ to denote adjoint. We 
emphasize that the matrices $\alpha$ and $\beta$ in {\bf (BC1)} are taken 
independent of $\lambda$, as is the matrix $\Theta$ in {\bf (BC2)}. These
restrictions on $\alpha$, $\beta$, and $\Theta$ are not necessary for 
the approach, and are assumed rather to streamline the statements of our 
main results. 

We will refer to a value $\lambda \in I$ as an eigenvalue of (\ref{hammy})
if there exists a function 
$y(\cdot; \lambda) \in \AC([0,1], \mathbb{C}^{2n}) \backslash \{0\}$
that solves (\ref{hammy}) along with prescribed boundary conditions of 
the form {\bf (BC1)} or {\bf (BC2)}. (Here, $AC (\cdot)$ denotes absolute 
continuity.) Given any pair of values $\lambda_1, \lambda_2 \in I$, 
$\lambda_1 < \lambda_2$, our goal is to obtain a count $\mathcal{N} ([\lambda_1, \lambda_2))$,
including multiplicity, 
of the number of eigenvalues that (\ref{hammy}) has on the interval 
$[\lambda_1, \lambda_2)$. (The choice of a closed endpoint on the left
and an open endpoint on the right is taken by convention associated 
with the general definition that the Morse index of an operator corresponds
with a 
count of the number of eigenvalues the operator has strictly 
below $0$.) Our tool for this analysis will be the Maslov index, and  
as a starting point for a discussion of this object, we define what 
we will mean by a {\it Lagrangian subspace} of $\mathbb{C}^{2n}$.

\begin{definition} \label{lagrangian_subspace}
We say $\ell \subset \mathbb{C}^{2n}$ is a Lagrangian subspace of $\mathbb{C}^{2n}$
if $\ell$ has dimension $n$ and
\begin{equation} 
(J u, v) = 0, 
\end{equation} 
for all $u, v \in \ell$. (Here, $(\cdot, \cdot)$ denotes
the standard inner product on $\mathbb{C}^{2n}$.) In addition, we denote by 
$\Lambda (n)$ the collection of all Lagrangian subspaces of $\mathbb{C}^{2n}$, 
and we will refer to this as the {\it Lagrangian Grassmannian}. 
\end{definition}

\begin{remark} Following the convention of Arnol'd's foundational 
paper \cite{Arnold67}, the notation $\Lambda (n)$ is often used to denote the 
Lagrangian Grassmanian associated with $\mathbb{R}^{2n}$. Our expectation 
is that it can be used in the current setting of $\mathbb{C}^{2n}$ without
confusion. We note that the Lagrangian Grassmannian associated with 
$\mathbb{C}^{2n}$ has been considered by a number of authors, including
(ordered by publication date)
Bott \cite{Bott56}, Kostrykin and Schrader \cite{KS99}, Arnol'd
\cite{Arnold00}, and Schulz-Baldes \cite{S-B07, S-B12}. It is shown
in all of these references that $\Lambda (n)$ is 
homeomorphic to the set of $n \times n$ unitary matrices $U (n)$,
and in \cite{S-B07, S-B12} the relationship is shown to be 
diffeomorphic. It is also shown in \cite{S-B07} that the fundamental
group of $\Lambda (n)$ is the integers $\mathbb{Z}$.
 \end{remark}

Any Lagrangian subspace of $\mathbb{C}^{2n}$ can be
spanned by a choice of $n$ linearly independent vectors in 
$\mathbb{C}^{2n}$. We will generally find it convenient to collect
these $n$ vectors as the columns of a $2n \times n$ matrix $\mathbf{X}$, 
which we will refer to as a {\it frame} for $\ell$. Moreover, we will 
often coordinatize our frames as $\mathbf{X} = {X \choose Y}$, where $X$ and $Y$ are 
$n \times n$ matrices. Following \cite{F} (p. 274), we specify 
a metric on $\Lambda (n)$ in terms of appropriate orthogonal projections. 
Precisely, let $\mathcal{P}_i$ 
denote the orthogonal projection matrix onto $\ell_i \in \Lambda (n)$
for $i = 1,2$. I.e., if $\mathbf{X}_i$ denotes a frame for $\ell_i$,
then $\mathcal{P}_i = \mathbf{X}_i (\mathbf{X}_i^* \mathbf{X}_i)^{-1} \mathbf{X}_i^*$.
We take our metric $d$ on $\Lambda (n)$ to be defined 
by 
\begin{equation*}
d (\ell_1, \ell_2) := \|\mathcal{P}_1 - \mathcal{P}_2 \|,
\end{equation*} 
where $\| \cdot \|$ can denote any matrix norm. We will say 
that a path of Lagrangian subspaces 
$\ell: \mathcal{I} \to \Lambda (n)$ is continuous provided it is 
continuous under the metric $d$. 

Suppose $\ell_1 (\cdot), \ell_2 (\cdot)$ denote continuous paths of Lagrangian 
subspaces $\ell_i: \mathcal{I} \to \Lambda (n)$, for some parameter interval 
$\mathcal{I}$. The Maslov index associated with these paths, which we will 
denote $\mas (\ell_1, \ell_2; \mathcal{I})$, is a count of the number of times
the paths $\ell_1 (\cdot)$ and $\ell_2 (\cdot)$ intersect, counted
with both multiplicity and direction. (In this setting, if we let 
$t_*$ denote the point of intersection (often referred to as a 
{\it crossing point}), then multiplicity corresponds with the dimension 
of the intersection $\ell_1 (t_*) \cap \ell_2 (t_*)$; a precise definition of what we 
mean in this context by {\it direction} will be
given in Section \ref{maslov_section}.) 

The key ingredient we will need for connecting Maslov index 
calculations with renormalized oscillation results is 
monotonicity. We say that the evolution of $\mathcal{L} = (\ell_1, \ell_2)$
is {\it monotonic} provided all intersections occur with 
the same direction. If the intersections all correspond with 
the positive direction, and if the crossing points are all 
discrete, then we can compute
\begin{equation*}
\mas (\ell_1, \ell_2; \mathcal{I}) 
= \sum_{t \in \mathcal{I}} \dim (\ell_1 (t) \cap \ell_2 (t)).
\end{equation*}
Suppose $\mathbf{X}_1 (t) = {X_1 (t) \choose Y_1 (t)}$ and 
$\mathbf{X}_2 (t) = {X_2 (t) \choose Y_2 (t)}$ respectively 
denote frames for Lagrangian subspaces of $\mathbb{C}^{2n}$,
$\ell_1 (t)$ and $\ell_2 (t)$. Then we can express
this last relation as 
\begin{equation*}
\mas (\ell_1, \ell_2; \mathcal{I}) 
= \sum_{t \in \mathcal{I}} \dim \ker (\mathbf{X}_1 (t)^* J \mathbf{X}_2 (t)).
\end{equation*}
(See Lemma \ref{intersection_lemma} below.)
The right-hand side of this final expression, expressed in terms
of the matrix Wronskian $W(t) = \mathbf{X}_1 (t)^* J \mathbf{X}_2 (t)$, 
has the form we associate with renormalized oscillation theory (see, e.g., 
\cite{GZ2017}), and we will sometimes adopt the notation 
of \cite{GZ2017} and use the {\it counting function} 
\begin{equation} \label{counting-function}
N_{\mathcal{I}} (\mathbf{X}_1 (\cdot)^* J \mathbf{X}_2 (\cdot))
:= \sum_{t \in \mathcal{I}} \dim \ker (\mathbf{X}_1 (t)^* J \mathbf{X}_2 (t)).
\end{equation}
 
\begin{remark} \label{gz-remark}
Renormalized oscillation theory was introduced 
in \cite{GST1996} in the context of single Sturm-Liouville
equations, and subsequently was developed in 
\cite{Teschl1996, Teschl1998} for Jacobi operators and 
Dirac operators. More recently, Gesztesy and Zinchenko
have extended these early results to the setting of 
(\ref{hammy}) with 
$\mathbb{B} (x; \lambda) = \lambda A(x) + B(x)$
(with suitable assumptions on $A$ and $B$) and for three
classes of domain: bounded, half-line, and $\mathbb{R}$
(see \cite{GZ2017}). This last reference served as the 
direct motivation for our analysis.
\end{remark}

In order to formulate our theorem regarding {\bf (BC1)}, we fix 
any pair $\lambda_1, \lambda_2 \in I$, with $\lambda_1 < \lambda_2$,
and let $\mathbf{X}_1 (x; \lambda)$ denote a $2n \times n$ matrix 
solving 
\begin{equation} \label{frame1}
\begin{aligned}
J \mathbf{X}_1' &= \mathbb{B} (x; \lambda) \mathbf{X}_1 \\
\mathbf{X}_1 (0; \lambda) &= J \alpha^*. 
\end{aligned}
\end{equation}
Under our assumptions {\bf (A)} on $\mathbb{B} (x; \lambda)$, 
we can conclude that for each $\lambda \in I$, 
$\mathbf{X}_1 (\cdot; \lambda) \in AC ([0,1], \mathbb{C}^{2n \times n})$,
and additionally that $\mathbf{X}_1 (x; \lambda)$ is differentiable 
in $\lambda$ with 
$\partial_{\lambda} \mathbf{X}_1 (\cdot; \lambda) \in AC ([0,1], \mathbb{C}^{2n \times n})$
and 
\begin{equation}\label{lambda-derivative-relation}
    J (\partial_{\lambda} \mathbf{X}_1 (x; \lambda))' 
    = \mathbb{B}_{\lambda} (x; \lambda) \mathbf{X}_1 (\cdot; \lambda)
    + \mathbb{B} (x; \lambda) \partial_{\lambda} \mathbf{X}_1 (x; \lambda),
\end{equation}
for a.e. $x \in (0, 1)$. 

\begin{remark} \label{lambda-derivative-remark} 
Regarding (\ref{lambda-derivative-relation}), the observation is simply that 
we can justify switching the order of differentiation of $\mathbf{X}_1 (x; \lambda)$
with respect to $x$ and $\lambda$. Our assumptions allow us to do this by 
integrating (\ref{frame1}) to 
\begin{equation*}
    J \mathbf{X}_1 (x; \lambda) = - \alpha^*
    + \int_0^x \mathbb{B} (\xi; \lambda) \mathbf{X}_1 (\xi; \lambda) d\xi,
\end{equation*}
and then differentiating through the integral in $\lambda$, followed by 
differentiation in $x$. These are straightforward calculations that follow 
the approach of Chapter 2 in \cite{Weidmann1987}.
\end{remark}

As shown in \cite{HJK2018}, for each pair $(x, \lambda) \in [0,1] \times I$,
$\mathbf{X}_1 (x; \lambda)$ is the frame for a Lagrangian subspace
$\ell_1 (x; \lambda)$. (In \cite{HJK2018}, the authors make slightly 
stronger assumptions on $\mathbb{B} (x; \lambda)$, but their proof
carries over immediately into our setting.)
Likewise, keeping in mind that $\lambda_2$ is fixed, we let 
$\ell_2 (x; \lambda_2)$ denote the map of Lagrangian subspaces
associated with frames $\mathbf{X}_2 (x; \lambda_2)$ solving 
\begin{equation} \label{frame2}
\begin{aligned}
J \mathbf{X}_2' &= \mathbb{B} (x; \lambda_2) \mathbf{X}_2 \\
\mathbf{X}_2 (1; \lambda_2) &= J \beta^*. 
\end{aligned}
\end{equation}
We emphasize that $\mathbf{X}_2 (x; \lambda_2)$ is initialized 
at $x=1$. 

In addition to Assumptions {\bf (A)}, we make the following 
{\it positivity} assumptions: 

\medskip
{\bf (B1)} For any $\lambda \in I$, the matrix
\begin{equation} \label{positive-definite}
\int_0^1 \mathbf{X}_1 (x; \lambda)^* \mathbb{B}_{\lambda} (x; \lambda) \mathbf{X}_1 (x; \lambda) dx 
\end{equation}
is positive definite. 

\medskip
{\bf (B2)} For the fixed values $\lambda_1, \lambda_2 \in I$, 
$\lambda_1 < \lambda_2$, the matrix $(\mathbb{B} (x; \lambda_2) - \mathbb{B} (x; \lambda_1))$ 
is non-negative for a.e. $x \in (0,1)$, and moreover 
there is no interval $[a, b] \subset [0, 1]$, $a < b$, so that 
\begin{equation*}
\dim (\ell_1 (x; \lambda_1) \cap \ell_2 (x; \lambda_2)) \ne 0
\end{equation*}
for all $x \in [a, b]$.

\begin{remark} \label{positivity-remark}
Assumption {\bf (B1)} is standard for ensuring 
that as $\lambda$ varies, with $x$ fixed, crossings of 
$\ell_1 (x; \lambda)$ and $\ell_2 (x; \lambda_2)$ will 
all occur in the same direction. The frame $\mathbf{X} (x; \lambda_2)$
does not appear in the assumption, because it does not vary 
with $\lambda$. For Assumption {\bf (B2)},
the a.e. non-negativity of 
$(\mathbb{B} (x; \lambda_2) - \mathbb{B} (x; \lambda_1))$
ensures that as $x$ varies, crossings of $\ell_1 (x; \lambda_1)$
and $\ell_2 (x; \lambda_2)$ all occur in the same direction, while 
the moreover part ensures that $\ell_1 (x; \lambda_1)$
and $\ell_2 (x; \lambda_2)$ cannnot get stuck at an intersection.  

We will verify in Section \ref{proofs} that the moreover part of 
{\bf (B2)} is implied by the following form of Atkinson 
positivity: for any $[a, b] \subset [0, 1]$, $a < b$, and any non-trivial 
solution $y(\cdot; \lambda_1) \in \AC ([0, 1], \mathbb{C}^{2n})$ 
of $J y' = \mathbb{B} (x; \lambda_1) y$, 
we must have 
\begin{equation} \label{difference-definite}
\int_a^b ((\mathbb{B} (x; \lambda_2) - \mathbb{B} (x; \lambda_1)) y(x; \lambda_1), y(x; \lambda_1)) dx > 0. 
\end{equation}
In the case that $\mathbb{B} (x; \lambda)$ is linear in $\lambda$ 
(as in Remark \ref{gz-remark}), 
non-negativity of $(\mathbb{B} (x; \lambda_2) - \mathbb{B} (x; \lambda_1))$ corresponds
with non-negativity of the matrix $A$, and the integral conditions
(\ref{positive-definite}) and (\ref{difference-definite})
are both equivalent to Atkinson positivity (see, e.g., Section 4 in \cite{LS2012}
and Section IV.4 in \cite{Krall2002}). For a development of renormalized 
oscillation theory under fewer such restrictions, we refer the reader 
to \cite{Elyseeva2021} and the references therein. 
\end{remark}

Suppose that for some value $\lambda \in I$ equation (\ref{hammy})
with specified boundary conditions admits one or more linearly independent
solutions. We denote the subspace spanned by these solutions by 
$\mathbb{E} (\lambda)$, noting that $\dim \mathbb{E} (\lambda) \le 2n$.
Given any two values $\lambda_1, \lambda_2 \in I$, with $\lambda_1 < \lambda_2$,
it is shown in \cite{HJK2018} that under positivity assumptions {\bf (B1)} 
the {\it spectral count} 
\begin{equation} \label{spectral_count}
\mathcal{N} ([\lambda_1, \lambda_2)) 
:= \sum_{\lambda \in [\lambda_1, \lambda_2)} \dim \mathbb{E} (\lambda),
\end{equation} 
is well-defined. It's clear that $\mathcal{N} ([\lambda_1, \lambda_2))$
is a count of the eigenvalues of (\ref{hammy}) on $[\lambda_1, \lambda_2)$,
counted with geometric multiplicity. In order to understand the nature 
of algebraic multiplicity in this setting, as well as the notion 
of essential spectrum, it's useful to frame our discussion in terms of 
the operator pencil 
\begin{equation*}
\mathcal{L} (\lambda) = J \frac{d}{dx} - \mathbb{B} (x; \lambda),
\end{equation*}
specified on the domain (independent of $\lambda$)  
\begin{equation*}
\begin{aligned}
\mathcal{D} &:= \{y \in L^2 ((0,1), \mathbb{C}^{2n}): 
y \in \AC ([0,1], \mathbb{C}^{2n}), \\
&\quad \quad \mathcal{L} y \in L^2 ((0,1), \mathbb{C}^{2n}),
\, \alpha y(0) = 0, \, \beta y(1) = 0 \}
\end{aligned}
\end{equation*}
(for boundary conditions {\bf (BC1)}, and with a similar specification 
for boundary conditions {\bf (BC2)}).
Using the methods of \cite{Weidmann1987}, we can readily verify that for each $\lambda \in I$, 
$\mathcal{L} (\lambda)$ is Fredholm and self-adjoint on $\mathcal{D}$, from which we can 
conclude that $\mathcal{L}$ has no essential spectrum on $I$. 
Moreover, under slightly
stronger assumptions on $\mathbb{B}$ 
(in particular, $\mathbb{B}(\cdot; \lambda) \in L^{2} ((0,1),\mathbb{C}^{2n \times 2n})$ 
for all $\lambda \in I$),
we can verify that $\mathcal{L}$ has no Jordan chains of length greater than one, 
implying that the algebraic and geometric multiplicities of its eigenvalues agree.
(See the appendix for further discussion, and also Section 1.2 
of \cite{HP2017}, in which the authors consider the same operator pencil
under slightly stronger assumptions on $\mathbb{B} (x; \lambda)$.)  

We will establish the following theorem.

\begin{theorem} \label{bc1_theorem} 
Fix $\lambda_1, \lambda_2 \in I$, $\lambda_1 < \lambda_2$.
For equation (\ref{hammy}), let Assumptions {\bf (A)} and {\bf (B1)} hold, 
and let $\mathbf{X}_1$ and $\mathbf{X}_2$ respectively denote
the Lagrangian frames specified in (\ref{frame1}) and 
(\ref{frame2}). If $\mathcal{N} ([\lambda_1, \lambda_2))$ denotes the spectral
count for (\ref{hammy}) with boundary conditions {\bf (BC1)},
then 
\begin{equation} \label{thm1-eqn1}
\mathcal{N} ([\lambda_1, \lambda_2))
= \mas (\ell_1 (\cdot; \lambda_1), \ell_2 (\cdot; \lambda_2); [0,1]).
\end{equation}
Moreover, under the additional assumption {\bf (B2)}, with 
$\lambda_1$ and $\lambda_2$ as in {\bf (B2)}, we have 
\begin{equation} \label{thm1-eqn2}
\mathcal{N} ([\lambda_1, \lambda_2))
= \sum_{x \in (0,1]} \dim \ker (\mathbf{X}_1 (x; \lambda_1)^* J \mathbf{X}_2 (x; \lambda_2)).
\end{equation}
\end{theorem}

We note that the final relation in Theorem \ref{bc1_theorem} could be stated with 
$[\lambda_1, \lambda_2)$ and $(0,1]$ replaced respectively 
with any of the following three combinations: 
$(\lambda_1, \lambda_2)$ and $(0,1)$, $[\lambda_1, \lambda_2]$ and $[0,1]$,
or $(\lambda_1, \lambda_2]$ and $[0,1)$. For the correspondence
between $\lambda = \lambda_2$ and $x = 0$ (i.e., the fact that 
the associated delimiters are either both round or both square),
we observe that $\lambda_2$ is an eigenvalue of (\ref{hammy})
if and only if $\ell_1 (0; \lambda_2)$ and $\ell_2 (0; \lambda_2)$
intersect, and since $\ell_1 (0; \lambda_1) = \ell_1 (0; \lambda_2)$,
we can conclude that $\lambda_2$ is an eigenvalue of (\ref{hammy})
if and only if $\ell_1 (0; \lambda_1)$ and $\ell_2 (0; \lambda_2)$
intersect. This final intersection is precisely the quantity 
detected at $x = 0$ in the sum. Likewise, for the correspondence
between $\lambda = \lambda_1$ and $x = 1$, $\lambda_1$ is an eigenvalue
of (\ref{hammy}) if and only if $\ell_1 (1; \lambda_1)$ and 
$\ell_2 (1; \lambda_1)$ intersect, and 
$\ell_2 (1; \lambda_1) = \ell_2 (1; \lambda_2)$. In addition, we observe
that by exchanging the pair $(\ell_1 (x; \lambda_1), \ell_2 (x; \lambda_2))$
in our analysis with the pair $(\ell_1 (x; \lambda_2), \ell_2 (x; \lambda_1))$,
we arrive at the alternative formulation 
\begin{equation*} 
\mathcal{N} ([\lambda_1, \lambda_2))
= - \mas (\ell_1 (\cdot; \lambda_2), \ell_2 (\cdot; \lambda_1); [0,1]),
\end{equation*}
which, under Assumption {\bf (B2)}, implies 
\begin{equation} \label{thm1-eqn3} 
\mathcal{N} ([\lambda_1, \lambda_2))
= \sum_{x \in [0,1)} \dim \ker (\mathbf{X}_1 (x; \lambda_2)^* J \mathbf{X}_2 (x; \lambda_1)).
\end{equation}
In particular, as opposed to (\ref{thm1-eqn2}), the sum on the right-hand side of (\ref{thm1-eqn3}) 
includes $x = 0$ and excludes $x = 1$. Moreover, in this case,
the pairing of $[\lambda_1, \lambda_2)$ with $[0, 1)$ can be replaced by 
any other pairing in which the delimiters of $\lambda_1$ and $0$ 
agree, and the delimiters of $\lambda_2$ and $1$ agree. 
 
Our Theorem \ref{bc1_theorem} is quite similar to Theorem 3.10
of \cite{GZ2017}, though we observe the following differences: in \cite{GZ2017}, the authors
use Theorem 3.10 as a statement about both the case of bounded intervals
and the case of half-lines, and we are comparing with the bounded-interval
statement. Keeping this in mind, our theorem is slightly less
restrictive, in that (1) it allows for a more general class
of matrices $\mathbb{B} (x; \lambda)$; and (2) it allows for the 
possibility that $\lambda_1$ and/or $\lambda_2$ is an eigenvalue  
of (\ref{hammy}). Regarding Item (1), our Theorem \ref{bc1_theorem}
allows for $\mathbb{B} (x; \lambda)$ to depend nonlinearly on $\lambda$.
This happens, for example, when (\ref{hammy}) arises from consideration 
of certain differential-algebraic systems; in Section \ref{applications}
we give one such example. Having drawn this comparison, we should emphasize 
that the primary goal of \cite{GZ2017} (and indeed the primary motivation 
for the original work of \cite{GST1996}) was to develop a theory that 
would allow the authors to count discrete eigenvalues between bands of 
essential spectrum. We do not consider that important case here. 

\begin{remark} \label{standard-v-renormalized} 
It's instructive to view Theorem \ref{bc1_theorem} in 
relation to the results of \cite{HJK2018, HS2016}, for which the authors 
specify $\mathbf{X}_1 (x; \lambda)$ precisely as here, but in lieu
of $\mathbf{X}_2 (x; \lambda_2)$, use the fixed frame 
$\tilde{\mathbf{X}}_2 = J \beta^*$ for a target space $\tilde{\ell}_2$. 
Under Assumptions {\bf (A)}
and {\bf (B1)}, and assuming boundary conditions {\bf (BC1)}, 
the methods of \cite{HJK2018, HS2016} can be used to establish the 
relation 
\begin{equation} \label{old_news}
\mathcal{N} ([\lambda_1, \lambda_2))
= - \mas (\ell_1 (\cdot; \lambda_1), \tilde{\ell}_2; [0, 1])
+ \mas (\ell_1 (\cdot; \lambda_2), \tilde{\ell}_2; [0, 1]).
\end{equation}
Critically, however, in the 
setting of \cite{HJK2018, HS2016}, intersections between 
$\ell_1 (\cdot; \lambda_i)$ ($i = 1,2$) and $\tilde{\ell}_2$
are not necessarily monotone, and so (\ref{old_news})
cannot generally be formulated as a simple count of nullities
as in (\ref{thm1-eqn2}) of Theorem \ref{bc1_theorem}. 
One way in which this difference in approaches manifests is 
in the nature of {\it spectral curves}, by which we mean 
continuous paths of crossing-point pairs $(x, \lambda)$ (assuming
for simplicity of the discussion that these curves are 
non-intersecting). In the setting of \cite{HJK2018, HS2016}, such 
curves can reverse direction as depicted on the left-hand 
sides of Figures \ref{dirac_figure}, \ref{sl_figure}, 
and \ref{fourth-figure} in Section \ref{applications}, 
but in the current setting these spectral curves 
are necessarily monotone, as depicted on the right-hand
sides of the same figures. This 
dynamic can be viewed as a graphical interpretation of why 
renormalized oscillation theory works in the elegant way 
that it does.  
\end{remark}

For the case of (\ref{hammy}) with boundary conditions {\bf (BC2)},
we follow \cite{HJK2018} and begin by defining a Lagrangian subspace in terms
of the ``trace" operator
\begin{equation} \label{trace_operator}
\mathcal{T}_x y := \mathcal{M} {y(0) \choose y(x)},
\end{equation}
where 
\begin{equation*}
\mathcal{M} = 
\begin{pmatrix}
I_n & 0 & 0 & 0 \\
0 & 0 & I_n & 0 \\
0 & -I_n & 0 & 0 \\
0 & 0 & 0 & I_n
\end{pmatrix}.
\end{equation*} 
In \cite{HJK2018}, the authors verify that the subspace 
\begin{equation} \label{lagrangian2}
\ell_3 (x; \lambda) := \{\mathcal{T}_x y: y (\cdot; \lambda) \in \AC ([0, 1], \mathbb{C}^{2n}), \,
J y' = \mathbb{B} (x; \lambda) y 
\text{ a.e } x \in (0,1) \}
\end{equation}
is a Lagrangian subspace of $\mathbb{C}^{2n}$ for all $(x,\lambda) \in [0,1] \times I$. 

In order to establish notation for the statement of our second theorem, we  
let $\Phi (x; \lambda)$ denote the $2n \times 2n$
fundamental matrix solution to 
\begin{equation} \label{cap_phi}
J \Phi' = \mathbb{B} (x; \lambda) \Phi; 
\quad \Phi (0; \lambda) = I_{2n},
\end{equation}
and write 
\begin{equation*}
\Phi (x; \lambda) = 
\begin{pmatrix}
\Phi_{11} (x; \lambda) & \Phi_{12} (x; \lambda) \\
\Phi_{21} (x; \lambda) & \Phi_{22} (x; \lambda) 
\end{pmatrix}.
\end{equation*}
With this notation, we can express the frame for $\ell_3 (x; \lambda)$
as 
\begin{equation} \label{frame3}
\mathbf{X}_3 (x, \lambda) = {X_3 (x, \lambda) \choose Y_3 (x, \lambda)}
=
\begin{pmatrix}
I_n & 0 \\
\Phi_{11} (x; \lambda) & \Phi_{12} (x; \lambda) \\
0 & -I_n \\
\Phi_{21} (x; \lambda) & \Phi_{22} (x; \lambda)
\end{pmatrix}.
\end{equation}

We see by direct calculation that $\mathbf{X}_3 (x; \lambda)$ 
can be interpreted as the frame associated with a linear
Hamiltonian system 
\begin{equation*}
J_{4n} \mathbf{X}_3' = \mathcal{B} (x; \lambda) \mathbf{X}_3,
\end{equation*}
where 
\begin{equation} \label{mathcal-B-defined}
\mathcal{B} (x; \lambda) 
= \begin{pmatrix}
0 & 0 & 0 & 0 \\
0 & B_{11} (x; \lambda) & 0 & B_{12} (x; \lambda) \\
0 & 0 & 0 & 0 \\
0 & B_{21} (x; \lambda) & 0 & B_{22} (x; \lambda)
\end{pmatrix}; \quad 
\text{using} \, \,
\mathbb{B} = 
\begin{pmatrix}
B_{11} & B_{12} \\
B_{21} & B_{22}
\end{pmatrix},
\end{equation}
and the flow is initialized by 
\begin{equation*}
\mathbf{X}_3 (0; \lambda) 
=
\begin{pmatrix}
I_n & 0 \\
I_n & 0 \\
0 & -I_n \\
0 & I_n
\end{pmatrix}.
\end{equation*}
Here, we use the notation $J_{4n}$ to designate the matrix $J$ with each $I_n$
replaced by $I_{2n}$.

For our second path of Lagrangian subspaces, we let $\mathbf{X}_4 (x; \lambda_2)$
solve
\begin{equation} \label{frame4}
\begin{aligned}
J_{4n} \mathbf{X}_4' &= \mathcal{B} (x; \lambda_2) \mathbf{X}_4 \\
\mathbf{X}_4 (1; \lambda_2) &= \mathcal{M} \mathcal{J}_{4n} \Theta^*.
\end{aligned}
\end{equation}
In \cite{HJK2018}, the authors verify that $\mathcal{M} \mathcal{J}_{4n} \Theta^*$
is the frame for a Lagrangian subspace of $\mathbb{C}^{4n}$, and that intersections of 
$\ell_3 (1; \lambda)$ with this Lagrangian subspace correspond with eigenvalues 
of (\ref{hammy})-{\bf (BC2)}. (Again, strictly speaking, the 
authors of \cite{HJK2018} are working with Lagrangian subspaces of 
$\mathbb{R}^{4n}$.)

In this case, we make the following positivity 
assumptions: 

\medskip
{\bf (B1)$\mathbf{'}$} For any $\lambda \in I$, the matrix
\begin{equation} \label{positive_definite2}
\int_0^1 \Phi (x; \lambda)^* \mathbb{B}_{\lambda} (x; \lambda) \Phi (x; \lambda) dx 
\end{equation}
is positive definite. 

\medskip
{\bf (B2)$\mathbf{'}$} For the fixed values $\lambda_1, \lambda_2 \in I$, 
$\lambda_1 < \lambda_2$, the matrix $(\mathbb{B} (x; \lambda_2) - \mathbb{B} (x; \lambda_1))$ 
is non-negative for a.e. $x \in (0,1)$, and moreover 
there is no interval $[a, b] \subset [0, 1]$, $a < b$, so that 
\begin{equation*}
\dim (\ell_3 (x; \lambda_1) \cap \ell_4 (x; \lambda_2)) \ne 0
\end{equation*}
for all $x \in [a, b]$.

\medskip

We are now in a position to state our second theorem. 

\begin{theorem} \label{bc2_theorem}
Fix $\lambda_1, \lambda_2 \in I$, $\lambda_1 < \lambda_2$.
For equation (\ref{hammy}), let Assumptions {\bf (A)}
and {\bf (B1)$\mathbf{'}$} hold,  
and let $\mathbf{X}_3$ and $\mathbf{X}_4$ respectively denote
the Lagrangian frames specified in (\ref{frame3}) and
(\ref{frame4}). If $\mathcal{N} ([\lambda_1, \lambda_2))$ 
denotes the spectral count for (\ref{hammy}) with boundary 
conditions {\bf (BC2)}, then
\begin{equation} \label{thm2-eqn1}
\mathcal{N} ([\lambda_1, \lambda_2))
= \mas (\ell_3 (\cdot; \lambda_1), \ell_4 (\cdot; \lambda_2); [0,1]).
\end{equation}
Moreover, under the additional assumption
{\bf (B2)$\mathbf{'}$}, with $\lambda_1$ and $\lambda_2$ as in 
{\bf (B2)$\mathbf{'}$}, we have,   
\begin{equation} \label{thm2-eqn2}
\mathcal{N} ([\lambda_1, \lambda_2))
= \sum_{x \in (0,1]} \dim \ker (\mathbf{X}_3 (x; \lambda_1)^* J \mathbf{X}_4 (x; \lambda_2)).
\end{equation}
\end{theorem}

\begin{remark} 
The considerations discussed in Remark \ref{positivity-remark}
carry over to the setting of Theorem \ref{bc2_theorem}, and 
in particular, we will verify in Section \ref{proofs} that the 
moreover part of {\bf (B2)$\mathbf{'}$} is implied by 
(\ref{difference-definite}) precisely as previously stated
(i.e., it's not necessary to replace $\mathbb{B}$ with 
$\mathcal{B}$). 

Although the proof of Theorem \ref{bc2_theorem} is 
essentially identical to the proof of Theorem \ref{bc1_theorem}
(at least in our development), we are not aware of a statement 
along the lines of Theorem \ref{bc2_theorem} in the literature. 
\end{remark}

The remainder of the paper is organized as follows. In 
Section \ref{maslov_section}, we develop the Maslov index 
framework in $\mathbb{C}^{2n}$, and in Section \ref{rotation_section}
we develop the tools we will need to verify the monotonicity 
that will be necessary to conclude the second statements 
in Theorems \ref{bc1_theorem}
and \ref{bc2_theorem}. In Section \ref{proofs}, we prove 
Theorems \ref{bc1_theorem} and \ref{bc2_theorem}, and in 
Section \ref{applications} we conclude by verifying that our 
assumptions hold for five example cases: Dirac systems, 
Sturm-Liouville systems, the class of systems analyzed in 
\cite{GZ2017}, a system associated with 
differential-algebraic Sturm-Liouville systems (for which 
$\mathbb{B} (x; \lambda)$ depends nonlinearly on $\lambda$), 
and a self-adjoint fourth-order equation for which we take
periodic boundary conditions.
In a short appendix, we discuss the interpretation of 
$\mathcal{L} (\lambda)$ as an operator pencil.

\section{The Maslov Index on $\mathbb{C}^{2n}$} \label{maslov_section}

In this section, we verify that the framework developed in \cite{HLS2017} for 
computing the Maslov index for Lagrangian pairs in $\mathbb{R}^{2n}$ extends 
to the case of Lagrangian pairs in $\mathbb{C}^{2n}$. A similar framework 
has been developed in \cite{S-B12}, and in particular, some of the 
results in this section correspond with results in Section 2 of that
reference. We include details here (1) for completeness; and 
(2) because we need some additional information that is not 
developed in \cite{S-B12}. 

As a starting point, we note the following direct relation between 
Lagrangian subspaces on $\mathbb{C}^{2n}$ and their associated frames.

\begin{proposition} \label{lagrangian_property} 
A $2n \times n$ matrix $\mathbf{X} = {X \choose Y}$
is a frame for a Lagrangian subspace of $\mathbb{C}^{2n}$ if and 
only if the columns of $\mathbf{X}$ are linearly independent, and additionally 
\begin{equation*}
X^* Y - Y^* X = 0.
\end{equation*} 
We refer to this relation as the Lagrangian property for frames.
\end{proposition}

\begin{remark} The straightforward proof of Proposition \ref{lagrangian_property}
is essentially the same as for the case of $\mathbb{R}^{2n}$ (see Proposition 2.1
in \cite{HLS2017}). We note that the Lagrangian property can also be
expressed as $\mathbf{X}^* J \mathbf{X} = 0$. According to the Fredholm 
Alternative, $\mathbb{C}^{2n} = \operatorname{ran} (\mathbf{X}) 
\oplus \operatorname{ker} (\mathbf{X}^*)$, and since 
$\dim \operatorname{ran} \mathbf{X} = n$,
we must have $\dim \ker \mathbf{X}^* = n$. I.e., Lagrangian subspaces 
on $\mathbb{C}^{2n}$ are {\it maximal}; 
no subspace of $\mathbb{C}^{2n}$ with dimension greater than $n$ can 
have the Lagrangian property.   
\end{remark}

\begin{proposition} If $\mathbf{X} = {X \choose Y}$ is the frame for a Lagrangian
subspace of $\mathbb{C}^{2n}$,
then the matrices $X \pm i Y$ are both invertible. 
\end{proposition}

\begin{proof} First, it's standard that $\mathbf{X}^* \mathbf{X} = X^* X + Y^* Y$
is invertible if and only if the columns of $\mathbf{X}$ are linearly 
independent. Now suppose $v \in \ker (X+iY)$ so that $(X+iY)v=0$. We can 
multiply this equation by $(X^*-iY^*)$ and use the Lagrangian
property to see that $(X^* X + Y^* Y) v = 0$. This implies
$v = 0$, from which we can conclude that $(X+iY)$ is invertible.
The case $(X-iY)$ is similar.
\end{proof} 

\begin{proposition} \label{useful1}
If $\mathbf{X} = {X \choose Y}$ is the frame for a Lagrangian subspace of $\mathbb{C}^{2n}$,
then
\begin{equation*}
(X \pm i Y)^{-1} = M^2 (X^* \mp iY^*),
\end{equation*}
where $M := (\mathbf{X}^* \mathbf{X})^{-1/2}$.
\end{proposition}

\begin{proof} Computing directly, we see that 
\begin{equation*}
\begin{aligned}
(X^* - i Y^*) (X+iY) &= X^* X + Y^* Y + i (X^*Y - Y^*X) \\
&= (M^2)^{-1},
\end{aligned}
\end{equation*}
where we have used the Lagrangian property for frames to see that the imaginary
part is $0$. The claim now follows upon multiplication
on the left by $M^2$ and on the right by $(X+iY)^{-1}$. The case 
$(X - iY)^{-1}$ is similar.
\end{proof}

For the next proposition, we set 
\begin{equation*}
\tilde{W}_D := (X + iY)(X-iY)^{-1},
\end{equation*}
noting that the subscript $D$ indicates (as we will check just below) 
that $\tilde{W}_D$ detects intersections of $\ell = \colspan(\mathbf{X})$ 
with the Dirichlet plane $\ell_D = \colspan({0 \choose I_n})$.

\begin{proposition} \label{unitary} 
If $\mathbf{X} = {X \choose Y}$ is the frame for a Lagrangian subspace of 
$\mathbb{C}^{2n}$, then $\tilde{W}_D$ is unitary.
\end{proposition}

\begin{proof} Computing directly, we find 
\begin{equation*}
\begin{aligned}
\tilde{W}_D^* \tilde{W}_D &= 
(X^* + iY^*)^{-1} (X^* - iY^*) (X + iY)(X-iY)^{-1} \\
&= (X^* + iY^*)^{-1} \Big{\{} X^*X + Y^* Y \Big{\}} (X-iY)^{-1} \\
&= (X^* + iY^*)^{-1} (M^2)^{-1} M^2 (X^* + iY^*)  = I, 
\end{aligned}
\end{equation*}
where for the second equality we used the Lagrangian property for frames, 
and for the third we used Proposition \ref{useful1}.
\end{proof}

\begin{definition} Let $\mathbf{X}_1 = {X_1 \choose Y_1}$ and 
$\mathbf{X}_2 = {X_2 \choose Y_2}$ denote frames for two 
Lagrangian subspaces of $\mathbb{C}^{2n}$. We define 
\begin{equation} \label{tildeW}
\tilde{W} := - (X_1 + iY_1)(X_1-iY_1)^{-1} (X_2 - iY_2)(X_2+iY_2)^{-1}. 
\end{equation}
According to Proposition \ref{unitary}, $\tilde{W}$ is the product of 
unitary matrices, and so is unitary. Consequently, the eigenvalues of 
$\tilde{W}$ will be confined to the unit circle in the complex plane, 
$S^1$.
\end{definition}

\begin{lemma} \label{dimensions_lemma}
Suppose $\mathbf{X}_1 = {X_1 \choose Y_1}$ and 
$\mathbf{X}_2 = {X_2 \choose Y_2}$ respectively denote frames for 
Lagrangian subspaces of $\mathbb{C}^{2n}$. Then
\begin{equation*}
\dim \ker (\mathbf{X}_1^* J \mathbf{X}_2) = \dim \ker (\tilde{W} + I).
\end{equation*}
More precisely, 
\begin{equation*}
\ker (\mathbf{X}_1^* J \mathbf{X}_2) 
= \ran \Big((X_2 + i Y_2)^{-1}\Big|_{\ker (\tilde{W} + I)}\Big).
\end{equation*}
\end{lemma}

\begin{proof} First, suppose 
\begin{equation*}
\dim \ker \mathbf{X}_1^* J \mathbf{X}_2 = m > 0,
\end{equation*}
and let $\{v_k\}_{k=1}^m$ denote a basis for this kernel. Then, 
in particular, 
\begin{equation} \label{CLPk}
(-X_1^* Y_2 + Y_1^* X_2) v_k = 0
\end{equation}
for all $k \in \{1, 2, \dots, m\}$. Set 
\begin{equation*}
w_k = (X_2 + i Y_2) v_k,
\end{equation*}
and notice that since $X_2 + i Y_2$ is invertible, $\{v_k\}_{k=1}^m$
comprises a linearly independent set of vectors if and only if 
$\{w_k\}_{k=1}^m$
comprises a linearly independent set of vectors.

Now, we compute directly, 
\begin{equation*}
\begin{aligned}
\tilde{W} w_k &= - (X_1 + iY_1)(X_1-iY_1)^{-1} (X_2 - iY_2) v_k \\
&= - (X_1 + iY_1) M_1^2 (X_1^* + i Y_1^*) (X_2 - iY_2) v_k \\
&= - (X_1 + iY_1) M_1^2 \Big{\{} X_1^* X_2 + Y_1^* Y_2 
+ i (Y_1^* X_2 - X_1^* Y_2) \Big{\}} v_k.
\end{aligned}
\end{equation*}
Using (\ref{CLPk}), we see that 
\begin{equation*}
\begin{aligned}
\tilde{W} w_k &= - (X_1 + iY_1) M_1^2 (X_1^* - i Y_1^*) (X_2 + iY_2) v_k \\
&= - (X_1 + iY_1) (X_1 + iY_1)^{-1} w_k = - w_k,
\end{aligned}
\end{equation*}
and so $w_k \in \ker (\tilde{W} + I)$. We can conclude that 
\begin{equation*}
\dim \ker \mathbf{X}_1^* J \mathbf{X}_2 \le \dim \ker (\tilde{W} + I).
\end{equation*}
For the second part, our calculation has established 
\begin{equation*}
\ker \mathbf{X}_1^* J \mathbf{X}_2 
\subset \ran \Big((X_2 + i Y_2)^{-1}\Big|_{\ker (\tilde{W} + I)}\Big).
\end{equation*}

Turning the argument around, we get the inequality and the associated
inclusion in the other direction, so we can conclude equality in both 
cases.
\end{proof}

\begin{lemma} \label{intersection_lemma}
Suppose $\mathbf{X}_1 = {X_1 \choose Y_1}$ and 
$\mathbf{X}_2 = {X_2 \choose Y_2}$ respectively denote frames for 
Lagrangian subspaces of $\mathbb{C}^{2n}$, $\ell_1$ and $\ell_2$. Then 
\begin{equation*}
\dim \ker (\mathbf{X}_1^* J \mathbf{X}_2) = \dim (\ell_1 \cap \ell_2).
\end{equation*}
More precisely, 
\begin{equation*}
\ran \Big(\mathbf{X}_2\Big|_{\ker (\mathbf{X}_1^* J \mathbf{X}_2)}\Big) = \ell_1 \cap \ell_2.
\end{equation*}
\end{lemma}

\begin{proof} First, suppose 
\begin{equation*}
\dim \ker (\mathbf{X}_1^* J \mathbf{X}_2) = m > 0,
\end{equation*}
and let $\{v_k\}_{k=1}^m$ denote a basis for this kernel. Then, 
in particular, 
\begin{equation} \label{CLPk2}
\mathbf{X}_1^* J \mathbf{X}_2 v_k = 0
\end{equation}
for all $k \in \{1, 2, \dots, m\}$.

Set $\zeta_k = \mathbf{X}_2 v_k$. Then $\zeta_k \in \ell_2$
(as a linear combination of the columns of $\mathbf{X}_2$),
and since $\mathbf{X}_1^* J \zeta_k = 0$, we must have 
$\zeta_k \in \ell_1$ by maximality. We see that 
we can associate with the $\{v_k\}_{k=1}^m$ a linearly 
independent set $\{\zeta_k\}_{k=1}^m \subset \ell_1 \cap \ell_2$,
and so 
\begin{equation*}
\dim \ker \mathbf{X}_1^* J \mathbf{X}_2 \le \dim (\ell_1 \cap \ell_2).
\end{equation*}
For the second part, our calculation has established 
\begin{equation*}
\ran \Big(\mathbf{X}_2\Big|_{\ker (\mathbf{X}_1^* J \mathbf{X}_2)}\Big) 
\subset \ell_1 \cap \ell_2.
\end{equation*}

Turning the argument around, we get the inequality and the associated
inclusion in the other direction, so we can conclude equality in both 
cases.
\end{proof}

Combining Lemmas \ref{dimensions_lemma} and \ref{intersection_lemma},
we see that 
\begin{equation} \label{detecting}
\dim \ker (\tilde{W} + I) = \dim (\ell_1 \cap \ell_2),
\end{equation}
which is the key relation in our computation of the Maslov 
index. Before properly defining the Maslov index, 
we note that the point $-1 \in S^1$ is chosen 
essentially at random, and any other point on 
$S^1$ would serve just as well. Indeed, we have the 
following proposition.

\begin{proposition} \label{X3}
Suppose $\mathbf{X}_1 = {X_1 \choose Y_1}$ and 
$\mathbf{X}_2 = {X_2 \choose Y_2}$ respectively denote frames for 
Lagrangian subspaces of $\mathbb{C}^{2n}$, $\ell_1$ and $\ell_2$, 
and let $\tilde{W}$ be as in (\ref{tildeW}). Given any value 
$\tilde{w} \in S^1$, set 
\begin{equation} \label{frame3rotation}
\mathbf{X}_3 = \mathbf{X}_3 (\tilde{w}) := i (1-\tilde{w}) \mathbf{X}_2 
- (1+\tilde{w}) J \mathbf{X}_2.
\end{equation} 
Then $\mathbf{X}_3$ is the frame for a Lagrangian subspace of $\mathbb{C}^{2n}$,
$\ell_3 = \ell_3 (\tilde{w})$, and 
\begin{equation*}
\dim \ker (\tilde{W} - \tilde{w}I) = \dim \ker (\mathbf{X}_1^* J \mathbf{X}_3)
= \dim (\ell_1 \cap \ell_3).
\end{equation*}
\end{proposition} 

\begin{proof} 
In order to check that $\mathbf{X}_3$ is Lagrangian, we note that by 
straightforward calculations we find 
\begin{equation*}
\mathbf{X}_3^* \mathbf{X}_3 = (|1-\tilde{w}|^2 + |1+\tilde{w}|^2) \mathbf{X}_2^* \mathbf{X}_2, 
\end{equation*}
and likewise $\mathbf{X}_3^* J \mathbf{X}_3 = 0$. It's clear from the first of 
these relations that $\dim \colspan (\mathbf{X}_3) = n$ and from the 
second that $\mathbf{X}_3$ satisfies the Lagrangian property.

Next, again by direct calculation, we find that 
\begin{equation*}
(X_3 - iY_3)(X_3+iY_3)^{-1} = -\frac{1}{\tilde{w}} (X_2 - iY_2)(X_2+iY_2)^{-1},
\end{equation*}
from which we see that $\tilde{w}$ is 
an eigenvalue of
\begin{equation*}
\tilde{W} = - (X_1 + iY_1)(X_1-iY_1)^{-1} (X_2 - iY_2)(X_2+iY_2)^{-1}
\end{equation*}
if and only if $-1$ is an eigenvalue of 
\begin{equation*}
\tilde{\mathcal{W}} = - (X_1 + iY_1)(X_1-iY_1)^{-1} (X_3 - iY_3)(X_3+iY_3)^{-1}.
\end{equation*}
In this way, 
\begin{equation*}
\dim (\ell_1 \cap \ell_3) = \dim \ker (\mathbf{X}_1^* J \mathbf{X}_3)
= \dim \ker (\tilde{\mathcal{W}}+I) = \dim \ker (\tilde{W}-\tilde{w}I). 
\end{equation*}
\end{proof}

\begin{remark} \label{two-ways}
We can interpret Proposition \ref{X3} in two useful ways. First, 
each eigenvalue of $\tilde{W}$, $\tilde{w} \in S^1$, indicates an 
intersection between $\ell_1$ and $\ell_3 (\tilde{w})$. In this way,
we can associate with any Lagrangian subspace $\ell_1$ a family of 
up to $n$ Lagrangian subspaces (depending on multiplicities) 
obtained through (\ref{frame3rotation}) as the Lagrangian subspaces
$\colspan (\mathbf{X}_3 (\tilde{w}))$ for some $\tilde{w}$ such 
that $\dim \ker (\tilde{W} - \tilde{w} I) \ne 0$. On the other hand, 
suppose we would like to move our 
spectral flow calculation from $-1$ to some other $\tilde{w} \in S^1$.
We let $\mathbf{X}_3$ denote our target (generally denoted $\mathbf{X}_2$),
and solve (\ref{frame3rotation}) for $\mathbf{X}_2$.
Then $\tilde{w} \in \sigma (\tilde{W})$ corresponds precisely
with intersections between $\ell_1$ and $\ell_3$, allowing us to 
compute the Maslov index as a spectral flow through $\tilde{w}$
(using $\tilde{W}$).  
\end{remark}

Turning now to our definition of the Maslov index, we note that 
since $\tilde{W}$ is unitary, we can define the Maslov index
in the $\mathbb{C}^{2n}$ setting precisely as in the $\mathbb{R}^{2n}$
setting in \cite{HLS2017}. For completeness, we sketch the 
development. 

Given two continuous maps $\ell_1 (t), \ell_2 (t)$ on a parameter
interval $\mathcal{I}$, we denote by $\mathcal{L}(t)$ the path 
\begin{equation*}
\mathcal{L} (t) = (\ell_1 (t), \ell_2 (t)).
\end{equation*} 
In what follows, we will define the Maslov index for the path 
$\mathcal{L} (t)$, which will be a count, including both multiplicity
and direction, of the number of times the Lagrangian paths
$\ell_1$ and $\ell_2$ intersect. In order to be clear about 
what we mean by multiplicity and direction, we observe that 
associated with any path $\mathcal{L} (t)$ we will have 
a path of unitary complex matrices as described in (\ref{tildeW}).
We have already noted that the Lagrangian subspaces $\ell_1$
and $\ell_2$ intersect at a value $t_* \in \mathcal{I}$ if and only 
if $\tilde{W} (t_*)$ has -1 as an eigenvalue. (Recall that we refer to the
value $t_*$ as a {\it crossing point}.) In the event of 
such an intersection, we define the multiplicity of the 
intersection to be the multiplicity of -1 as an eigenvalue of 
$\tilde{W} (t_*)$ (since $\tilde{W} (t_*)$ is unitary the algebraic and geometric
multiplicites are the same). When we talk about the direction 
of an intersection, we mean the direction the eigenvalues of 
$\tilde{W} (t_*)$ are moving (as $t$ increases) along the unit circle 
$S^1$ when they cross $-1$ (we take counterclockwise as the 
positive direction). We note that we will need to take care with 
what we mean by a crossing in the following sense: we must decide
whether to increment the Maslov index upon arrival or 
upon departure. Indeed, there are several different approaches 
to defining the Maslov index (see, for example, \cite{CLM, rs93}), 
and they often disagree on this convention. 

Following \cite{BF98, F, P96} (and in particular Definition 1.5 
from \cite{BF98}), we proceed by choosing a 
partition $a = t_0 < t_1 < \dots < t_n=b$ of $\mathcal{I} = [a,b]$, along 
with numbers $\{\epsilon_j\}_{j=1}^n \subset (0,\pi)$ so that 
$\ker\big(\tilde{W} (t) - e^{i (\pi \pm \epsilon_j)} I\big)=\{0\}$
for $t_{j-1} \le t \le t_j$; 
that is, $e^{i(\pi \pm \epsilon_j)} \in \mathbb{C} \setminus \sigma(\tilde{W} (t))$, 
for $t_{j-1} \le t \le t_j$ and $j=1,\dots,n$. 
Moreover, we notice that for each $j=1,\dots,n$ and any 
$t \in [t_{j-1},t_j]$ there are only 
finitely many values $\beta \in [0,\epsilon_j)$ 
for which $e^{i(\pi+\beta)} \in \sigma(\tilde{W} (t))$.

Fix some $j \in \{1, 2, \dots, n\}$ and consider the value
\begin{equation} \label{kdefined}
k (t,\epsilon_j) := 
\sum_{0 \leq \beta < \epsilon_j}
\dim \ker \big(\tilde{W} (t) - e^{i(\pi+\beta)}I \big).
\end{equation} 
for $t_{j-1} \leq t \leq t_j$. This is precisely the sum, along with multiplicity,
of the number of eigenvalues of $\tilde{W} (t)$ that lie on the arc 
\begin{equation*}
A_j := \{e^{i t}: t \in [\pi, \pi+\epsilon_j)\}.
\end{equation*}
The stipulation that 
$e^{i(\pi\pm\epsilon_j)} \in \mathbb{C} \setminus \sigma(\tilde{W} (t))$, for 
$t_{j-1} \le t \le t_j$
ensures that no eigenvalue can enter $A_j$ in the clockwise direction 
or exit in the counterclockwise direction during the interval $t_{j-1} \le t \le t_j$.
In this way, we see that $k(t_j, \epsilon_j) - k (t_{j-1}, \epsilon_j)$ is a 
count of the number of eigenvalues that enter $A_j$ in the counterclockwise 
direction (i.e., through $-1$) minus the number that leave in the clockwise direction
(again, through $-1$) during the interval $[t_{j-1}, t_j]$.


In dealing with the catenation of paths, it's particularly important to 
understand the difference $k(t_j, \epsilon_j) - k (t_{j-1}, \epsilon_j)$ if an 
eigenvalue resides at $-1$ at either $t = t_{j-1}$
or $t = t_j$ (i.e., if an eigenvalue begins or ends at a crossing). If an eigenvalue 
moving in the counterclockwise direction 
arrives at $-1$ at $t = t_j$, then we increment the difference forward, while if 
the eigenvalue arrives at -1 from the clockwise direction we do not (because it
was already in $A_j$ prior to arrival). On
the other hand, suppose an eigenvalue resides at -1 at $t = t_{j-1}$ and moves
in the counterclockwise direction. The eigenvalue remains in $A_j$, and so we do not increment
the difference. However, if the eigenvalue leaves in the clockwise direction 
then we decrement the difference. In summary, the difference increments forward upon arrivals 
in the counterclockwise direction, but not upon arrivals in the clockwise direction,
and it decrements upon departures in the clockwise direction, but not upon 
departures in the counterclockwise direction.      

We are now ready to define the Maslov index.

\begin{definition} \label{dfnDef3.6}  
Let $\mathcal{L} (t) = (\ell_1 (t), \ell_2 (t))$, where 
$\ell_1, \ell_2:\mathcal{I} \to \Lambda (n)$ 
are continuous paths in the Lagrangian--Grassmannian. 
The Maslov index $\mas(\mathcal{L}; \mathcal{I})$ is defined by
\begin{equation}
\mas(\mathcal{L}; \mathcal{I})=\sum_{j=1}^n(k(t_j,\epsilon_j)-k(t_{j-1},\epsilon_j)).
\end{equation}
\end{definition}

\begin{remark}
As we did in the introduction, we will typically refer explicitly to 
the individual paths with the notation $\mas (\ell_1, \ell_2; \mathcal{I})$.
\end{remark}

\begin{remark} As discussed in \cite{BF98}, the Maslov index does not depend
on the choices of $\{t_j\}_{j=0}^n$ and $\{\epsilon_j\}_{j=1}^n$, so long as 
these choices follow the specifications described above. Also, we emphasize that
Phillips' specification of the spectral flow allows for an infinite
number of crossings. In such cases, all except a finite number are 
necessarily transient (i.e., an eigenvalue crosses $-1$, but then crosses back, 
yielding no net contribution to the Maslov index).
\end{remark}

One of the most important features of the Maslov index is homotopy invariance, 
for which we need to consider continuously varying families of Lagrangian 
paths. To set some notation, we denote by $\mathcal{P} (\mathcal{I})$ the collection 
of all paths $\mathcal{L} (t) = (\ell_1 (t), \ell_2 (t))$, where 
$\ell_1, \ell_2: \mathcal{I} \to \Lambda (n)$ are continuous paths in the 
Lagrangian--Grassmannian. We say that two paths 
$\mathcal{L}, \mathcal{M} \in \mathcal{P} (\mathcal{I})$ are homotopic provided 
there exists a family $\mathcal{H}_s$ so that 
$\mathcal{H}_0 = \mathcal{L}$, $\mathcal{H}_1 = \mathcal{M}$, 
and $\mathcal{H}_s (t)$ is continuous as a map from $(t,s) \in \mathcal{I} \times [0,1]$
into $\Lambda (n) \times \Lambda (n)$. 
 
The Maslov index has the following properties. 

\medskip
\noindent
{\bf (P1)} (Path Additivity) If $\mathcal{L} \in \mathcal{P} (\mathcal{I})$
and $a, b, c \in \mathcal{I}$, with $a < b < c$, then 
\begin{equation*}
\mas (\mathcal{L};[a, c]) = \mas (\mathcal{L};[a, b]) + \mas (\mathcal{L}; [b, c]).
\end{equation*}

\medskip
\noindent
{\bf (P2)} (Homotopy Invariance) If $\mathcal{I} = [a, b]$,
$a < b$,
and $\mathcal{L}, \mathcal{M} \in \mathcal{P} (\mathcal{I})$ 
are homotopic with $\mathcal{L} (a) = \mathcal{M} (a)$ and  
$\mathcal{L} (b) = \mathcal{M} (b)$ (i.e., if $\mathcal{L}, \mathcal{M}$
are homotopic with fixed endpoints) then 
\begin{equation*}
\mas (\mathcal{L};[a, b]) = \mas (\mathcal{M};[a, b]).
\end{equation*} 

Straightforward proofs of these properties appear in \cite{HLS2017}
for Lagrangian subspaces of $\mathbb{R}^{2n}$, and proofs in the current setting of 
Lagrangian subspaces of $\mathbb{C}^{2n}$ are essentially identical.

\section{Direction of Rotation} \label{rotation_section}

As noted in the previous section, the direction we associate with a 
crossing point is determined by the direction in which eigenvalues
of $\tilde{W}$ rotate through $-1$ (counterclockwise is positive, 
while clockwise is negative). When analyzing the Maslov index, we 
need a convenient framework for analyzing this direction, and 
the development of such a framework is the goal of this section. 

First, in order to understand monotonicity as the spectral parameter $\lambda$
evolves, we can use the following lemma from \cite{HLS2017}. 
(In \cite{HLS2017}, the statement takes
the frames to be $C^1$, but the proof only requires differentiability,
as asserted here.)

\begin{lemma} \label{monotonicity-lemma1}
Suppose $\ell_1, \ell_2: \mathcal{I} \to \Lambda (n)$ denote paths of 
Lagrangian subspaces of $\mathbb{C}^{2n}$ with respective frames 
$\mathbf{X}_1 = {X_1 \choose Y_1}$ and $\mathbf{X}_2 = {X_2 \choose Y_2}$
that are differentiable at $t_0 \in \mathcal{I}$. If the matrices 
\begin{equation*}
- \mathbf{X}_1 (t_0)^* J \mathbf{X}_1' (t_0) = X_1 (t_0)^* Y_1' (t_0) - Y_1 (t_0)^* X_1'(t_0)
\end{equation*}
and (noting the sign change)
\begin{equation*}
\mathbf{X}_2 (t_0)^* J \mathbf{X}_2' (t_0) = - (X_2 (t_0)^* Y_2' (t_0) - Y_2 (t_0)^* X_2'(t_0))
\end{equation*}
are both non-negative, and at least one is positive definite, then the eigenvalues of 
$\tilde{W} (t)$ rotate in the counterclockwise direction as $t$ increases through $t_0$. 
Likewise, if both of these matrices are non-positive, and at least one is 
negative definite, then the eigenvalues of $\tilde{W} (t)$ rotate in the clockwise direction as 
$t$ increases through $t_0$.
\end{lemma}

In order to analyze monotonicity as the independent variable $x$
evolves, we require a more detailed analysis, and for this we'll 
consider a general linear Hamiltonian system 
\begin{equation} \label{general-linear-hammy}
    Jy' = \mathpzc{B} (t) y, 
    \quad t \in (0,1), \quad y(t) \in \mathbb{C}^{2n},
\end{equation}
for which we assume $\mathpzc{B} \in L^1 ((0,1), \mathbb{C}^{2n})$,
and that $\mathpzc{B} (t)$ is self-adjoint for a.e. $t \in (0, 1)$. 
Throughout this discussion, we will let $\mathbf{X} (t) = {X (t) \choose Y(t)}$ 
denote a $2n \times n$ matrix solution of (\ref{general-linear-hammy}), 
and we will assume that for some $t_0 \in [0, 1]$ $\mathbf{X} (t_0)$
is a frame for a Lagrangian subspace, from which it follows that 
$\mathbf{X} (t)$ is the frame for a Lagrangian subspace $\ell (t)$ 
for all $t \in [0, 1]$. In addition, we will denote by 
$\tilde{\ell}$ a fixed Lagrangian target with frame 
$\tilde{\mathbf{X}} = {\tilde{X} \choose \tilde{Y}}$, and consider the Maslov
index 
\begin{equation}
    \mas (\ell (\cdot), \tilde{\ell}; [0, 1]),
\end{equation}
which can be computed as described above with the matrix 
\begin{equation} \label{mathcal-tilde-W}
\tilde{\mathcal{W}} (t) := - (X (t) + i Y (t)) (X (t) - i Y (t))^{-1}
(\tilde{X} - i \tilde{Y}) (\tilde{X} + i \tilde{Y})^{-1}.
\end{equation}

As discussed in \cite{HS2016}, if we fix any $t_0 \in [0, 1]$, 
then we can write 
\begin{equation*}
\tilde{\mathcal{W}} (t) = \tilde{\mathcal{W}} (t_0) e^{i R(t)}
\end{equation*}
for $t$ sufficiently close to $t_0$. Here, $i R(t)$ is
the logarithm of $\tilde{\mathcal{W}}(t_0)^{-1} \tilde{\mathcal{W}} (t)$, and 
we clearly have $R(t_0) = 0$. (Since $\tilde{\mathcal{W}}(t_0)^{-1} \tilde{\mathcal{W}} (t_0) = I$,
it's clear that $\tilde{\mathcal{W}}(t_0)^{-1} \tilde{\mathcal{W}} (t)$ has a unique
logarithm for $t$ sufficiently close to $t_0$.) For any $t \in [0, 1]$, we set 
$\tilde{\Omega} (t) := -i \tilde{\mathcal{W}} (t)^{-1} \tilde{\mathcal{W}}'(t)$,
so that 
\begin{equation*}
\tilde{\mathcal{W}}'(t) = i \tilde{\mathcal{W}} (t) \tilde{\Omega} (t),
\end{equation*}
for all $t \in [0, 1]$. Comparing expressions, we see that 
$R' (t_0) = \tilde{\Omega} (t_0)$. 

As discussed particularly in the proof of Lemma 3.11 in \cite{HS2016}, 
the direction of rotation for eigenvalues of $\tilde{\mathcal{W}} (t)$ as 
they cross $\tilde{w}_0 \in \sigma (\tilde{\mathcal{W}} (t_0))$ is determined by 
the nature of $R(t)$ for $t$ near $t_0$. For the current analysis, we 
will require an extension of Lemma 3.11 in \cite{HS2016}, and 
in developing this we will repeat part of the argument 
from \cite{HS2016}. 

First, following \cite{F}[p. 306], we fix any $\theta \in [0,2 \pi)$
so that $e^{i\theta} \notin \ker (\tilde{\mathcal{W}} (t_0))$ and, for $t$ 
sufficiently close to $t_0$, define the auxiliary matrix 
\begin{equation} \label{Adefined}
\tilde{A} (t) := i (e^{i \theta} I - \tilde{\mathcal{W}} (t))^{-1} (e^{i \theta} I + \tilde{\mathcal{W}} (t)).
\end{equation} 
It is straightforward to check that $\tilde{A} (t)$ is self-adjoint, and this 
allows us to conclude that its eigenvalues will all be real-valued.
If we denote the eigenvalues of $\tilde{\mathcal{W}} (t)$ by $\{\tilde{w}_j (t)\}_{j=1}^n$,
and the eigenvalues of $\tilde{A}(t)$ by $\{\tilde{a}_j (t)\}_{j=1}^n$ then by 
spectral mapping we have the correspondence 
\begin{equation} \label{a-w-connection}
\tilde{a}_j (t) = i \frac{e^{i \theta} + \tilde{w}_j (t)}{e^{i \theta} - \tilde{w}_j (t)}.
\end{equation}
By a short argument in the proof of Lemma 3.11 in \cite{HS2016}, 
the authors find that if the eigenvalue $\tilde{a}_j (t)$ is increasing
as $t$ increases through $t_0$ then the corresponding eigenvalue 
$\tilde{w}_j (t)$ of $\tilde{W} (t)$ will rotate in the 
counterclockwise direction along 
$S^1$ as $t$ increases through $t_0$. This means that the rotation of 
the eigenvalues of $\tilde{\mathcal{W}} (t)$ can be determined by the 
linear motion of the eigenvalues of $\tilde{A} (t)$. 

Suppose $\tilde{w}_0$ is an eigenvalue of $\tilde{\mathcal{W}} (t_0)$
with multiplicity $m$. Since $\tilde{\mathcal{W}} (t_0)$ is unitary, the 
algebraic and geometric multiplicities of its eigenvalues agree, 
so the eigenspace associated with $\tilde{w}_0$, which we denote $\tilde{V}_0$, 
has dimension $m$. From Theorem II.5.4 in \cite{Kato}, we know there exists a 
corresponding eigenvalue group 
$\{\tilde{w}_j (t)\}_{j=1}^m \subset \sigma (\tilde{\mathcal{W}} (t))$
so that $\tilde{w}_j (t_0) = \tilde{w}_0$ for $j = 1, 2, ..., m$. 
By a natural choice 
of indexing, each such $\tilde{w}_j (t)$ will have a corresponding eigenvalue 
$\tilde{a}_j (t) \in \sigma (\tilde{A}(t))$, and the eigenspace associated with 
$\{\tilde{a}_j (t_0)\}_{j=1}^m \subset \sigma (\tilde{A}(t_0))$ will be $\tilde{V}_0$. 
In particular, it follows from the discussion above that 
the rotation of the eigenvalues   
$\{\tilde{w}_j (t)\}_{j=1}^m \subset \sigma (\tilde{\mathcal{W}} (t))$ 
through $\tilde{w}_0$ will be determined by the linear motion 
of the eigenvalues $\{\tilde{a}_j (t)\}_{j=1}^m \subset \sigma (\tilde{A}(t))$ through 
$\tilde{a}_0$.

In order to apply these observations in the current setting, we will 
use the following lemma, which is based on Theorem II.5.4 from \cite{Kato}. 

\begin{lemma} \label{eigenvalue-lemma}
Let $A \in AC ([0, 1], \mathbb{C}^{n \times n})$, with 
$A(t)$ self-adjoint for each $t \in [0, 1]$. Fix $t_0 \in [0, 1]$,
and suppose there exists $\delta > 0$ sufficiently small so that 
$A' (t)$ is non-negative (resp. non-positive) for a.e. 
$t \in (t_0 - \delta, t_0 + \delta) \cap [0, 1]$.
Then the $n$ eigenvalues of $A(t)$ must be non-decreasing 
(resp. non-increasing)
on $(t_0 - \delta, t_0 + \delta) \cap [0, 1]$. 
\end{lemma}

\begin{proof} Since $A \in AC ([0, 1], \mathbb{C}^{n \times n})$, 
we have that $A$ is differentiable a.e. in $(0, 1)$. Suppose $A$ 
is differentiable at a value $\tau \in (0, 1)$, and let $a (\tau)$ 
denote any eigenvalue of $A (\tau)$. If $m$ denotes the multiplicity 
of $a (\tau)$ as an eigenvalue of $A (\tau)$, then 
we have from Theorem II.5.4
in \cite{Kato} that there will correspond an eigenvalue 
group $\{a_j (t)\}_{j=1}^m \subset \sigma (A(t))$ that can 
be expressed as 
\begin{equation} \label{expansion1}
a_j (t) = a (\tau) + \alpha_j^{\tau} (t - \tau) 
+ \mathbf{o} (|t-\tau|),
\quad \forall \, j \in \{1, 2, \dots, m\},
\end{equation}
for $t$ sufficiently close to $\tau$, and where the values 
$\{\alpha_j^{\tau}\}_{j=1}^m$ are eigenvalues of 
$P_\tau A'(\tau) P_\tau$ in the space $P_{\tau} \mathbb{C}^{2n}$,
with $P_\tau$ denoting projection onto the eigenspace
of $A(\tau)$ associated to $a (\tau)$. 

Since $A' (t)$ is non-negative 
for a.e. $t \in (t_0 - \delta, t_0 + \delta) \cap [0, 1]$, we can 
conclude that the values $\{\alpha_j^{\tau}\}_{j=1}^m$ must be 
non-negative for a.e. $\tau \in (t_0 - \delta, t_0 + \delta) \cap [0, 1]$. 
In particular, we see that 
\begin{equation*}
a_j' (\tau) = \alpha_j^{\tau} \ge 0; 
\quad \text{a.e. } \tau \in (t_0 - \delta, t_0 + \delta) \cap [0, 1].
\end{equation*} 
Upon integrating this last relation on any 
interval $(t_1, t_2) \subset  (t_0 - \delta, t_0 + \delta) \cap [0, 1]$,
we see that 
\begin{equation*}
    a_j (t_2) \ge a_j (t_1),
    \quad \forall \, j \in \{1, 2, \dots, m\}.
\end{equation*}
Since this is true for all eigenvalue groups of $A(t)$, 
the proof is complete for the case in which $A' (t)$ is non-negative  
for a.e. $t \in (t_0 - \delta, t_0 + \delta) \cap [0, 1]$. The 
case in which $A' (t)$ is non-positive 
for a.e. $t \in (t_0 - \delta, t_0 + \delta) \cap [0, 1]$ can be 
established similarly.
\end{proof}

In our analysis, monotonicity as the independent varariable $x$ evolves
will be a consequence of the following lemma. 

\begin{lemma} \label{monotonicity-in-x-lemma}
For $\tilde{\mathcal{W}} (t)$ as specified in (\ref{mathcal-tilde-W}),
suppose $\mathpzc{B} (t)$ is non-negative (resp. non-positive) for a.e. $t \in (0, 1)$. 
Then no eigenvalue of $\tilde{\mathcal{W}} (t)$ 
can rotate in the counterclockwise (resp. clockwise) direction on 
any interval $[a, b] \subset [0, 1]$, $a < b$.
\end{lemma}

\begin{proof}
In Section 4 of \cite{HLS2017}, the authors show that if 
$\tilde{A} (t)$ from (\ref{Adefined}) is differentiable at $t$, then  
\begin{equation} \label{Aprime}
\tilde{A}'(t) = 2 \Big( (e^{i \theta} I - \tilde{\mathcal{W}} (t))^{-1} \Big)^*
\tilde{\Omega} (t) (e^{i \theta} I - \tilde{\mathcal{W}} (t))^{-1},
\end{equation}
with (under the assumptions of the current lemma)
\begin{equation} \label{omega-tilde}
    \tilde{\Omega} (t)
    = - 2 ((X (t) - i Y (t))^{-1} \tilde{\mathpzc{W}})^* 
    \mathbf{X} (t)^* \mathpzc{B} (t)  \mathbf{X} (t) ((X (t) - i Y (t))^{-1} \tilde{\mathpzc{W}}),
\end{equation}
where $\tilde{\mathpzc{W}}$ denotes the constant matrix
\begin{equation*}
    \tilde{\mathpzc{W}} :=  (\tilde{X} - i \tilde{Y}) (\tilde{X} + i \tilde{Y})^{-1}.
\end{equation*}
(In \cite{HLS2017}, the authors take $\mathbf{X} (t), \tilde{\mathbf{X}} \in \mathbb{R}^{2n \times n}$,
but the calculation in our setting simply replaces transpose with adjoint where appropriate.)

If $\mathpzc{B} (t)$ is non-negative for a.e. $t \in (0, 1)$, then $\tilde{\Omega} (t)$ is
non-positive for a.e. $t \in (0,1)$, and consequently $\tilde{A}'(t)$ is non-positive for 
a.e. $t \in (0, 1)$. It follows from Lemma \ref{eigenvalue-lemma} that the eigenvalues
of $\tilde{A} (t)$ must be non-increasing on $(0, 1)$, and using (\ref{a-w-connection})
we can conclude that no eigenvalue of $\tilde{\mathcal{W}} (t)$ 
can rotate in the counterclockwise direction on any interval $[a, b] \subset [0, 1]$,
$a < b$.
\end{proof}

Since $\mathpzc{B} (t)$ need not be strictly positive for all $t \in (0,1)$, 
Lemma \ref{monotonicity-in-x-lemma} leaves open the possibility that there
exists an interval $[a, b] \subset [0, 1]$, $a < b$, so that 
\begin{equation} \label{intersection}
    \dim (\ell (t) \cap \tilde{\ell}) \ne 0
\end{equation}
for all $t \in [a, b]$. (E.g., $\mathpzc{B} (t) \equiv 0$ would be a trivial example 
allowing this possibility.) In our applications we will have a slightly stronger condition
than (\ref{intersection}), namely that there 
exists some $m \in \{1, 2, \dots, n\}$ so that 
\begin{equation} \label{intersection-m}
\dim (\ell (t) \cap \tilde{\ell}) = m
\end{equation}
for all $t \in [a, b]$. In particular, we will make use of the 
following lemma.

\begin{lemma} \label{continuation-lemma1}
Let $\mathbf{X} (t) = {X (t) \choose Y (t)}$ be a continuous
matrix function $\mathbf{X}: [0, 1] \to \mathbb{C}^{2n \times n}$
such that for each $t \in [0, 1]$, $\mathbf{X} (t)$ is the frame
for a Lagrangian subspace $\ell (t)$, and let 
$\tilde{\mathbf{X}} \in \mathbb{C}^{2n \times n}$ 
denote any fixed Lagrangian frame. Suppose that as $t$ increases
from $0$ to $1$ any eigenvalue of $\tilde{\mathcal{W}} (t)$ 
(from (\ref{mathcal-tilde-W}), except under the current 
assumptions on $\mathbf{X} (t)$ and $\tilde{\mathbf{X}}$) 
that crosses $-1$ makes the crossing in the counterclockwise 
direction. If there exists an interval $[a, b] \subset [0, 1]$,
$a < b$, so that $\dim \ker (\tilde{\mathcal{W}}(t) + I) \ne 0$ 
for all $t \in [a, b]$, then there exists an integer 
$m \in \{1, 2, \dots, n \}$ and a subinterval 
$[c, d] \subset [a, b]$, $c < d$, so that 
$\dim \ker (\tilde{\mathcal{W}}(t) + I) = m$ for all $t \in [c, d]$.
The same conclusion holds true under the alternative assumption that 
as $t$ increases from $0$ to $1$ any eigenvalue of $\tilde{\mathcal{W}} (t)$
that crosses $-1$ makes the crossing in the clockwise 
direction. 
\end{lemma}

\begin{proof} Fix any $t_0 \in (a, b)$, and observe
that we necessarily have
\begin{equation*}
\dim \ker (\tilde{\mathcal{W}} (t_0) + I) = m_0
\end{equation*}
for some $m_0 \in \{1, 2, \dots, n\}$. Using continuity 
of $\tilde{\mathcal{W}} (t)$, fix $\delta_0 > 0$
sufficiently small so that none of the eigenvalues of 
$\tilde{\mathcal{W}} (\cdot)$ arrive at 
$-1$ during the interval $[t_0, t_0 + \delta_0]$. (Note
that, by assumption, arrivals would necessarily occur
from the counterclockwise direction, and also that we don't count
an eigenvalue that has remained at $-1$ for the full 
interval $[t_0, t_0 + \delta_0]$ as arriving at $-1$.)

As $t$ increases from $t_0$, one or more eigenvalues 
residing at $-1$ could rotate away from $-1$ in the 
counterclockwise direction, though the total 
number of eigenvalues residing at $-1$ could not be
reduced to $0$ (by assumption). Since there are only 
a finite number of eigenvalues, these departure times 
must be separated by intervals, and so we can take 
$[c, d]$ to be any interval between departures.
\end{proof}  

According to Lemma \ref{dimensions_lemma}, we can 
re-state (\ref{intersection-m}) in the following way: 
the matrix Wronskian 
\begin{equation*}
\mathcal{W} (t) := \tilde{\mathbf{X}}^* J \mathbf{X} (t)
\end{equation*}
satisfies $\dim \ker \mathcal{W} (t) = m$ for all $t \in [a, b]$.
In particular, $0$ is an eigenvalue of $\mathcal{W} (t)$ with 
multiplicity exactly $m$ for all $t \in [a, b]$. 
In the terminology of \cite{Kato}, there is an eigenvalue
group associated with $0$, with each eigenvalue in 
the group identically $0$ for all $t \in [a, b]$,
and such that the associated projection $P(t)$ 
projects onto the $m$-dimensional space $\ker \mathcal{W} (t)$ for all 
$t \in [a, b]$. According to a slight extension 
of Theorem II.5.4 in \cite{Kato}, if $\mathbf{X} (t)$ 
is absolutely continuous on $(a, b)$ then $P(t)$ 
will be absolutely continuous on this interval as well. 

\begin{lemma} \label{technical-lemma}
For $\tilde{\mathcal{W}} (t)$ as specified in (\ref{mathcal-tilde-W}),
suppose either that $\mathpzc{B} (t)$ is non-negative for a.e. $t \in (0, 1)$
or that $\mathpzc{B} (t)$ is non-positive for a.e. $t \in (0, 1)$. If 
there exists an interval $[a, b] \subset [0, 1]$ and an 
integer $m \in \{1, 2, \dots, n\}$ so that 
$\dim \ker (\tilde{\mathcal{W}} (t) + I) = m$ for all 
$t \in [a, b]$ then the following hold.

\medskip
\noindent
(i) There exists a function $v \in AC ([a, b], \mathbb{C}^n)$ so that 
$\mathbf{X} (t) v (t) \in \ell (t) \cap \tilde{\ell}$ for all 
$t \in [a, b]$; 

\medskip
\noindent
(ii) Given any $v (t)$ such that $\mathbf{X} (t) v (t) \in \ell (t) \cap \tilde{\ell}$ for all 
$t \in [a, b]$ (not necessarily the choice of $v (t)$ from Item (i)), there 
exists $w \in AC ([a, b], \mathbb{C}^n)$ to that 
\begin{equation*}
    \mathbf{X} (t) v (t) = \tilde{\mathbf{X}} w(t),
    \quad \forall \,\, t \in (a, b),
\end{equation*}
and additionally
\begin{equation*}
    \mathbf{X} (t) v' (t) = \tilde{\mathbf{X}} w'(t),
    \quad \forall \,\, t \in (a, b).
\end{equation*}
\end{lemma}

\begin{proof}
For Item (i), we fix some $t_* \in (a,b)$ and observe
that, by assumption, 
$\dim (\ell (t_*) \cap \tilde{\ell}) = m$. It follows
that there exists a vector $v_* \in \mathbb{C}^n$
so that 
\begin{equation*}
    \mathbf{X} (t_*) v_* \in \ell (t_*) \cap \tilde{\ell}.
\end{equation*}
Letting $P (t)$ denote orthogonal projection 
onto $\ker \mathcal{W} (t)$, we set 
$v (t) := P (t) v_*$. Since $P \in \AC ([a, b], \mathbb{C}^n)$,
we can conclude that $v \in \AC ([a, b], \mathbb{C}^n)$.
In addition, $\mathbf{X} (t) v(t)$ is clearly in 
$\ell (t)$ for all $t \in [a, b]$, and 
since $\mathcal{W} (t) v(t) = 0$ for all 
$t \in [a, b]$, we have $\tilde{\mathbf{X}}^* J \mathbf{X} (t) v(t) = 0$
for all $t \in [a, b]$, from which we can conclude by 
Lagrangian maximality that 
$\mathbf{X} (t) v(t) \in \tilde{\ell}$ for all $t \in [a, b]$.
This complete the proof of Item (i). 

Turning to Item (ii), let $v$ denote any vector function 
$v \in AC ([a, b], \mathbb{C}^n)$ so that 
$\mathbf{X} (t) v (t) \in \ell (t) \cap \tilde{\ell}$ for all 
$t \in [a, b]$. For each $t \in [a, b]$ there must 
exist some $w = w(t)$ so that 
\begin{equation*}
    \mathbf{X} (t) v (t) = \tilde{\mathbf{X}} w (t),
\end{equation*}
and we can use the Moore-Penrose pseudoinverse
of $\tilde{\mathbf{X}}$, to write 
\begin{equation*}
w (t) = (\tilde{\mathbf{X}}^* \tilde{\mathbf{X}})^{-1} \tilde{\mathbf{X}}^* \mathbf{X} (t) v(t).
\end{equation*}
We see from this expression that $w (t)$ must also be 
absolutely continuous on $[a, b]$. 

Next, by direct calculation we find 
\begin{equation*}
\mathcal{W}'(t) = \tilde{\mathbf{X}}^* J \mathbf{X}' (t)
= \tilde{\mathbf{X}}^* \mathpzc{B} (t) \mathbf{X} (t).
\end{equation*}
We are justified in differentiating the relation 
$\mathcal{W} (t) v(t) = 0$, and we obtain 
\begin{equation*}
\mathcal{W}'(t) v(t) + \mathcal{W}(t) v' (t) = 0,
\quad {\rm a.e.} \,\, t \in (a, b),
\end{equation*} 
and this relation can be expressed as 
\begin{equation} \label{weaker}
\tilde{\mathbf{X}}^* \mathpzc{B} (t) \mathbf{X} (t) v(t)
+ \tilde{\mathbf{X}}^* J \mathbf{X} (t) v' (t) = 0,
\quad {\rm a.e.} \,\, t \in (a, b).
\end{equation}
If we take a $\mathbb{C}^n$ inner product of this equation with $w(t)$, 
we find that 
\begin{equation} \label{inner-product}
(\tilde{\mathbf{X}}^* \mathpzc{B} (t) \mathbf{X} (t) v(t), w(t))
+ (\tilde{\mathbf{X}}^* J \mathbf{X} (t) v' (t), w(t)) = 0,
\quad {\rm a.e.} \,\, t \in (a, b).
\end{equation}
For the second in this last expression, we see that 
\begin{equation*}
\Big(\tilde{\mathbf{X}}^* J \mathbf{X} (t) v'(t), w(t) \Big) 
= - \Big(v'(t), \mathbf{X} (t)^* J \tilde{\mathbf{X}} w(t) \Big)
= 0,
\end{equation*}
because $\mathbf{X} (t)^* J \tilde{\mathbf{X}} (t) w(t) = 0$ 
(since $\tilde{\mathbf{X}} (t) w(t) \in \ell (t)$). 
This leaves only a single summand in (\ref{inner-product}),
and we can write 
\begin{equation*}
    0 = (\tilde{\mathbf{X}}^* \mathpzc{B} (t) \mathbf{X} (t) v(t), w(t))
    = (\mathpzc{B} (t) \mathbf{X} (t) v(t), \tilde{\mathbf{X}} w(t)). 
\end{equation*}
Using the relation $\mathbf{X} (t) v (t) = \tilde{\mathbf{X}} w (t)$, 
we can express this as 
\begin{equation*}
    (\mathpzc{B} (t) \mathbf{X} (t) v(t), \mathbf{X} (t) v(t)) = 0,
    \quad  {\rm a.e.} \,\, t \in (a, b). 
\end{equation*}
Since $\mathpzc{B} (t)$ is non-negative for a.e. $t \in (a, b)$ 
(or, alternatively, non-positive for a.e. $t \in (a, b)$), 
we must have $\mathbf{X} (t) v(t) \in \ker \mathpzc{B} (t)$ for 
a.e. $t \in (a, b)$. This allows us to additionally compute
\begin{equation*}
    J \mathbf{X}' (t) v (t) = \mathpzc{B} (t) \mathbf{X}' (t) v(t) = 0,
\end{equation*}
from which we can conclude that $\mathbf{X}' (t) v (t) = 0$ for 
a.e. $t \in (a, b)$, and consequently 
$\mathbf{X} (t) v' (t) = \tilde{\mathbf{X}} w' (t)$ for a.e. 
$t \in (a, b)$, establishing the final claim of the lemma. 
\end{proof}

\begin{lemma} \label{continuation-lemma2}
For $\tilde{\mathcal{W}} (t)$ as specified in (\ref{mathcal-tilde-W}),
suppose either that $\mathpzc{B} (t)$ is non-negative for a.e. $t \in (0, 1)$
or that $\mathpzc{B} (t)$ is non-positive for a.e. $t \in (0, 1)$. 
If there exists an interval $[a, b] \subset [0, 1]$ and an 
integer $m \in \{1, 2, \dots, n\}$ so that 
$\dim \ker (\tilde{\mathcal{W}} (t) + I) = m$ for all 
$t \in [a, b]$ then there exists an interval 
$[c, d] \subset [a, b]$, $c < d$, and a constant vector 
$v_0 \in \mathbb{C}^n \backslash \{0\}$ so that 
$\mathbf{X} (t) v_0 \in \ell (t) \cap \tilde{\ell}$
for all $t \in [c, d]$. 
\end{lemma}

\begin{proof} First, if $m = n$, then we must 
have $\ell (t) = \tilde{\ell}$ for all $t \in [a, b]$, 
and so for {\it any} $v_0 \in \mathbb{C}^n \backslash \{0\}$, 
we have $\mathbf{X} (t) v_0 \in \ell (t) \cap \tilde{\ell}$
for all $t \in [a, b]$. In particular, in this case, the claim
holds for $[c, d] = [a, b]$.

Next, suppose $m \in \{1, 2, \cdots, n-1\}$. 
If $\dim (\ell (t) \cap \tilde{\ell}) = m$ for all 
$t \in [a, b]$, then we can fix some $t_* \in (a, b)$ and let 
$\{v_j^*\}_{j=1}^m$ denote a basis for 
$\ker \tilde{\mathbf{X}}^* J \mathbf{X} (t_*)$. 
If, as in the proof of Lemma \ref{technical-lemma}, we let  
$P(t)$ denote projection onto $\ker \tilde{\mathbf{X}}^* J \mathbf{X} (t)$,
then the elements $v_j (t) := P(t) v_j^*$, $j = 1, 2, \dots, m$, comprise a 
collection of vector functions $v_j \in \AC ([a, b], \mathbb{C}^n)$  
that are linearly independent for $t$ sufficiently close to $t_*$.
We denote this interval of linear independence $I_*$. 
It follows immediately that the collection of vector functions 
$\{\mathbf{X} (t) v_j (t)\}_{j=1}^m$ forms a basis for
$\ell (t) \cap \tilde{\ell}$ for each $t \in I_*$. (Here, 
$\mathbf{X} (t) v_j (t) \in \tilde{\ell}$ because 
$\tilde{\mathbf{X}}^* J \mathbf{X} (t) v_j (t) = 0$.)

We observe that if we let $V(t)$ denote the $n \times m$ matrix
with columns $\{v_j (t)\}_{j=1}^m$, then the columns of 
$\mathbf{X} (t) V(t)$ are precisely the basis elements for 
$\ell (t) \cap \tilde{\ell}$ selected above. In the usual way, we can 
make a change of basis by multiplying $\mathbf{X} (t) V(t)$
on the right by any non-singular $m \times m$ matrix 
$M(t)$. We claim that by restricting to a smaller 
interval if necessary, we can choose $M(t)$ so that 
$\tilde{V}(t) := V(t) M(t)$ is in column reduced echelon form for 
all $t$ in the smaller interval. 

In order to understand this claim, we think as follows. 
We begin with the first components of the vectors 
$\{v_j (t)\}_{j=1}^m$, which we will denote
$\{v_{1j} (t)\}_{j=1}^m$. If each of these is 0
on the entirety of $I_*$, then the entire first row 
of $\tilde{V} (t)$ will be 0. Otherwise, there exists
some index $j \in \{1, 2, \dots, m\}$ and some value
$t_1 \in I_*$ so that 
$v_{1j} (t_1) \ne 0$. For notational convenience, let's
suppose $j = 1$, so that $v_{11} (t_1) \ne 0$. Then 
by continuity there exists some interval 
$I_1 \subset I_*$ so that 
$v_{11} (t) \ne 0$ for all $t \in I_1$. This allows 
us to divide each entry in the column $v_1 (t)$ by the entry  
$v_{11} (t)$, and consquently we can perform column
operations to eliminate the first component of each
column $\{v_j (t)\}_{j=2}^m$. 

Following the preceding operations, the second 
component in the second column of the resulting matrix
can be expressed as 
\begin{equation} \label{entry22}
v_{22} (t) - v_{12} (t) \frac{v_{21} (t)}{v_{11} (t)}.
\end{equation}
If this quantity is $0$ for all $t \in I_1$, then this entry will 
be $0$ in $\tilde{V} (t)$, and we move to the
third entry in the second column. If (\ref{entry22})
is not $0$ on the entirety of $I_1$, then by continuity
there exists an interval $I_2 \subset I_1$ so that 
(\ref{entry22}) is non-zero on the entirety of $I_2$.
In this case, we divide the second column of our matrix
by (\ref{entry22}) and use this pivot to eliminate 
the remaining entries in the second row of our new 
matrix. Continuing in this way, we obtain a matrix 
$\tilde{V} (t)$ on some final (smallest) interval 
$I$, whose columns span the same space as the columns of 
$V(t)$, and which has at least $m$ rows with a single 
$1$ as the only non-zero entry. 

Let $\{\tilde{v}_j (t)\}_{j=1}^m$ denote the columns of 
this new matrix, and set 
\begin{equation*}
\mathbf{v} (t) := \sum_{j=1}^m \tilde{v}_j (t). 
\end{equation*}  
Then, in particular, we have 
$\mathbf{X} (t) \mathbf{v} (t) \in \ell (t) \cap \tilde{\ell}$
for all $t \in I$. Since $\mathpzc{B} (t)$ is non-negative
for a.e. $t \in I$,  
we can conclude from Lemma \ref{technical-lemma} that 
there exists some $w \in AC (I, \mathbb{C}^n)$ so 
that $\mathbf{X} (t) v (t) = \tilde{\mathbf{X}} w(t)$
and $\mathbf{X} (t) v' (t) = \tilde{\mathbf{X}} w'(t)$
both hold for a.e. $t \in I$. So, in particular, 
$\mathbf{X} (t) \mathbf{v}' (t) \in \ell (t) \cap \tilde{\ell}$
for a.e. $t \in I$. If $\mathbf{v}'(t) = 0$ for a.e. 
$t \in I$, then we are done, because we can take 
$v_0 := \mathbf{v} (t)$, which is constant in $t$. (We
recall that $\mathbf{v} (t)$
is absolutely continuous.) If 
$\mathbf{v}'(t)$ is not $0$ for a.e. $t \in I$, then it 
must be linearly independent of the set $\{\tilde{v}_j (t)\}_{j=1}^m$,
because each of the $\{\tilde{v}_j (t)\}_{j=1}^m$ will have
a non-zero entry in a row where $\mathbf{v}'(t)$ has only zeros. 
But this means $\ell (t) \cap \tilde{\ell}$ has dimension at least $m+1$,
a contradiction. 
\end{proof}

\begin{lemma} \label{continuation-lemma3}
For $\tilde{\mathcal{W}} (t)$ as specified in (\ref{mathcal-tilde-W}),
suppose either that $\mathpzc{B} (t)$ is non-negative for a.e. $t \in (0, 1)$
or that $\mathpzc{B} (t)$ is non-positive for a.e. $t \in (0, 1)$. 
If there exists an interval $[a, b] \subset [0, 1]$, $a < b$, 
so that $\mathbf{X} (t) v_0 \in \ell (t) \cap \tilde{\ell}$ 
for all $t \in [a, b]$, then $\mathbf{X} (t) v_0$ is constant
on $[a, b]$. 
\end{lemma}

\begin{proof}
We know from Lemma \ref{technical-lemma} that there exists 
$w \in \AC ([a, b], \mathbb{C}^n)$ so that 
\begin{equation*}
\mathbf{X} (t) v_0 = \tilde{\mathbf{X}} w(t), \quad \forall \, t \in [a, b],
\end{equation*}
with also 
\begin{equation*}
0 = \mathbf{X} (t) \frac{d}{dt} v_0 = \tilde{\mathbf{X}} w'(t).
\end{equation*}
We can conclude that $w' (t) = 0$ for a.e. $t \in (a, b)$, 
and so $w (t) = w_0$ for some fixed $w_0 \in \mathbb{C}^n \backslash \{0\}$.
But then 
\begin{equation*}
\mathbf{X} (t) v_0 = \tilde{\mathbf{X}} w_0, \quad \forall \, t \in [a, b],
\end{equation*}
giving the claim. 
\end{proof}

\section{Proofs of Theorems \ref{bc1_theorem} and \ref{bc2_theorem}} \label{proofs}

In this section, we use our Maslov index framework to prove our two 
main theorems. 

\subsection{Proof of Theorem \ref{bc1_theorem}}

Fix any pair $\lambda_1, \lambda_2 \in I$, with $\lambda_1 < \lambda_2$,
and let $\ell_1 (x; \lambda)$ denote the map of Lagrangian subspaces
associated with the frames $\mathbf{X}_1 (x; \lambda)$ specified 
in (\ref{frame1}). Keeping in mind that $\lambda_2$ is fixed, let 
$\ell_2 (x; \lambda_2)$ denote the map of Lagrangian subspaces
associated with the frames $\mathbf{X}_2 (x; \lambda_2)$ specified 
in (\ref{frame2}). We emphasize that $\mathbf{X}_2 (x; \lambda_2)$ is 
initialized at $x=1$. Effectively, this 
means that we are looking sideways at the usual Maslov Box, 
setting the target as the right shelf $\lambda = \lambda_2$,
rather than the top shelf. 

By Maslov Box in this setting, we mean the following sequence of contours:
(1) fix $x = 0$ and let $\lambda$ increase from $\lambda_1$ to $\lambda_2$ 
(the {\it bottom shelf}); 
(2) fix $\lambda = \lambda_2$ and let $x$ increase from $0$ to $1$ 
(the {\it right shelf}); (3) fix $x = 1$ and let $\lambda$
decrease from $\lambda_2$ to $\lambda_1$ (the {\it top shelf}); and (4) fix
$\lambda = \lambda_1$ and let $x$ decrease from $1$ to $0$ (the {\it left shelf}). 
(See Figure \ref{box-figure}.)

\begin{figure}[ht]
\begin{center}
\begin{tikzpicture}
\draw[<->] (-8,-1) -- (1,-1);	
\draw[<->] (0,-2) -- (0,4.5);	
\node at (.5,4.3) {$x$};
\node at (-7.5,-1.5) {$\lambda$};
\node at (-6,-1.5) {$\lambda_1$};
\node at (-1,-1.5) {$\lambda_2$};
\node at (.5,3) {$1$};
\draw[-] (-.1,3) -- (.1,3);
\node at (.5,-.7) {$0$};
\draw[-] (-.1,-1) -- (.1,-1);
%
\draw[thick, ->] (-6,-1) -- (-3.5,-1);
\draw[thick] (-3.5,-1) -- (-1,-1);	
\draw[thick, ->] (-1,-1) -- (-1,1);
\draw[thick] (-1,1) -- (-1, 3);
\draw[thick,->] (-1,3) -- (-3.5, 3);
\draw[thick] (-3.5,3) -- (-6, 3);
\draw[thick,->] (-6,3) -- (-6, 1);
\draw[thick] (-6,1) -- (-6, -1);
%
\node[scale = .75] at (-3.5, -.55) {$\mas (\ell_1 (0; \cdot), \ell_2 (0; \lambda_2); [\lambda_1, \lambda_2])$};
\node[scale = .75, rotate=90] at (-.55, 1.1) {$\mas (\ell_1 (\cdot; \lambda_2), \ell_2 (\cdot; \lambda_2); [0, 1])$};
\node[scale = .75] at (-3.6, 3.45) {$- \mas (\ell_1 (1; \cdot), \ell_2 (1; \lambda_2); [\lambda_1, \lambda_2])$};
\node[scale = .75, rotate=90] at (-6.45, 1.1) {$- \mas (\ell_1 (\cdot; \lambda_1), \ell_2 (\cdot; \lambda_2); [0, 1])$};
\end{tikzpicture}
\end{center}
\caption{The Maslov Box.} \label{box-figure}
\end{figure}

Following the general framework discussion in Section \ref{maslov_section},
we will compute the Maslov index along the contours specified
in the Maslov box as a spectral flow of the unitary matrix function
\begin{equation} \label{tilde-W-defined-1}
\begin{aligned}
\tilde{W} (x; \lambda) &: = - (X_1 (x; \lambda) + i Y_1 (x; \lambda)) (X_1 (x; \lambda) - i Y_1 (x; \lambda))^{-1} \\
& \quad \quad \times (X_2 (x; \lambda_2) - i Y_2 (x; \lambda_2)) (X_2 (x; \lambda_2) + i Y_2 (x; \lambda_2))^{-1}.
\end{aligned}
\end{equation}

{\it The bottom shelf.}
We begin our analysis with the bottom shelf. Since 
$\mathbf{X}_1 (0; \lambda) = J \alpha^*$ for all 
$\lambda \in [\lambda_1, \lambda_2]$ (in particular, 
is independent of $\lambda$), and $\mathbf{X}_2 (0; \lambda_2)$ does 
not vary with $\lambda$, we see that in fact the matrix 
$\tilde{W} (0; \lambda)$ is constant as $\lambda$ varies
from $\lambda_1$ to $\lambda_2$, and so 
\begin{equation} \label{bottom_shelf}
\mas (\ell_1 (0; \cdot), \ell_2 (0; \lambda_2); [\lambda_1, \lambda_2]) = 0.
\end{equation} 
This does not necessarily mean that $-1$ is not an eigenvalue of 
$\tilde{W} (0; \lambda)$; rather, if $-1$ is an eigenvalue of
$\tilde{W} (0;\lambda)$ with multiplicity $m$ for some 
$\lambda \in [\lambda_1, \lambda_2]$, then it remains fixed 
as an eigenvalue of $\tilde{W} (0; \lambda)$ with multiplicity 
$m$ for all $\lambda \in [\lambda_1, \lambda_2]$.  

{\it The right shelf.}
For the right shelf, $\lambda$ is fixed at $\lambda_2$ for both 
$\mathbf{X}_1$ and $\mathbf{X}_2$. By construction, 
$\ell_1 (\cdot; \lambda_2)$ will intersect $\ell_2 (\cdot; \lambda_2)$
at some $x = x_*$ with dimension $m$ if and only if $\lambda_2$ is an eigenvalue 
of (\ref{hammy}) with multiplicity $m$. In the event that $\lambda_2$
is not an eigenvalue of (\ref{hammy}), there will be no 
crossing points along the right shelf. On the other hand, if
$\lambda_2$ is an eigenvalue of (\ref{hammy}) with multiplicity
$m$, then $\tilde{W} (x; \lambda_2)$ will have $-1$ as an eigenvalue
with multiplicity $m$ for all $x \in [0, 1]$. In either case,
\begin{equation} \label{right_shelf}
\mas (\ell_1 (\cdot; \lambda_2), \ell_2 (\cdot; \lambda_2); [0, 1]) = 0.
\end{equation} 

{\it The top shelf.}
For the top shelf, we know from Lemma \ref{monotonicity-lemma1} that monotonicity
in $\lambda$ is determined by 
$- \mathbf{X}_1 (1; \lambda)^* J \partial_{\lambda} \mathbf{X}_1 (1; \lambda)$, 
and we readily compute 
\begin{equation*}
\begin{aligned}
\frac{\partial}{\partial x} \mathbf{X}_1^* (x; \lambda) J \partial_{\lambda} \mathbf{X}_1 (x; \lambda)
&= (\mathbf{X}_1^{\prime})^* J \partial_{\lambda} \mathbf{X}_1 
+ \mathbf{X}_1^* J \partial_{\lambda} \mathbf{X}_1^{\prime} \\
&= - (\mathbf{X}_1^{\prime})^* J^t \partial_{\lambda} \mathbf{X}_1 
+ \mathbf{X}_1^* \partial_{\lambda} J \mathbf{X}_1^{\prime} \\
&= - \mathbf{X}_1^* \mathbb{B} (x; \lambda) \partial_{\lambda} \mathbf{X}_1
+ \mathbf{X}_1^* \partial_{\lambda} (\mathbb{B} (x; \lambda) \mathbf{X}_1) 
= \mathbf{X}_1^* \mathbb{B}_{\lambda} \mathbf{X}_1.
\end{aligned}
\end{equation*}
Integrating on $[0,1]$, and noting that $\partial_{\lambda} \mathbf{X}_1 (0;\lambda) = 0$,
we see that   
\begin{equation*}
\mathbf{X}_1 (1; \lambda)^* J \partial_{\lambda} \mathbf{X}_1 (1; \lambda)
= \int_0^1 \mathbf{X}_1 (x;\lambda)^* \mathbb{B}_{\lambda} (x; \lambda) \mathbf{X}_1 (x; \lambda) dy. 
\end{equation*}
In this way, we see that condition {\bf (B1)} ensures that 
as $\lambda$ increases the eigenvalues of $\tilde{W} (1; \lambda)$ will 
rotate in the clockwise direction. Since each crossing along the top shelf
corresponds with an eigenvalue, we can conclude that 
\begin{equation*}
\mathcal{N} ([\lambda_1, \lambda_2)) = 
- \mas (\ell_1 (1; \cdot), \ell_2 (1; \lambda_2); [\lambda_1, \lambda_2]).
\end{equation*}
We note that $\lambda_1$ is included in the count, because in the event
that $(1, \lambda_1)$ is a crossing point, eigenvalues of 
$\tilde{W} (1; \lambda)$ will rotate away from $-1$ in the clockwise 
direction as $\lambda$ increases from $\lambda_1$ (thus decrementing the 
Maslov index). Likewise, $\lambda_2$ is not included in the count, because in the event
that $(1, \lambda_2)$ is a crossing point, eigenvalues of 
$\tilde{W} (1; \lambda)$ will rotate into $-1$ in the clockwise 
direction as $\lambda$ increases to $\lambda_2$ (thus leaving the 
Maslov index unchanged). 

{\it The left shelf.}
Our analysis so far leaves only the left shelf to consider, and 
we observe that the Maslov index along the left shelf can be expressed as 
\begin{equation*}
- \mas (\ell_1 (\cdot; \lambda_1), \ell_2 (\cdot; \lambda_2); [0, 1]).
\end{equation*} 
Using path additivity and homotopy invariance, we can sum the Maslov
indices on each shelf of the Maslov Box to arrive at the relation 
\begin{equation} \label{box_sum}
\mathcal{N} ([\lambda_1, \lambda_2)) 
= \mas (\ell_1 (\cdot; \lambda_1), \ell_2 (\cdot; \lambda_2); [0, 1]),
\end{equation}
which is (\ref{thm1-eqn1}) from Theorem \ref{bc1_theorem}.

In order to get from (\ref{box_sum}) to the second claim in 
Theorem \ref{bc1_theorem}, we need to verify 
that crossings along the left shelf occur monotonically in the
counterclockwise direction as $x$ increases. In this case we will have
\begin{equation*}
\begin{aligned}
\mas (\ell_1 (\cdot; \lambda_1), \ell_2 (\cdot; \lambda_2); [0, 1])
& = \sum_{0 < x \le 1} \dim (\ell_1 (x; \lambda_1) \cap \ell_2 (x; \lambda_2)) \\
&= \sum_{0 < x \le 1} \dim \ker (\mathbf{X}_1 (x; \lambda_1)^* J \mathbf{X}_2 (x; \lambda_2)).
\end{aligned}
\end{equation*} 
Here, $x = 0$ is not included in the sum, because if $x = 0$ is a crossing
point, then as $x$ increases from 0, the eigenvalues of $\tilde{W} (x; \lambda_1)$
will rotate away from $-1$ in the counterclockwise direction, and so will 
not increment the Maslov index.  

In order to address this question of monotonicity, we adapt an 
idea from Remark 4.2 in \cite{Elyseeva2021}. For this, we 
begin by introducing a $2n \times 2n$ fundamental 
matrix $\Psi (x; \lambda)$ for (\ref{hammy}), specified by 
\begin{equation} \label{Psi-specified}
\begin{aligned}
J \Psi' &= \mathbb{B} (x; \lambda_2) \Psi \\
\Psi (1; \lambda_2) &= I_{2n}. 
\end{aligned}
\end{equation}
With this notation, we can express $\mathbf{X}_2 (x; \lambda_2)$
as 
\begin{equation*}
    \mathbf{X}_2 (x; \lambda_2) = \Psi (x; \lambda_2) J\beta^*.
\end{equation*}
We now set 
\begin{equation} \label{bold-X}
    \mathbf{X} (x; \lambda, \lambda_2) := \Psi (x; \lambda_2)^{-1} \mathbf{X}_1 (x; \lambda),
\end{equation}
and note that
\begin{equation*}
    J \mathbf{X}' (x; \lambda, \lambda_2) 
    = J\{ (\Psi (x; \lambda_2)^{-1})' \mathbf{X}_1 (x; \lambda) 
    +  \Psi (x; \lambda_2)^{-1} \mathbf{X}_1' (x; \lambda) \}. 
\end{equation*}
By a straightforward calculation, we can check that 
\begin{equation} \label{useful2}
    \Psi (x; \lambda)^* J \Psi (x; \lambda) = J,
    \quad \, \forall \, x \in [0, 1],
\end{equation}
from which we see that 
\begin{equation*}
    \Psi (x; \lambda_2)^{-1} = - J  \Psi (x; \lambda_2)^* J.
\end{equation*}
This allows us to write 
\begin{equation*}
\begin{aligned}
    J \mathbf{X}' (x; \lambda, \lambda_2) 
    &= J\{ - J  \Psi' (x; \lambda_2)^* J \mathbf{X}_1 (x; \lambda) 
    - J  \Psi (x; \lambda_2)^* J \mathbf{X}_1' (x; \lambda) \} \\
    &= J\{J  (J \Psi' (x; \lambda_2))^* \mathbf{X}_1 (x; \lambda) 
    - J  \Psi (x; \lambda_2)^* \mathbb{B} (x; \lambda) \mathbf{X}_1 (x; \lambda) \}.
\end{aligned}
\end{equation*}
Next, we use (\ref{bold-X}) to re-write this relation as
\begin{equation*}
\begin{aligned}
    J \mathbf{X}' (x; \lambda, \lambda_2) 
    &= \{- (J \Psi' (x; \lambda_2))^* \Psi (x; \lambda_2) \mathbf{X} (x; \lambda, \lambda_2)
    + \Psi (x; \lambda_2)^* \mathbb{B} (x; \lambda) \Psi (x; \lambda_2) \mathbf{X} (x; \lambda, \lambda_2)\}  \\
    &= \{- (\mathbb{B} (x; \lambda_2) \Psi (x; \lambda_2))^* \Psi (x; \lambda_2) \mathbf{X} (x; \lambda, \lambda_2)
    + \Psi (x; \lambda_2)^* \mathbb{B} (x; \lambda) \Psi (x; \lambda_2) \mathbf{X} (x; \lambda, \lambda_2)\}  \\
    &= \Psi (x; \lambda_2)^* \{\mathbb{B} (x; \lambda) - \mathbb{B} (x; \lambda_2) \} \Psi (x; \lambda_2) \mathbf{X} (x; \lambda, \lambda_2). 
\end{aligned}
\end{equation*}
In summary, we can write 
\begin{equation} \label{script-B-defined}
    J \mathbf{X}' = \mathpzc{B} (x; \lambda, \lambda_2) \mathbf{X}, 
    \quad \mathpzc{B} (x; \lambda, \lambda_2) := \Psi (x; \lambda_2)^* \{\mathbb{B} (x; \lambda) - \mathbb{B} (x; \lambda_2) \} \Psi (x; \lambda_2),
\end{equation}
initialized with 
\begin{equation*}
    \mathbf{X} (1; \lambda, \lambda_2) = \Psi (1; \lambda_2)^{-1} \mathbf{X}_1 (1; \lambda)
    = \mathbf{X}_1 (1; \lambda).
\end{equation*}
Since $\mathpzc{B} (x; \lambda, \lambda_2)$ is self-adjoint for a.e. $x \in (0,1)$ and 
$\mathbf{X}_1 (1; \lambda)$ is the frame for a Lagrangian subspace, we 
can conclude that $\mathbf{X} (x; \lambda, \lambda_2)$ is the frame for a Lagrangian 
subspace for each $(x, \lambda) \in [0, 1] \times [\lambda_1, \lambda_2]$. 
In addition, according to Assumption {\bf (B2)}
$\mathpzc{B} (x; \lambda_1, \lambda_2)$ is a non-positive matrix for a.e. $x \in (0,1)$. 

\begin{lemma} \label{frame-lemma}
For each $x \in [0, 1]$,
\begin{equation*}
    \mathbf{X} (x; \lambda, \lambda_2)^* J (J \beta^*)
    = \mathbf{X}_1 (x; \lambda)^* J \mathbf{X}_2 (x; \lambda_2).
\end{equation*}
\end{lemma}

\begin{proof}
To see this, we simply use (\ref{useful2}) to write 
\begin{equation*}
\begin{aligned}
    \mathbf{X} (x; \lambda, \lambda_2)^* J (J \beta^*)
    &= (\Phi (x; \lambda_2)^{-1} \mathbf{X}_1 (x; \lambda))^* J J \beta^*
    = - \mathbf{X}_1 (x; \lambda)^* (\Phi (x; \lambda_2)^{-1})^* \beta^* \\
    &= - \mathbf{X}_1 (x; \lambda)^* (- J  \Phi (x; \lambda_2)^* J)^* \beta^* 
    = \mathbf{X}_1 (x; \lambda)^* J \Phi (x; \lambda_2) J \beta^* \\
    &= \mathbf{X}_1 (x; \lambda)^* J \mathbf{X}_2 (x; \lambda_2).
\end{aligned}
\end{equation*}
\end{proof}

Lemma \ref{frame-lemma} allows us to detect intersections between 
$\ell_1 (x; \lambda)$ and $\ell_2 (x; \lambda_2)$ by instead 
detecting intersections between $\ell (x; \lambda, \lambda_2) := \colspan (\mathbf{X} (x; \lambda, \lambda_2))$ 
and the fixed target space $\ell_{\beta} := \colspan (J \beta^*)$. 
For these latter intersections, 
the associated rotation matrix is 
(with $\mathbf{X} (x; \lambda, \lambda_2) = {X (x; \lambda, \lambda_2) \choose Y (x; \lambda, \lambda_2)}$)
\begin{equation} \label{tilde-W-fixed-target}
\begin{aligned}
    \tilde{\mathcal{W}} (x; \lambda, \lambda_2) 
    &:= - (X (x; \lambda, \lambda_2) + i Y (x; \lambda, \lambda_2)) (X (x; \lambda, \lambda_2) - i Y (x; \lambda, \lambda_2))^{-1} \\
    & \quad \quad \quad \quad \times 
    (-\beta^* - i \beta_1^*) (-\beta_2^* + i \beta_1^*)^{-1}.
\end{aligned}
\end{equation}

At this point, we can proceed with a Maslov-box argument for 
intersections between $\ell (x; \lambda, \lambda_2)$
and $\ell_{\beta}$. As before, we get no contributions from the bottom and 
right shelves, so all that's left to verify is monotonicity on 
the top and left shelves. 

For the top shelf, the rotation is determined as usual by 
\begin{equation*}
\begin{aligned}
    \mathbf{X} (1; \lambda, \lambda_2)^* J \partial_{\lambda} \mathbf{X} (1; \lambda, \lambda_2)
    &= \int_0^1 \mathbf{X} (x; \lambda, \lambda_2)^* \mathpzc{B}_{\lambda} (x; \lambda, \lambda_2) \mathbf{X} (x; \lambda, \lambda_2) d x \\
    &= \int_0^1 \mathbf{X} (x; \lambda, \lambda_2)^* \Psi (x; \lambda_2)^* \mathbb{B}_{\lambda} (x; \lambda) \Psi (x; \lambda_2) \mathbf{X} (x; \lambda, \lambda_2) d x \\
    &= \int_0^1 \mathbf{X}_1 (x; \lambda)^* \mathbb{B}_{\lambda} (x; \lambda) \mathbf{X}_1 (x; \lambda) d x,
\end{aligned}
\end{equation*}
and so we obtain positivity as previously from {\bf (B1)}. 

For the left shelf, we've seen that $\mathpzc{B} (x; \lambda_1, \lambda_2)$ is non-positive
for a.e. $x \in (0, 1)$, and it follows immediately from Lemma \ref{monotonicity-in-x-lemma}
that no eigenvalue of $\tilde{\mathcal{W}} (x; \lambda_1, \lambda_2)$ can rotate in the clockwise 
direction as $x$ increases along $[0, 1]$. This still leaves open the 
possibility that an eigenvalue of $\tilde{\mathcal{W}} (x; \lambda_1, \lambda_2)$ rotates into 
$-1$ at $x = x_*$ and remains there for some interval $[x_*, x_* + \delta]$,
$\delta > 0$, but this is precisely the event that is ruled out by
the second part of
Assumption {\bf (B2)}. We can conclude that the Maslov index for the left
shelf is a monotonic (positive) count of intersections between $\ell (x; \lambda_1, \lambda_2)$
and $\ell_{\beta}$, and this can be expressed as 
\begin{equation} \label{box_sum2}
\begin{aligned}
\mas (\ell (\cdot; \lambda_1, \lambda_2), \ell_{\beta}; [0, 1]) 
&= \sum_{0 < x \le 1} \dim (\ell (x; \lambda_1, \lambda_2) \cap \ell_{\beta}) \\
&= \sum_{0 < x \le 1} \dim (\ell_1 (x; \lambda_1) \cap \ell_2 (x; \lambda_2)) \\
&= \sum_{0 < x \le 1} \dim \ker (\mathbf{X}_1 (x; \lambda_1)^* J \mathbf{X}_2 (x; \lambda_2)).
\end{aligned}
\end{equation}

Using again path additivity and 
homotopy invariance, we can sum the Maslov
indices on each shelf of the Maslov Box to arrive at the relation 
\begin{equation} \label{box_sum3}
\begin{aligned}
\mathcal{N} ([\lambda_1, \lambda_2)) 
&= \mas (\ell (\cdot; \lambda_1, \lambda_2), \ell_{\beta}; [0, 1]) \\
&= \sum_{0 < x \le 1} \dim (\ell_1 (x; \lambda_1) \cap \ell_2 (x; \lambda_2)) \\
&= \sum_{0 < x \le 1} \dim \ker (\mathbf{X}_1 (x; \lambda_1)^* J \mathbf{X}_2 (x; \lambda_2)),
\end{aligned}
\end{equation}
which is (\ref{thm1-eqn2}) in Theorem \ref{bc1_theorem}.

\begin{remark}
    Comparing the relations 
    \begin{equation*}
    \mathcal{N} ([\lambda_1, \lambda_2)) 
    = \mas (\ell_1 (\cdot; \lambda_1), \ell_2 (\cdot; \lambda_2); [0, 1]),    
    \end{equation*}
    and 
    \begin{equation*}
    \mathcal{N} ([\lambda_1, \lambda_2)) 
    = \mas (\ell (\cdot; \lambda_1, \lambda_2), \ell_{\beta}; [0, 1]),    
    \end{equation*}
    we see that 
    \begin{equation} \label{symplectic-invariance}
    \mas (\ell_1 (\cdot; \lambda_1), \ell_2 (\cdot; \lambda_2); [0, 1])
    = \mas (\ell (\cdot; \lambda_1, \lambda_2), \ell_{\beta}; [0, 1]).
    \end{equation}
    Since these two Maslov indices detect precisely the same 
    intersections with the same multiplicities, and the index on 
    the right-hand side is a monotonic count of crossings, it 
    must be the case that the index on the left-hand side is 
    also a monotonic count of crossings. It's interesting to 
    note, however, that while the eigenvalues of $\tilde{\mathcal{W}} (x; \lambda_1, \lambda_2)$
    rotate monotonically as $x$ increases along $(0, 1)$, 
    the eigenvalues of $\tilde{W} (x; \lambda_1)$ are only necessarily monotonic at 
    crossing points. In this way, the computation of 
    $\mas (\ell (\cdot; \lambda_1, \lambda_2), \ell_{\beta}; [0, 1])$ is 
    more tractable than the corresponding 
    calculation of $\mas (\ell_1 (\cdot; \lambda_1), \ell_2 (\cdot; \lambda_2); [0, 1])$. 
    Nonetheless, our preference has been to cast our analysis
    primarily in terms of the latter, because it is in no
    way tied to the specific application considered 
    here, and so generalizes more readily to cases such 
    as singular Hamiltonian systems.  
    
    The authors are grateful to the referee for pointing out that 
    (\ref{symplectic-invariance}) can also be viewed as a consequence of symplectic 
    invariance (Property V in Section 1 of \cite{CLM}); namely,
    we have the relations 
    \begin{equation*}
        \mathbf{X} (x; \lambda_1, \lambda_2) = \Psi (x; \lambda_2)^{-1} \mathbf{X}_1 (x; \lambda_1)
        \quad \text{and} \quad
        J\beta^* = \Psi (x; \lambda_2)^{-1} \mathbf{X}_1 (x; \lambda_1), 
    \end{equation*}
    and (\ref{symplectic-invariance}) follows immediately by symplectic 
    invariance and the observation that 
    $\Psi (x; \lambda_2)^{-1}$ is a symplectic matrix for 
    all $x \in [0, 1]$. Moreover, if we specify the symplectic matrix 
    \begin{equation*}
        \mathfrak{B} 
        := \begin{pmatrix}
        \beta \\ - (\beta \beta^*) \beta J
        \end{pmatrix}, 
    \end{equation*}
    we see that $\mathfrak{B} J \beta^* = (0 \,\,\, I)^T$, so that 
    symplectic invariance additionally implies 
    \begin{equation} \label{kratz-setting}
    \mas (\ell (\cdot; \lambda_1, \lambda_2), \ell_{\beta}; [0, 1])
    = \mas (\ell_{\mathfrak{B}} (\cdot; \lambda_1, \lambda_2), \ell_D; [0, 1]),
    \end{equation}
    where $\ell_{\mathfrak{B}} (\cdot; \lambda_1, \lambda_2)$ is the Lagrangian subspace
    with frame $\mathfrak{B} \mathbf{X} (x; \lambda_1, \lambda_2)$, and $\ell_D$ denotes
    the usual Dirichlet subspace. Restricted to the setting of $\mathbb{R}^{2n}$, 
    the right-hand side of (\ref{kratz-setting}) now fits into the  
    framework of \cite{Kratz1995}, with 
    our crossing points corresponding precisely with the {\bf focal points}
    of Definition 1.1.1(ii) in that reference. In addition, our Assumption 
    {\bf (B2)} implies the {\bf controllability} property specified in 
    Definition 4.1.1 of \cite{Kratz1995} (see also Theorem 4.1.3 of \cite{Kratz1995}).
\end{remark}

We conclude this section by establishing a convenient criterion for 
verifying Assumption {\bf (B2)}. 

\begin{lemma} \label{B2verification}
Let Assumptions {\bf (A)} hold, and also assume there
exist values $\lambda_1, \lambda_2 \in I$, $\lambda_1 < \lambda_2$, 
so that $\mathbb{B} (x; \lambda_2) - \mathbb{B} (x; \lambda_1)$ is non-negative
for a.e. $x \in (0, 1)$.  Suppose that for any $[a, b] \subset [0, 1]$,
$a < b$, and any non-trivial 
solution $y(x; \lambda_1)$ of $J y' = \mathbb{B} (x; \lambda_1) y$, we must have 
\begin{equation} \label{difference-definite-again}
\int_a^b ((\mathbb{B} (x; \lambda_2) - \mathbb{B} (x; \lambda_1)) y(x; \lambda_1), y(x; \lambda_1)) dx > 0. 
\end{equation}
We have the following:

\smallskip
\noindent
(i) Let Assumption {\bf (B1)} hold, and let 
$\mathbf{X}_1 (x; \lambda_1)$ and $\mathbf{X}_2 (x; \lambda_2)$
be the Lagrangian frames specified in (\ref{frame1}) and 
(\ref{frame2}), with respective Lagrangian subspaces $\ell_1 (x; \lambda_1)$
and $\ell_2 (x; \lambda_2)$. Then there is no interval $[a, b] \subset [0, 1]$, $a < b$, 
so that 
\begin{equation*}
\dim (\ell_1 (x; \lambda_1) \cap \ell_2 (x; \lambda_2)) \ne \{0\}
\end{equation*}
for all $x \in [a, b]$. 

\smallskip
\noindent
(ii) Let Assumption {\bf (B1)$\mathbf{^{\prime}}$} hold, and let 
$\mathbf{X}_3 (x; \lambda_1)$ and $\mathbf{X}_4 (x; \lambda_2)$
be the Lagrangian frames specified in (\ref{frame3}) and 
(\ref{frame4}), with respective Lagrangian subspaces $\ell_3 (x; \lambda_1)$
and $\ell_4 (x; \lambda_2)$. Then there is no interval $[a, b] \subset [0, 1]$, $a < b$,
so that 
\begin{equation*}
\dim (\ell_3 (x; \lambda_1) \cap \ell_4 (x; \lambda_2)) \ne \{0\}
\end{equation*}
for all $x \in [a, b]$.
\end{lemma}

\begin{proof} 
We proceed by contradiction. First, for (i),
suppose there exists an interval $[a,b] \subset [0,1]$, $a<b$,
so that 
\begin{equation*}
\ell_1 (x;\lambda_1) \cap \ell_2 (x; \lambda_2) \ne \{0\}
\end{equation*} 
for all $x \in [a, b]$. Then from Lemma \ref{frame-lemma},
we have 
\begin{equation*}
\ell (x;\lambda_1, \lambda_2) \cap \ell_{\beta} \ne \{0\}
\end{equation*} 
for all $x \in [a, b]$. For this latter relation, we've 
seen that any eigenvalue of $\tilde{\mathcal{W}} (x; \lambda_1, \lambda_2)$ 
that crosses $-1$ as $x$ increases through a crossing point $x_*$ 
must cross in the counterclockwise direction. This allows us to apply 
Lemma \ref{continuation-lemma1} to conclude that there exists
some interval $[c, d] \subset [a, b]$, 
$c < d$, along with some integer $m \in \{1, 2, \dots, n\}$,
so that $\dim \ker (\tilde{\mathcal{W}} (x; \lambda_1, \lambda_2) + I) = m$ for 
all $x \in [c, d]$. According, then, to 
Lemma \ref{continuation-lemma2}, we can conclude that
there exists $v_0 \in \mathbb{C}^n \backslash \{0\}$ so that 
$\mathbf{X} (x; \lambda_1, \lambda_2) v_0 \in \ell (x; \lambda_1, \lambda_2) \cap \tilde{\ell}$ for all 
$x$ in some possibly smaller interval 
$[\tilde{a}, \tilde{b}] \subset [c, d]$, $\tilde{a} < \tilde{b}$,
and subsequently we can conclude from Lemma \ref{continuation-lemma3}
that $\mathbf{X} (x; \lambda_1, \lambda_2) v_0$ is independent of $x$.
It follows from the relation $J \mathbf{X}' = \mathpzc{B} (x; \lambda_1, \lambda_2) \mathbf{X}$
that $\mathpzc{B} (x; \lambda_1, \lambda_2) \mathbf{X} v_0 = 0$ for a.e.  
$x \in [\tilde{a}, \tilde{b}]$, and consequently we see that 
\begin{equation} \label{integral1}
    \int_{\tilde{a}}^{\tilde{b}} (\mathpzc{B} (x; \lambda_1, \lambda_2) 
    \mathbf{X} (x; \lambda_1, \lambda_2) v_0, \mathbf{X} (x; \lambda_1, \lambda_2) v_0) dx = 0.
\end{equation}
Recalling (\ref{script-B-defined}), we can express (\ref{integral1}) as 
\begin{equation*}
\begin{aligned}
     \int_{\tilde{a}}^{\tilde{b}}&
     (\Psi (x; \lambda_2)^* \{\mathbb{B} (x; \lambda_1) - \mathbb{B} (x; \lambda_2) \} \Psi (x; \lambda_2) 
     \mathbf{X} (x; \lambda_1, \lambda_2) v_0, \mathbf{X} (x; \lambda_1, \lambda_2) v_0) dx \\
     &=  \int_{\tilde{a}}^{\tilde{b}}
    (\{\mathbb{B} (x; \lambda_1) - \mathbb{B} (x; \lambda_2) \} \Psi (x; \lambda_2) \mathbf{X} (x; \lambda_1, \lambda_2) v_0,
     \Psi (x; \lambda_2) \mathbf{X} (x; \lambda_1, \lambda_2) v_0) dx \\
     &= \int_{\tilde{a}}^{\tilde{b}}
    (\{\mathbb{B} (x; \lambda_1) - \mathbb{B} (x; \lambda_2) \} \mathbf{X}_1 (x; \lambda_1) v_0,
    \mathbf{X}_1 (x; \lambda_1) v_0) dx = 0. 
\end{aligned}     
\end{equation*}
But $\mathbf{X}_1 (x; \lambda_1) v_0$ is a non-trivial solution to $Jy' = \mathbb{B} (x; \lambda_1) y$,
so this contradicts (\ref{difference-definite-again}), giving Item (i). 

Turning to Item (ii), we now suppose 
there exists an interval $[a,b] \subset [0,1]$, $a<b$,
so that 
\begin{equation} \label{contradict1}
\ell_3 (x;\lambda_1) \cap \ell_4 (x; \lambda_2) \ne \{0\}
\end{equation} 
for all $x \in [a, b]$. Arguing similarly as for Item (i), we find 
that there exists $v_0 \in \mathbb{C}^{2n} \backslash \{0\}$
so that $z (x; \lambda_1) := \mathbf{X}_3 (x; \lambda_1) v_0$
satisfies the relation 
\begin{equation*}
\int_{\tilde{a}}^{\tilde{b}}
(\{\mathcal{B} (x; \lambda_1) - \mathcal{B} (x; \lambda_2) \} z (x; \lambda_1),
z (x; \lambda_1)) dx = 0, 
\end{equation*}
for some interval $[\tilde{a}, \tilde{b}] \subset [0, 1]$, not necessarily the same 
as for Item (i). (Here, we recall that $\mathcal{B} (x; \lambda)$ denotes the matrix
specified in (\ref{mathcal-B-defined}) associated with (\ref{hammy})--{\bf (BC2)}.) 
If we write $z = (z_1, z_2, z_3, z_4)^T$, and set $w := (z_2, z_4)^T$ 
then we can express this last integral relation as 
\begin{equation*} 
\int_{\tilde{a}}^{\tilde{b}} ((\mathbb{B} (x; \lambda_1) - \mathbb{B} (x; \lambda_2)) w(x; \lambda_1), w(x; \lambda_1)) dx = 0, 
\end{equation*}
which contradicts (\ref{difference-definite-again}) unless $z_2 (x; \lambda_1)$ and 
$z_4 (x; \lambda_1)$ are identically 0 on $[\tilde{a}, \tilde{b}]$. In this case, 
it's clear from the specification of $\mathbf{X}_3 (x; \lambda_1)$ 
in (\ref{frame3}) that there must exist constant vectors $c_1, c_3 \in \mathbb{C}^n$ 
so that $z (x; \lambda_1) = (c_1, 0, c_3, 0)^T$ for all $x \in [\tilde{a}, \tilde{b}]$. 
But $z (x; \lambda_1) \in \ell_3 (x; \lambda_1)$ as specified in (\ref{lagrangian2}),
from which we see that $z (x; \lambda_1)$ must have the form 
$z (x; \lambda_1) = (y_1 (0; \lambda_1), y_1 (x; \lambda_1), - y_2 (0; \lambda_1), y_2 (x; \lambda_1))$
for some $y (x; \lambda_1)$ satisfying $J y' = \mathbb{B} (x; \lambda_1) y$. I.e., we
must have $y_1 (x; \lambda_1) = 0$ and $y_2 (x; \lambda_1) = 0$ for all $x \in [\tilde{a}, \tilde{b}]$,
and this can only happen if $y_1 (0; \lambda_1) = 0$ and $y_2 (0; \lambda_1) = 0$, in 
which case $z (x; \lambda_1) = 0$ for all $x \in [\tilde{a}, \tilde{b}]$, 
contradicting the implication of (\ref{contradict1}). 
\end{proof}

\subsection{Proof of Theorem \ref{bc2_theorem}}

Given the framework developed in the Introduction, 
the proof of Theorem \ref{bc2_theorem} is almost identical 
to the proof of Theorem \ref{bc1_theorem}, and we omit the 
details.

\section{Applications} \label{applications}

In this section, we verify that our assumptions are satisfied 
by five example cases: Dirac systems, Sturm-Liouville systems, 
the family of linear Hamiltonian systems considered in \cite{GZ2017},
a system associated with differential-algebraic 
Sturm-Liouville systems, and a fourth-order self-adjoint equation. 
In the fourth case, $\mathbb{B} (x; \lambda)$ is nonlinear in $\lambda$, 
and in the fifth we will consider periodic boundary conditions. In 
all cases, we also demonstrate our theory with a numerical calculation, 
though our numerical calculation for Sturm-Liouville systems serves
also as the numerical calculation for the family of systems 
from \cite{GZ2017} (which includes Sturm-Liouville systems). 
The particular coefficient functions and 
boundary conditions chosen for these numerical calculations are
not taken from physical applications, but rather have been selected
to correspond with illustrative spectral curves (see Figures 
\ref{dirac_figure}, \ref{sl_figure}, \ref{da_figure}, and 
\ref{fourth-figure}). 

\subsection{Dirac Systems} \label{dirac}

The canonical systems in which our assumptions clearly 
hold are Dirac systems, by which we mean equations of the general form 
\begin{equation} \label{dirac-equation}
J y' = (\lambda Q(x) + V(x))y, 
\end{equation}
where $Q, V \in L^1 ((0,1), \mathbb{C}^{2n})$ are self-adjoint
matrices and $Q (x)$ is positive definite for a.e. $x \in (0,1)$. 
For this example, we will assume separated boundary conditions {\bf (BC1)}.

We can think of this system in terms of the operator
\begin{equation*}
\mathcal{L}_d := Q(x)^{-1} (J \frac{d}{dx} - V(x)),
\end{equation*}
with which we associate the domain 
\begin{equation*}
\mathcal{D} (\mathcal{L}_d) = \{y \in L^2 ((0,1), \mathbb{C}^{2n}): y \in AC ([0,1], \mathbb{C}^{2n}),
\mathcal{L}_d y \in L^2 ((0,1), \mathbb{C}^{2n}); \text{{\bf (BC1)} holds} \},
\end{equation*}
and the inner product
\begin{equation*}
\langle f, g \rangle_Q := \int_0^1 (Q f, g) dx.
\end{equation*}
With this choice of domain and inner product, $\mathcal{L}_d$ is 
densely defined, closed, and self-adjoint, (see, e.g., \cite{Krall2002}), 
so $\sigma (\mathcal{L}_d) \subset \mathbb{R}$. I.e., we can take the 
interval $I$ associated with (\ref{hammy}) to be $I = \mathbb{R}$. 

In this case, $\mathbb{B} (x; \lambda) = \lambda Q(x) + V(x)$,
and we see immediately that our Assumptions {\bf (A)} hold.
For {\bf (B1)}, $\mathbb{B}_{\lambda} (x; \lambda) = Q(x)$, so 
that 
\begin{equation*} 
\int_0^1 \mathbf{X}_1 (x; \lambda)^* \mathbb{B}_{\lambda} (x; \lambda) \mathbf{X}_1 (x; \lambda) dx
=
\int_0^1 \mathbf{X}_1 (x; \lambda)^* Q(x) \mathbf{X}_1 (x; \lambda) d x, 
\end{equation*}
which is positive definite (since $Q (x)$ is positive
definite for a.e. $x \in (0, 1)$).  

For {\bf (B2)}, given any $\lambda_1, \lambda_2 \in \mathbb{R}$, 
$\lambda_1 < \lambda_2$, we have 
$\mathbb{B} (x; \lambda_2) - \mathbb{B} (x; \lambda_1) = (\lambda_2 - \lambda_1) Q (x)$,
which is clearly non-negative a.e. (in fact, positive definite).  
For the moreover part, we notice that 
\begin{equation*}
\mathbb{B} (x; \lambda_2) - \mathbb{B} (x; \lambda_1) 
= (\lambda_2 - \lambda_1) \mathbb{B}_{\lambda} (x; \lambda),
\end{equation*}
which allows us to use the argument establishing {\bf (B1)}
above to show that the condition in Claim \ref{B2verification} is satisfied.

We conclude from Theorem \ref{bc1_theorem} that if $\mathcal{N} ([\lambda_1, \lambda_2); \mathcal{L}_d)$
denotes the spectral count for $\mathcal{L}_d$, we have 
\begin{equation*}
\mathcal{N} ([\lambda_1, \lambda_2); \mathcal{L}_d) 
= \mathcal{N}_{(0,1]} (\mathbf{X}_1 (\cdot; \lambda_1)^* J \mathbf{X}_2 (\cdot; \lambda_2)), 
\end{equation*}
where $\mathbf{X}_1 (x; \lambda_1)$ and $\mathbf{X}_2 (x; \lambda_2)$ denote
the frames specified respectively in (\ref{frame1}) and (\ref{frame2}), 
with $\mathbb{B} (x; \lambda)$
as in this section, and the notation $\mathcal{N}_{(0,1]} (\cdot)$ is 
as in (\ref{counting-function}). 

In order to illustrate the difference between the approach taken in 
\cite{HJK2018, HS2016} and the renormalized approach taken here, we consider a 
specific example with $Q = I_{4}$,
\begin{equation} \label{example1}
V (x) := 
\begin{pmatrix}
.13 + \frac{.7*\cos(6\pi x)}{2+\cos(6\pi x)} & \frac{\cos(\pi x)}{2+\cos(4\pi x)} & 0 & 0 \\
\frac{\cos(\pi x)}{2+\cos(4\pi x)} & 1 & 0 & 0 \\
0 & 0 & 0 & 0 \\
0 & 0 & 0 & 0
\end{pmatrix},
\end{equation}
and Neumann boundary conditions specified by 
\begin{equation} \label{eg1bc}
\alpha = 
\begin{pmatrix}
0_2 & I_2
\end{pmatrix} \quad
\beta = 
\begin{pmatrix}
0_2 & I_2
\end{pmatrix}. 
\end{equation}

\begin{figure}[ht] 
\begin{center}\includegraphics[%
  width=6.5cm,
  height=5cm]{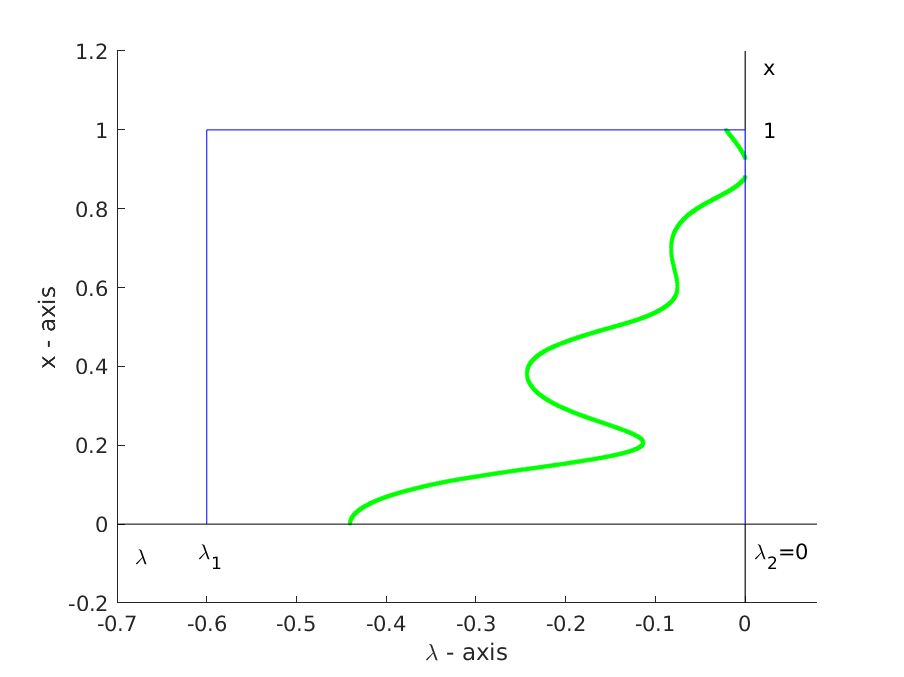}
\includegraphics[%
  width=6.5cm,
  height=5cm]{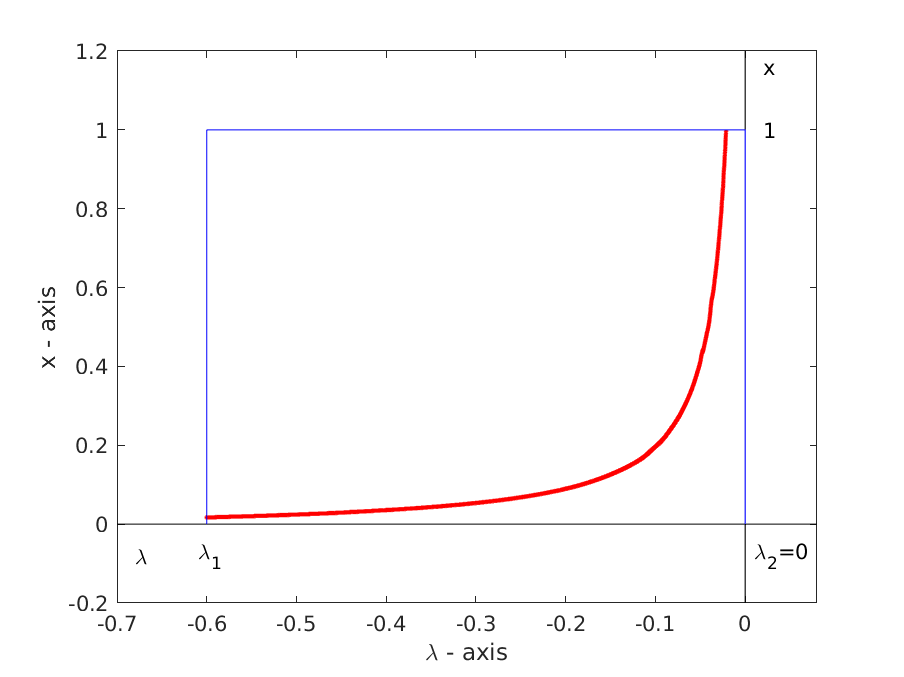}
\end{center}
\caption{Spectral curves for the Dirac equation example: approach of \cite{HJK2018, HS2016} on left; 
renormalized approach on right. \label{dirac_figure}}
\end{figure}

As noted in Remark \ref{standard-v-renormalized}, the authors 
of \cite{HJK2018, HS2016} specify $\mathbf{X}_1 (x; \lambda)$ precisely
as here, but in lieu
of $\mathbf{X}_2 (x; \lambda_2)$, use the fixed target space 
$\tilde{\mathbf{X}}_2 = J \beta^*$. In particular, in the setting
of \cite{HJK2018, HS2016} the Maslov index is computed via the unitary 
matrix 
\begin{equation*}
\tilde{\mathcal{W}}_{\beta} (x; \lambda)
:= - (X_1 (x; \lambda) + i Y_1 (x; \lambda)) (X_1 (x; \lambda) - i Y_1 (x; \lambda))^{-1}
(-\beta_2^* - i \beta_1^*) (-\beta_2^* + i \beta_1^*)^{-1}. 
\end{equation*}
The spectral curves discussed in Remark \ref{standard-v-renormalized}
can be computed in the case of \cite{HJK2018, HS2016} as the pairs $(x, \lambda)$ for which 
$\dim \ker (\tilde{\mathcal{W}}_{\beta} (x; \lambda) + I) \ne 0$,
and likewise can be computed in the current setting as the 
pairs $(x, \lambda)$ for which 
$\dim \ker (\tilde{W} (x; \lambda) + I) \ne 0$. Spectral curves 
for (\ref{dirac-equation}) with $Q = I_4$, $V$ specified in 
(\ref{example1}), and boundary conditions (\ref{eg1bc}) are 
depicted in Figure \ref{dirac_figure}, with the approach of 
\cite{HJK2018, HS2016} on the left and the renormalized approach 
on the right. Several things are worth noting about this 
comparison of images: (1) the difference between the 
non-monotonic curve on the left and the monotonic curve 
on the right is striking and illustrates precisely the main 
difference in the two approaches; (2) while the spectral curve 
on the left emerges from the bottom shelf, the
spectral curve on the right enters the Maslov box
through the left shelf; and (3) since crossings along the 
top shelf correspond with eigenvalues in both cases, 
the spectral curves in the left and right sides of 
Figure \ref{dirac_figure} both cross the top shelf at 
the same value of $\lambda$. Regarding Item (2), 
in the setting of \cite{HJK2018, HS2016}, spectral curves can enter 
through any of the three shelves---left, bottom, or right---while 
in the renormalized setting, spectral curves can only enter through 
the left shelf. 

\subsection{Linear Hamiltonian Systems with Block Matrix Coefficients} 
\label{block-matrix}

In \cite{GZ2017}, the authors consider linear Hamiltonian systems
\begin{equation} \label{bm-equation}
J y' = (\lambda Q(x) + V(x)) y,
\end{equation}
where 
\begin{equation*}
Q(x) 
= \begin{pmatrix}
R(x) & 0 \\
0 & 0
\end{pmatrix}, 
\end{equation*}
for some $r \times r$ matrix $R(x)$, $1 \le r \le 2n$. 
The matrices $R(x)$ and $V(x)$ are taken to be self-adjoint
for a.e. $x \in (0, 1)$, with 
$R(x)$ additionally positive definite for a.e. $x \in (0,1)$, and in 
the bounded-interval case, the authors assume
$Q, V \in L^1 ((0,1), \mathbb{C}^{2n \times 2n})$. (In \cite{GZ2017},
the authors work on a general bounded interval $(a, b)$, but this 
can always be scaled for convenience to $(0,1)$.) 

In order to accommodate the form of $Q$, the authors of \cite{GZ2017}
introduce a Hilbert space
\begin{equation*}
L_R^2 ((0, 1), \mathbb{C}^r) := 
\{f: (0, 1) \to \mathbb{C}^{r} \text{ measurable}, \, 
\|f\|_{L_R^2 ((0,1), \mathbb{C}^r)} < \infty \},
\end{equation*} 
where $\| \cdot \|_{L_R^2 ((0,1), \mathbb{C}^r)}$ denotes the weighted norm 
\begin{equation*}
\|f\|_{L^2_R ((0,1), \mathbb{C}^r)}^2 := \int_0^1 (R(x) f(x), f(x)) dx.
\end{equation*}
In addition, denoting the natural restriction operator 
$\hat{E}_r: \mathbb{C}^{2n} \to \mathbb{C}^r$, the authors 
introduce 
\begin{equation*}
L_Q^2 ((0, 1), \mathbb{C}^{2n}) := 
\{g: (0, 1) \to \mathbb{C}^{2n} \text{ measurable}, \, 
\hat{E}_r g \in L^2_R ((0,1), \mathbb{C}^r) \},
\end{equation*}  
along with the seminorm 
\begin{equation*}
\|g\|_{L^2_Q ((0,1), \mathbb{C}^{2n})} 
:= \|\hat{E}_r g\|_{L^2_R ((0,1), \mathbb{C}^r)}.
\end{equation*}
Finally, the authors assume Atkinson's definiteness condition, 
described as follows: assume that for all $a, b \in (0,1)$ with 
$a < b$, any nonzero solution $y \in AC([0,1], \mathbb{C}^{2n})$
of (\ref{bm-equation}) satisfies 
\begin{equation*}
\|\chi_{[a, b]} y \|_{L^2_Q ((0,1), \mathbb{C}^{2n})} > 0,
\end{equation*}
where $\chi_{[a, b]}$ denotes the usual characteristic 
function on $[a, b]$. 

Under these assumptions, the authors of \cite{GZ2017} are 
able to express (\ref{bm-equation}) in terms of the operator
\begin{equation} \label{bm-operator}
\mathcal{L}_{b} := C(x) (J \frac{d}{dx} - V(x)),
\end{equation}
where 
\begin{equation*}
C(x) = 
\begin{pmatrix}
R(x)^{-1} & 0 \\
0 & I_{2n-r}
\end{pmatrix},
\end{equation*}
and the domain of $\mathcal{L}_{b}$ is specified as 
\begin{equation*}
\mathcal{D} (\mathcal{L}_{b}) := 
\{y \in L^2_Q ((0,1), \mathbb{C}^{2n}): y \in AC ([0, 1], \mathbb{C}^{2n}), 
\mathcal{L}_{b} y \in E_r L^2_Q ((0,1), \mathbb{C}^{2n}), 
\text{{\bf (BC1)} holds}\},
\end{equation*}
with
\begin{equation*}
E_r = 
\begin{pmatrix}
I_r & 0 \\
0 & 0
\end{pmatrix}.
\end{equation*}
In \cite{GZ2017} the authors verify in Section 2 that 
$\sigma(\mathcal{L}_{b}) \subset \mathbb{R}$. 

We see directly from these specifications that our assumptions 
{\bf (A)} hold in this case. To check {\bf (B1)}, we compute 
$\mathbb{B}_{\lambda} (x; \lambda) = Q(x)$, from which we see 
that 
\begin{equation*}
\int_0^1 \mathbf{X}_1 (x; \lambda)^* \mathbb{B}_{\lambda} (x; \lambda) \mathbf{X}_1 (x; \lambda) dx 
= \int_0^1 \mathbf{X}_1 (x; \lambda)^* Q (x) \mathbf{X}_1 (x; \lambda) dx.
\end{equation*}
Since $Q (x)$ is non-negative for a.e. $x \in (0,1)$, this integral is certainly non-negative,
and moreover, it can only be zero if there exists a vector $v \in \mathbb{C}^n$
so that $Q(x) \mathbf{X}_1 (x; \lambda) v = 0$ for a.e. $x \in (0, 1)$.
By definition of $\mathbf{X}_1$, $\psi (x) := \mathbf{X}_1 (x; \lambda) v$
solves $J \psi' = \mathbb{B} (x; \lambda) \psi$. If we write $\psi = {\psi_1 \choose \psi_2}$,
with $\psi_1 (x; \lambda) \in \mathbb{C}^r$ and $\psi_2 (x; \lambda) \in \mathbb{C}^{2n - r}$,
we see that since $Q(x) \psi (x; \lambda) = 0$ for a.e. $x \in (0, 1)$, we must have
$R(x) \psi_1 (x; \lambda) = 0$ for a.e. $x \in (0, 1)$, and so $\psi_1 (x; \lambda) = 0$
for a.e. $x \in (0, 1)$. But then 
\begin{equation*}
\|\chi_{[0, 1]} \psi \|_{L^2_Q ((0,1), \mathbb{C}^{2n})} 
= \int_0^1 (R(x) \psi_1 (x; \lambda), \psi_1 (x; \lambda)) dx = 0,
\end{equation*}
and this contradicts Atkinson's positivity assumption. 

For {\bf (B2)}, we fix $\lambda_1, \lambda_2 \in I$, $\lambda_1 < \lambda_2$,
and observe that 
\begin{equation*}
\mathbb{B} (x; \lambda_2) - \mathbb{B} (x; \lambda_1) 
= (\lambda_2 - \lambda_1) \mathbb{B}_{\lambda} (x; \lambda).
\end{equation*}
We see immediately that $\mathbb{B} (x; \lambda_2) - \mathbb{B} (x; \lambda_1)$
is non-negative, and in addition, the same argument used to verify 
{\bf (B1)} shows that the condition assumed in Claim \ref{B2verification}
holds. Assumption {\bf (B2)} follows.

We conclude from Theorem \ref{bc1_theorem} that if 
$\mathcal{N} ([\lambda_1, \lambda_2); \mathcal{L}_b)$
denotes the spectral count for $\mathcal{L}_b$, we have 
\begin{equation*}
\mathcal{N} ([\lambda_1, \lambda_2); \mathcal{L}_b) 
= \mathcal{N}_{(0,1]} (\mathbf{X}_1 (\cdot; \lambda_1)^* J \mathbf{X}_2 (\cdot; \lambda_2)), 
\end{equation*}
where $\mathbf{X}_1 (x; \lambda_1)$ and $\mathbf{X}_2 (x; \lambda_2)$ denote
the frames (\ref{frame1}) and (\ref{frame2}), with $\mathbb{B} (x; \lambda)$
as specified in this section.

\subsection{Sturm-Liouville Systems} \label{sturm-liouville}

As an important special case of the general 
family of systems discussed in Section \ref{block-matrix},
we consider the Sturm-Liouville system
\begin{equation} \label{sl-eqn}
- (P(x) \phi')' + V(x) \phi = \lambda Q (x) \phi,
\end{equation} 
with boundary conditions
\begin{equation} \label{sl-bc}
\begin{aligned}
\alpha_1 \phi(0) + \alpha_2 P(0) \phi'(0) &= 0 \\
\beta_1 \phi(1) + \beta_2 P(1) \phi'(1) &= 0. 
\end{aligned}
\end{equation} 
Here, $\phi(x) \in \mathbb{C}^n$, and our notational
convention is to take $\alpha = (\alpha_1 \,\, \alpha_2) \in \mathbb{C}^{2n \times n}$
and $\beta = (\beta_1 \,\, \beta_2) \in \mathbb{C}^{2n \times n}$, with 
$\alpha$ and $\beta$ satisfying {\bf (BC1)}. We assume 
$P \in AC ([0,1], \mathbb{C}^{n \times n})$,
$V, Q \in L^1 ((0,1), \mathbb{C}^{n \times n})$, and that all 
three matrices are self-adjoint for a.e. $x \in (0, 1)$. In addition, we assume that 
$P(x)$ is invertible for each $x \in [0,1]$, and that  
$Q (x)$ is positive definite for a.e. $x \in (0,1)$. 

We can think of this system in terms of the operator
\begin{equation*}
\mathcal{L}_s \phi := Q(x)^{-1} \{- (P(x) \phi')' + V(x) \phi \},
\end{equation*}
with which we associate the domain
\begin{equation*}
\mathcal{D} (\mathcal{L}_s) 
= \{\phi \in L^2 ((0,1), \mathbb{C}^{n}): \phi, \phi' \in AC ([0,1], \mathbb{C}^n), 
\mathcal{L}_s \phi \in L^2 ((0,1), \mathbb{C}^n), (\ref{sl-bc}) \, \text{holds}\},
\end{equation*}
and the inner product 
\begin{equation*}
\langle \phi, \psi \rangle_{Q} := \int_0^1 (Q(x) \phi (x), \psi (x))_{\mathbb{C}^n} dx.
\end{equation*}
With this choice of domain and inner product, $\mathcal{L}_s$ is 
densely defined, closed, and self-adjoint, 
so $\sigma (\mathcal{L}_s) \subset \mathbb{R}$. I.e., we can take the 
interval $I$ associated with (\ref{hammy}) to be $I = \mathbb{R}$. 

For each $x \in [0,1]$, we define a new vector $y(x) \in \mathbb{C}^{2n}$
so that $y (x) = (y_1 (x) \, \, y_2 (x))^t$, with $y_1 (x) = \phi (x)$
and $y_2 (x) = P(x) \phi' (x)$. In this way, we express (\ref{sturm-liouville})
in the form 
\begin{equation*} \label{sturm-liouville2}
\begin{aligned}
&y' = \mathbb{A} (x; \lambda) y; \quad
\mathbb{A} (x; \lambda) = 
\begin{pmatrix}
0 & P(x)^{-1} \\
V(x) - \lambda Q(x) & 0 
\end{pmatrix}, \\
&\alpha y (0) = 0; \quad \beta y(1) = 0.
\end{aligned}
\end{equation*}
Upon multiplying both sides of this equation by $J$, 
we obtain (\ref{hammy}) with 
\begin{equation*}
\mathbb{B} (x; \lambda) = 
\begin{pmatrix}
\lambda Q(x) - V(x) & 0 \\
0 & P(x)^{-1}
\end{pmatrix}.
\end{equation*}
It is clear that $\mathbb{B} (x; \lambda)$ satifies our basic
assumptions {\bf (A)}. We check that $\mathbb{B} (x; \lambda)$ also 
satisfies Assumptions {\bf (B1)} and {\bf (B2)}. 

First, for {\bf (B1)}, we compute  
\begin{equation*}
\mathbb{B}_{\lambda} (x; \lambda) = 
\begin{pmatrix}
Q(x) & 0 \\
0 & 0
\end{pmatrix},
\end{equation*}
so that 
\begin{equation*}
\int_0^1 \mathbf{X}_1 (x;\lambda)^* \mathbb{B}_{\lambda} (x; \lambda) \mathbf{X}_1 (x; \lambda) dx 
= \int_0^1 X_1 (x;\lambda)^* Q(x) X_1 (x; \lambda) dx, 
\end{equation*}
which is clearly non-negative (since $Q$ is positive definite), 
and moreover it cannot have $0$ as an eigenvalue, because the associated
eigenvector $v \in \mathbb{C}^n$ would necessarily satisfy 
$X_1 (x; \lambda) v = 0$ for all $x \in [0,1]$, and this would contradict
linear independence of the columns of $X_1 (x; \lambda)$ (as solutions of 
(\ref{sl-eqn})). 

For {\bf (B2)}, we fix any $\lambda_1, \lambda_2 \in \mathbb{R}$,
$\lambda_1 < \lambda_2$, and notice that  
\begin{equation*}
\mathbb{B} (x; \lambda_2) - \mathbb{B} (x; \lambda_1)
= (\lambda_2 - \lambda_1) \mathbb{B}_{\lambda} (x; \lambda),
\end{equation*}
which is clearly non-negative. In addition, the same argument used to verify 
{\bf (B1)} shows that the condition assumed in Claim \ref{B2verification}
holds. Assumption {\bf (B2)} follows.

We conclude from Theorem \ref{bc1_theorem} that if 
$\mathcal{N} ([\lambda_1, \lambda_2); \mathcal{L}_s)$
denotes the spectral count for $\mathcal{L}_s$, we have 
\begin{equation*}
\mathcal{N} ([\lambda_1, \lambda_2); \mathcal{L}_s) 
= \mathcal{N}_{(0,1]} (\mathbf{X}_1 (\cdot; \lambda_1)^* J \mathbf{X}_2 (\cdot; \lambda_2)), 
\end{equation*}
where $\mathbf{X}_1 (x; \lambda_1)$ and $\mathbf{X}_2 (x; \lambda_2)$ denote
the frames specified in (\ref{frame1}) and (\ref{frame2}), with $\mathbb{B} (x; \lambda)$
as specified in this section. 

\begin{figure}[ht] 
\begin{center}\includegraphics[%
  width=6.5cm,
  height=5cm]{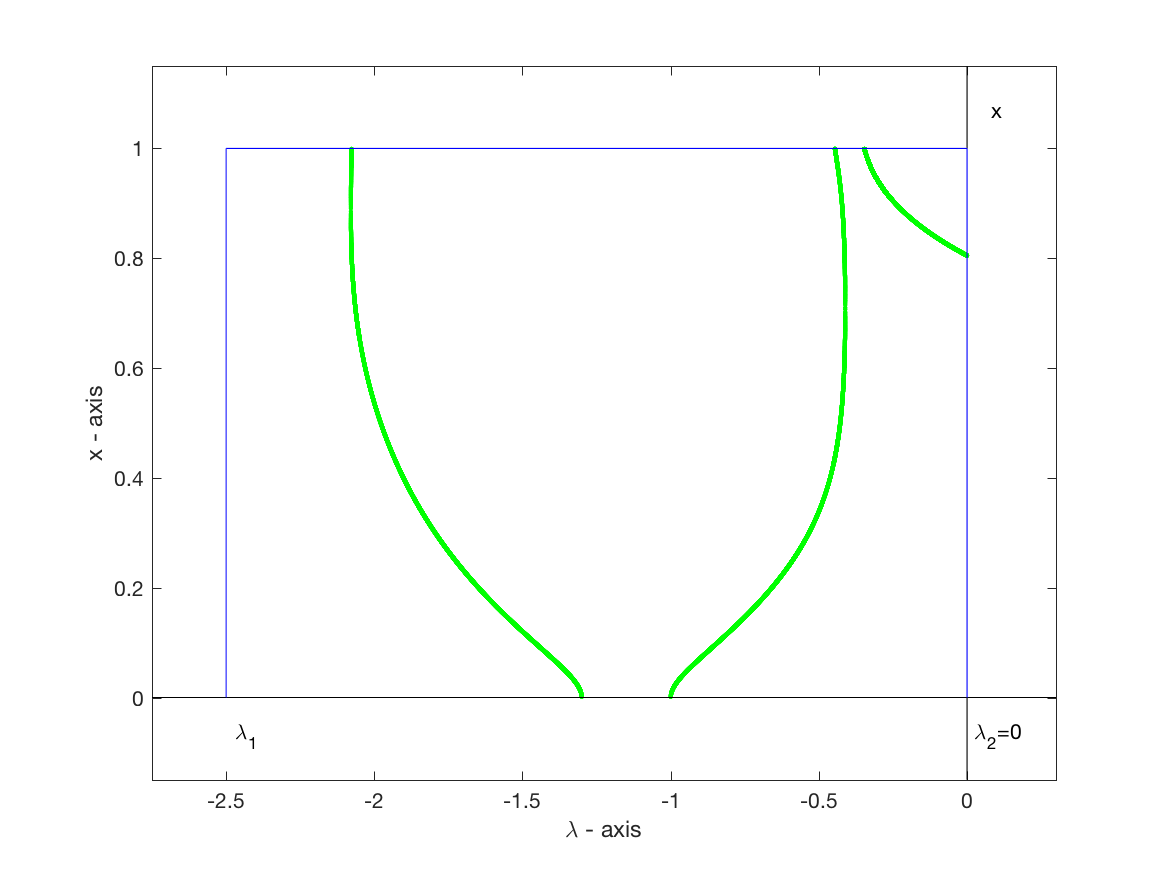}
\includegraphics[%
  width=6.5cm,
  height=5cm]{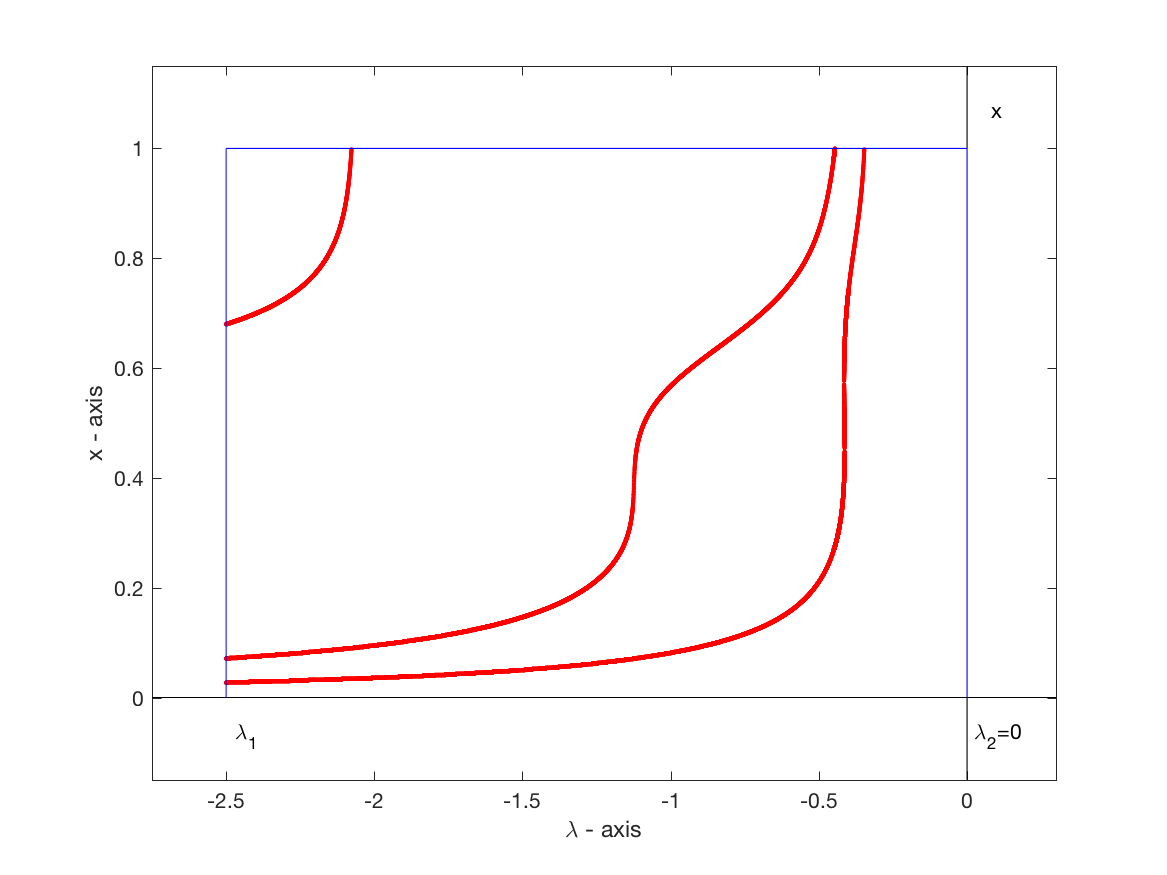}
\end{center}
\caption{Spectral curves for the Sturm-Liouville system example: approach of \cite{HJK2018, HS2016} on left; 
renormalized approach on right. \label{sl_figure}}
\end{figure}

As a specific example in this case, we consider (\ref{sl-eqn}) 
with $P = I_2$, $Q = 9 I_2$, 
\begin{equation} \label{example2}
V(x) = 
\begin{pmatrix}
-2.7 & - 18\sin(3x) + .0081x^2  \\
- 18\sin(3x) + .0081x^2 & 0 
\end{pmatrix},
\end{equation}
and boundary conditions {\bf (BC1)} specified by 
$\alpha = (\frac{1}{\sqrt{2}}I_2 \,\, \frac{1}{3\sqrt{2}}I_2)$ and 
$\beta = (\frac{1}{\sqrt{2}}I_2 \,\, \frac{1}{3\sqrt{2}}I_2)$.
Spectral curves for this equation are depicted in Figure \ref{sl_figure},
with the approach of \cite{HJK2018, HS2016} on the left and the renormalized 
approach on the right. As in our example for Dirac equations, 
several things are worth noting about this comparison of images:
(1) for the figure on the left, we see that the middle spectral
curve is non-monotonic, while for the figure on the right, 
all three spectral curves are monotonic; (2) as in Figure \ref{dirac_figure},
we see that in the setting of \cite{HJK2018, HS2016} spectral curves can 
emerge from any of the lower three shelves (bottom and right
in this case), while in the renormalized setting they
can only emerge from the left shelf; and (3) for the figure on the right, 
we see that spectral curves in the renormalized setting can almost 
become vertical. Regarding
Item (3), we verify in our proof of Theorem \ref{bc1_theorem} that 
the spectral curves cannot be vertical along any interval $\lambda \times (a,b)$
for $0 \le a < b \le 1$, and indeed this is also clear from the 
analysis of \cite{GBT1996} in which the authors show that there can only 
be a finite number of crossing points along any vertical shelf 
(see also \cite{GZ2017} for the same result in a more general setting).   

\subsection{Differential-Algebraic Sturm-Liouville Systems} \label{da_section}

We consider systems 
\begin{equation} \label{da-equation}
\mathcal{L}_a \phi 
= - (P(x) \phi')' + V(x) \phi = \lambda \phi, 
\end{equation} 
with degenerate matrices 
\begin{equation*}
P (x) = 
\begin{pmatrix}
P_{11} (x) & 0 \\
0 & 0
\end{pmatrix}.
\end{equation*}
Here, for some $0 < m < n$, $P_{11} \in AC ([0,1], \mathbb{C}^{m \times m})$ 
is a map into the space of self-adjoint matrices. We assume $P_{11} (x)$ is 
invertible for all $x \in [0, 1]$, and additionally that
$V \in C ([0,1], \mathbb{C}^{n \times n})$.
For notational convenience, we will write 
\begin{equation*}
V (x) = 
\begin{pmatrix}
V_{11} (x) & V_{12} (x) \\
V_{12} (x)^* & V_{22} (x)
\end{pmatrix},
\end{equation*}
where for each $x \in [0,1]$, $V_{11} (x)$ is a 
self-adjoint $m \times m$ matrix, 
$V_{12} (x)$ is an $m \times (n-m)$
matrix, and $V_{22} (x)$ is a self-adjoint $(n-m) \times (n-m)$
matrix. We will write 
\begin{equation*}
\phi = 
\begin{pmatrix}
\phi_1 \\
\phi_2
\end{pmatrix};
\quad \phi_1 (x; \lambda) \in \mathbb{C}^m; 
\quad \phi_2 (x; \lambda) \in \mathbb{C}^{n-m}, 
\end{equation*}
allowing us to express the system as 
\begin{equation*}
\begin{aligned}
- (P_{11} (x) \phi_1')' + V_{11} (x) \phi_1 + V_{12} (x) \phi_2 &= \lambda \phi_1 \\
V_{12} (x)^* \phi_1 + V_{22} (x) \phi_2 &= \lambda \phi_2.  
\end{aligned}
\end{equation*}

We place separated, self-adjoint boundary conditions on the 
first $m$ components, 
\begin{equation} \label{da-bc}
\begin{aligned}
\alpha_1 \phi_1 (0) + \alpha_2 P_{11} (0) \phi_1' (0) &= 0 \\
\beta_1 \phi_1 (1) + \beta_2 P_{11} (1) \phi_1' (1) &= 0,
\end{aligned}
\end{equation}
with $\alpha = (\alpha_1 \,\, \alpha_2)$ and $\beta = (\beta_1 \,\, \beta_2)$
satisfying the same assumptions as in Section \ref{sturm-liouville},
except with $n$ replaced by $m$. We specify the domain
\begin{equation*}
\begin{aligned}
\mathcal{D} (\mathcal{L}_a) &= \{ \phi = (\phi_1, \phi_2) \in L^2 ((0,1), \mathbb{C}^m) \times L^2 ((0,1), \mathbb{C}^{n-m}): 
\phi_1, \phi_1' \in AC ([0,1], \mathbb{C}^m), \\ 
&(\ref{da-bc}) \, \text{holds}, \, \mathcal{L}_a \phi \in L^2 ((0,1),\mathbb{C}^m) \times L^2 ((0,1),\mathbb{C}^{n-m})\},
\end{aligned}
\end{equation*}
and note that with this domain, it is straightfoward to verify that $\mathcal{L}_a$ 
is densely defined (in $L^2 ((0,1), \mathbb{C}^m) \times L^2 ((0,1), \mathbb{C}^{n-m})$), 
closed, and self-adjoint. 
In addition, we can apply Theorem 2.2 of \cite{ALMS1994} to see that 
the essential spectrum of $\mathcal{L}_a$ is precisely the union of 
the ranges of the eigenvalues of $V_{22} (x)$ as $x$ ranges over 
$[0,1]$. More precisely, let $\{\nu_k (x)\}_{k=1}^{n-m}$ denote the 
eigenvalues of $V_{22} (x)$, and let $\mathcal{R}_k$ denote the 
range of $\nu_k : [0,1] \to \mathbb{R}$. Then
\begin{equation*}
\sigma_{\operatorname{ess}} (\mathcal{L}_a) = 
\bigcup_{k=1}^{n-m} \mathcal{R}_k.
\end{equation*}

Now, fix any $\lambda_1 < \lambda_2$ so that 
$[\lambda_1, \lambda_2] \cap \sigma_{\operatorname{ess}} (\mathcal{L}_a) = \emptyset$,
and take any $\lambda \in [\lambda_1, \lambda_2]$. Then we can write 
\begin{equation*}
\phi_2 (x; \lambda) = (\lambda I - V_{22} (x))^{-1} V_{12} (x)^* \phi_1 (x; \lambda).
\end{equation*}
Upon substitution of $\phi_2$ into the equations for $\phi_1$, we obtain 
\begin{equation*}
- (P_{11} (x) \phi_1')' + V_{11} (x) \phi_1 + V_{12} (x) (\lambda I - V_{22} (x))^{-1} V_{12} (x)^* \phi_1  
= \lambda \phi_1. 
\end{equation*} 
We can express this system as a first-order system in the usual way, writing
$y_1 = \phi_1$ and $y_2 = P_{11} \phi_1'$, so that 
\begin{equation*}
y' = \mathbb{A} (x; \lambda) y,
\end{equation*}
with 
\begin{equation*}
\mathbb{A} (x; \lambda) =
\begin{pmatrix}
0 & P_{11} (x)^{-1} \\
\mathbf{V} (x; \lambda) - \lambda I & 0 
\end{pmatrix};
\quad 
\mathbf{V} (x; \lambda) = V_{11} + V_{12} (\lambda I - V_{22})^{-1} V_{12}^*. 
\end{equation*}
We multiply by $J$ to obtain the usual Hamiltonian form (\ref{hammy}) with 
\begin{equation*}
\mathbb{B} (x; \lambda) =
\begin{pmatrix}
\lambda I - \mathbf{V} (x; \lambda) & 0 \\
0 & P_{11} (x)^{-1} 
\end{pmatrix}.
\end{equation*}

It's clear from our Assumptions on $P(x)$ and $V(x)$ that the 
first part of Assumptions {\bf (A)} (addressing 
$\mathbb{B} (x; \lambda)$, not $\mathbb{B}_{\lambda} (x; \lambda)$)
holds in this case for any closed interval $I$ so that 
$I \cap \sigma_{\ess} (\mathcal{L}_a) = \emptyset$. 
In order to verify the second part of Assumptions {\bf (A)} and 
also Assumption {\bf (B1)}, we first compute 
\begin{equation*}
\mathbb{B}_{\lambda} (x; \lambda)
= 
\begin{pmatrix}
I - \mathbf{V}_{\lambda} (x; \lambda) & 0 \\
0 & 0
\end{pmatrix},
\end{equation*}
where  
\begin{equation*}
\mathbf{V}_{\lambda} (x; \lambda) = - V_{12} (x) (\lambda I - V_{22} (x))^{-2} V_{12} (x)^*.
\end{equation*} 
We see immediately that the second part of Assumptions {\bf (A)} holds for any
closed interval $I$ as described just above. In addition, 
since $(\lambda I - V_{22} (x))^{-1}$ is self-adjoint for all $x \in [0,1]$, 
we can express $\mathbf{V}_{\lambda} (x; \lambda)$ as 
\begin{equation*}
\mathbf{V}_{\lambda} (x; \lambda) = -  ((\lambda I - V_{22} (x))^{-1} V_{12} (x)^*)^* ((\lambda I - V_{22} (x))^{-1} V_{12} (x)^*),
\end{equation*} 
which is negative definite as long as $V_{12} (x)$ has trivial kernel and non-positive 
in any case. We see that $I - \mathbf{V}_{\lambda} (x; \lambda)$ is positive definite, and 
monotonicity (i.e., {\bf (B1)}) now follows in almost precisely the same way as in 
Section \ref{sturm-liouville}. 

Turning to {\bf (B2)}, we fix any $\lambda_1$ and $\lambda_2$ as described above,
and observe that 
\begin{equation*}
\mathbb{B} (x; \lambda_2) - \mathbb{B} (x; \lambda_1)
= 
\begin{pmatrix}
(\lambda_2 - \lambda_1) I - \mathbf{V} (x; \lambda_2) + \mathbf{V} (x; \lambda_1) & 0 \\
0 & 0 
\end{pmatrix}.
\end{equation*} 
Computing directly, we see that 
\begin{equation*}
\begin{aligned}
- \mathbf{V} (x; \lambda_2) &+ \mathbf{V} (x; \lambda_1)
= V_{12} (x) \{- (\lambda_2 I - V_{22})^{-1} + (\lambda_1 I - V_{22})^{-1} \} V_{12} (x)^* \\
&= (\lambda_2 - \lambda_1) V_{12} (x) \{(\lambda_1 I - V_{22})^{-1} (\lambda_2 I - V_{22})^{-1}\} V_{12} (x)^*. 
\end{aligned}
\end{equation*}
The matrix in curved brackets is self-adjoint, and by spectral mapping, its 
eigenvalues are 
\begin{equation} \label{daev}
\Big{\{} \frac{1}{(\lambda_1 - \nu_k (x)) (\lambda_2 - \nu_k (x))} \Big{\}}_{k=1}^{n-m}.
\end{equation}
By our assumption that $[\lambda_1, \lambda_2] \cap \sigma_{\operatorname{ess}} = \emptyset$, 
we see that for each $k \in \{1, 2, \dots, n-m\}$ and each $x \in [0,1]$, we either have 
$\lambda_1 < \lambda_2 < \nu_k (x)$ or we have 
$\nu_k (x) < \lambda_1 < \lambda_2$. (The idea is simply that $\nu_k (x)$ cannot 
lie between $\lambda_1$ and $\lambda_2$.) In either case, 
the eigenvalues (\ref{daev}) are all positive, verifying that the self-adjoint
matrix $(\lambda_1 I - V_{22})^{-1} (\lambda_2 I - V_{22})^{-1}$ is positive 
definite. It follows immediately that 
$(\lambda_2 - \lambda_1) I - \mathbf{V} (x; \lambda_2) + \mathbf{V} (x; \lambda_1)$
is positive definite, and $\mathbb{B} (x; \lambda_2) - \mathbb{B} (x; \lambda_1)$
is non-negative. In addition, an argument similar to our verification of 
{\bf (B1)} in this case serves to verify that the assumptions of Claim
\ref{B2verification} hold in this case, and the moreover part of 
{\bf (B2)} follows. 

We conclude from Theorem \ref{bc1_theorem} that if 
$\mathcal{N} ([\lambda_1, \lambda_2); \mathcal{L}_a)$
denotes the spectral count for $\mathcal{L}_a$, we have 
\begin{equation*}
\mathcal{N} ([\lambda_1, \lambda_2); \mathcal{L}_a) 
= \mathcal{N}_{(0,1]} (\mathbf{X}_1 (\cdot; \lambda_1)^* J \mathbf{X}_2 (\cdot; \lambda_2)), 
\end{equation*}
where $\mathbf{X}_1 (x; \lambda_1)$ and $\mathbf{X}_2 (x; \lambda_2)$ denote
the frames (\ref{frame1}) and (\ref{frame2}), with $\mathbb{B} (x; \lambda)$
as specified in this section. 

As a specific example in this case, we consider (\ref{da-equation}) with 
\begin{equation*}
P = 
\begin{pmatrix}
I_2 & 0_2 \\
0_2 & 0_2 
\end{pmatrix},
\end{equation*}
\begin{equation*}
V(x) =
\begin{pmatrix}
-8 - \frac{.7 \cos(6 \pi x)}{2+\cos(6 \pi x)} & - \frac{\cos (\pi x)}{2 + \cos (4 \pi x)} & 1 & 0 \\
- \frac{\cos (\pi x)}{2 + \cos (4 \pi x)} & 1 & 0 & 1 \\
1 & 0 & 1 - .8x \sin (x) & 0 \\
0 & 1 & 0 & 1-.8x \sin (x)
\end{pmatrix},
\end{equation*}
and Neumann boundary conditions specified by 
$\alpha = (0_2 \,\, I_2)$ and $\beta = (0_2 \,\, I_2)$.

\begin{figure}[ht] 
\begin{center}\includegraphics[%
  width=6.5cm,
  height=5cm]{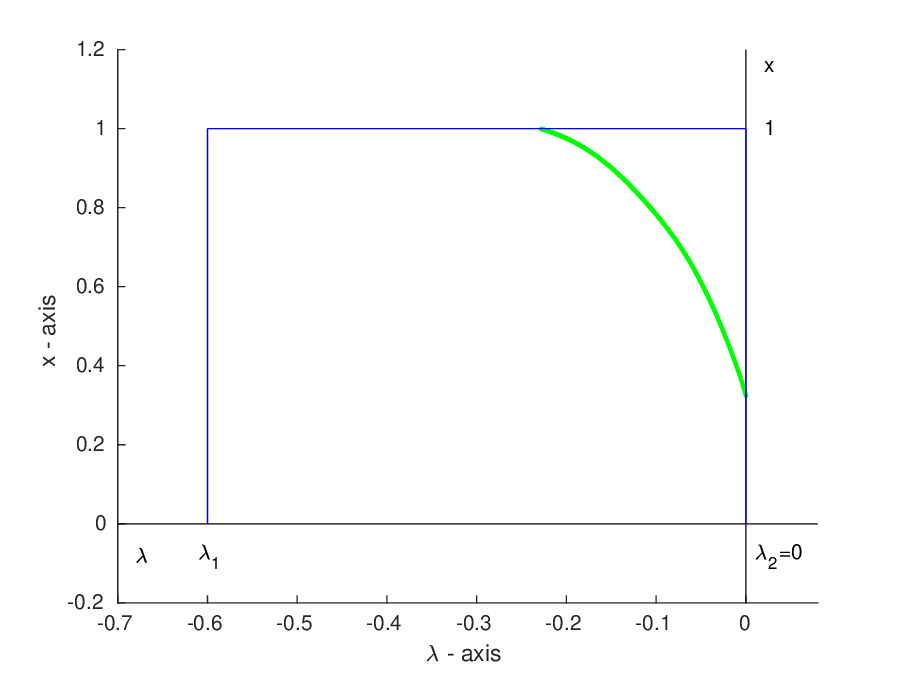}
\includegraphics[%
  width=6.5cm,
  height=5cm]{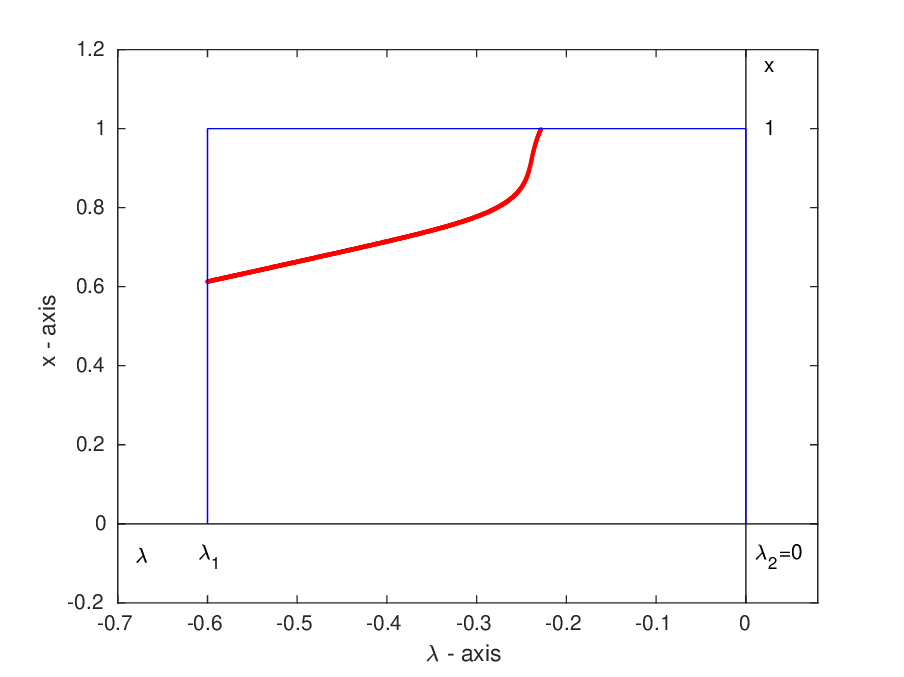}
\end{center}
\caption{Spectral curves for the differential-algebraic equation example: approach of \cite{HJK2018, HS2016} on left; 
renormalized approach on right. \label{da_figure}}
\end{figure}

Spectral curves are depicted in Figure \ref{da_figure} for this example, 
with the approach of \cite{HJK2018, HS2016} on the left and the renormalized approach 
on the right. In this case, 
\begin{equation*}
V_{22} (x) = 
\begin{pmatrix}
1 - .8x \sin (x) & 0 \\
0 & 1 - 8x \sin (x)
\end{pmatrix},
\end{equation*}
so the essential spectrum is confined to the range of 
$(1 - .8x \sin (x)) |_{[0,1]} = [1-.8\sin(1), 1]$,
well to the right of our depicted Maslov boxes.

\subsection{Fourth-Order Equations}
\label{fourth-section}

We consider fourth-order systems 
\begin{equation} \label{fourth-equation}
    \phi'''' - (V_2 (x) \phi')' + V_0 (x) \phi = \lambda \phi,
\end{equation}
with periodic boundary conditions 
\begin{equation} \label{fourth-bc}
    \phi^{(k)} (1) = \phi^{(k)} (0), \quad k = 0, 1, 2, 3.
\end{equation}
Here, $\phi(x) \in \mathbb{C}^n$, and we take 
$V_0 \in L^1 ((0,1), \mathbb{C}^{n \times n})$
and 
$V_2 \in AC ([0,1], \mathbb{C}^{n \times n})$,
with both matrix functions self-adjoint for a.e.
$x \in (0,1)$, and $V_2 (0) = V_2 (1)$. (This final 
condition isn't necessary, and is used only to 
retain periodicity of the boundary conditions for 
the first-order system in the form we'll express 
it.)

We can think of this system in terms of the operator
\begin{equation*}
\mathcal{L}_f \phi := \phi'''' - (V_2 (x) \phi')' + V_0 (x) \phi,
\end{equation*}
with which we associate the domain
\begin{equation*}
\begin{aligned}
\mathcal{D} (\mathcal{L}_f) 
&= \{\phi \in L^2 ((0,1), \mathbb{C}^{n}): \phi^{(k)} \in AC ([0,1], \mathbb{C}^n), k = 0, 1, 2, 3, \\
& \quad \quad \quad \mathcal{L}_f \phi \in L^2 ((0,1), \mathbb{C}^n), (\ref{fourth-bc}) \,\, \text{holds}\}.
\end{aligned}
\end{equation*}
With this choice of domain and inner product, $\mathcal{L}_f$ is 
densely defined, closed, and self-adjoint, 
so $\sigma (\mathcal{L}_f) \subset \mathbb{R}$. I.e., we can take the 
interval $I$ associated with (\ref{hammy}) to be $I = \mathbb{R}$. 
(See, e.g., \cite{Weidmann1987}.)

For each $x \in [0,1]$, we define a new vector $y(x) \in \mathbb{C}^{4n}$
so that 
\begin{equation*}
y (x) = (y_1 (x) \, \, y_2 (x) \,\, y_3 (x) \, \, y_4 (x))^T,     
\end{equation*}
with $y_1 (x) = \phi (x)$, $y_2 (x) = \phi'' (x)$, $y_3 (x) = - \phi'''(x) + V_2 (x) \phi'(x)$,
and $y_4 (x) = - \phi' (x)$. In this way, we express (\ref{fourth-equation})
in the form 
\begin{equation*} \label{fourth-equation2}
\begin{aligned}
&y' = \mathbb{A} (x; \lambda) y; \quad
\mathbb{A} (x; \lambda) = 
\begin{pmatrix}
0 & 0 & 0 & -I_n \\
0 & 0 & -I_n & -V_2 (x) \\
V_0 (x) - \lambda I_n & 0 & 0 & 0 \\
0 & -I_n & 0 & 0 \\
\end{pmatrix}, \\
&\Theta \begin{pmatrix}
y(0) \\ y(1)
\end{pmatrix} = 0, \quad \Theta = (I_{4n} \,\,\, -I_{4n}).
\end{aligned}
\end{equation*}
(See \cite{HJK2018} for a discussion of the rationale for 
this choice in defining $y$.)
Upon multiplying both sides of this equation by $J$, 
we obtain (\ref{hammy}) with 
\begin{equation*}
\mathbb{B} (x; \lambda) = 
\begin{pmatrix}
\lambda I_n - V_0 (x) & 0 & 0 & 0 \\
0 & I_n & 0 & 0 \\
0 & 0 & 0 & -I_n \\
0 & 0 & -I_n & -V_2 (x) \\
\end{pmatrix}.
\end{equation*}
It is clear that $\mathbb{B} (x; \lambda)$ satifies our basic
assumptions {\bf (A)}. We check that $\mathbb{B} (x; \lambda)$ also 
satisfies Assumptions {\bf (B1)$\mathbf{'}$} and {\bf (B2)$\mathbf{'}$}. 

First, for {\bf (B1)$\mathbf{'}$}, we compute  
\begin{equation*}
\mathbb{B}_{\lambda} (x; \lambda) = 
\begin{pmatrix}
I_n & 0 & 0 & 0 \\
0 & 0 & 0 & 0 \\
0 & 0 & 0 & 0 \\
0 & 0 & 0 & 0 \\
\end{pmatrix}.
\end{equation*}
so that 
\begin{equation*}
\int_0^1 \Phi (x; \lambda)^* \mathbb{B}_{\lambda} (x; \lambda) \Phi (x; \lambda) dx 
= \int_0^1 \Phi_1 (x;\lambda)^* \Phi_1 (x; \lambda) dx, 
\end{equation*}
where $\Phi_1 (x; \lambda)$ denotes the $n \times 4n$ matrix comprising 
the first $n$ rows of the $4n \times 4n$ fundamental matrix $\Phi (x; \lambda)$.
The matrix $\int_0^1 \Phi_1 (x;\lambda)^* \Phi_1 (x; \lambda) dx$ 
is clearly non-negative, 
and moreover it cannot have $0$ as an eigenvalue, because the associated
eigenvector $v \in \mathbb{C}^{4n}$ would necessarily satisfy 
$\Phi_1 (x; \lambda) v = 0$ for all $x \in [0,1]$, and this would contradict
linear independence of the columns of $\Phi_1 (x; \lambda)$ (as solutions of 
(\ref{fourth-equation})). 

For {\bf (B2)$\mathbf{'}$}, we fix any $\lambda_1, \lambda_2 \in \mathbb{R}$,
$\lambda_1 < \lambda_2$, and notice that  
\begin{equation*}
\mathbb{B} (x; \lambda_2) - \mathbb{B} (x; \lambda_1)
= (\lambda_2 - \lambda_1) \mathbb{B}_{\lambda} (x; \lambda),
\end{equation*}
which is clearly non-negative. In addition, the same argument used to verify 
{\bf (B1)$\mathbf{'}$} shows that the condition assumed in Claim \ref{B2verification}
holds. Assumption {\bf (B2)$\mathbf{'}$} follows.

We conclude from Theorem \ref{bc2_theorem} that if 
$\mathcal{N} ([\lambda_1, \lambda_2); \mathcal{L}_f)$
denotes the spectral count for $\mathcal{L}_f$, we have 
\begin{equation*}
\mathcal{N} ([\lambda_1, \lambda_2); \mathcal{L}_f) 
= \mathcal{N}_{(0,1]} (\mathbf{X}_3 (\cdot; \lambda_1)^* J \mathbf{X}_4 (\cdot; \lambda_2)), 
\end{equation*}
where $\mathbf{X}_3 (x; \lambda_1)$ and $\mathbf{X}_4 (x; \lambda_2)$ denote
the frames specified respectively in (\ref{frame3}) and (\ref{frame4}), 
with $\mathbb{B} (x; \lambda)$ as in this section. 

\begin{figure}[ht] 
\begin{center}\includegraphics[%
  width=6.5cm,
  height=5cm]{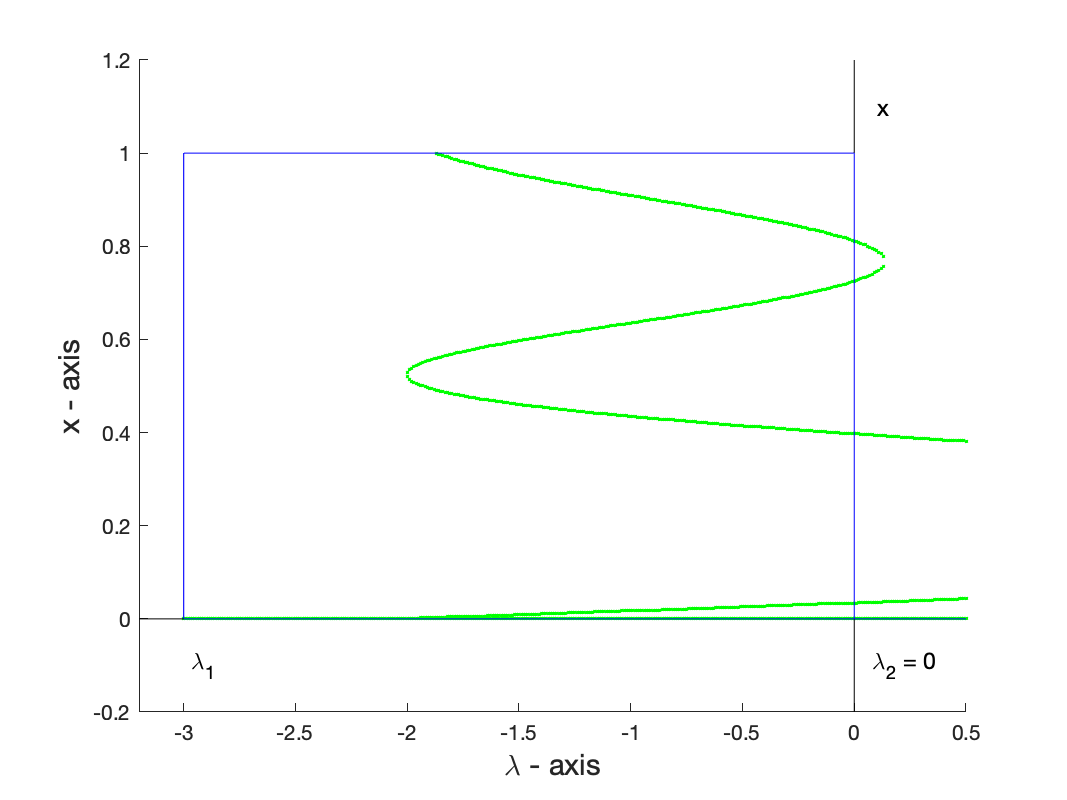}
\includegraphics[%
  width=6.5cm,
  height=5cm]{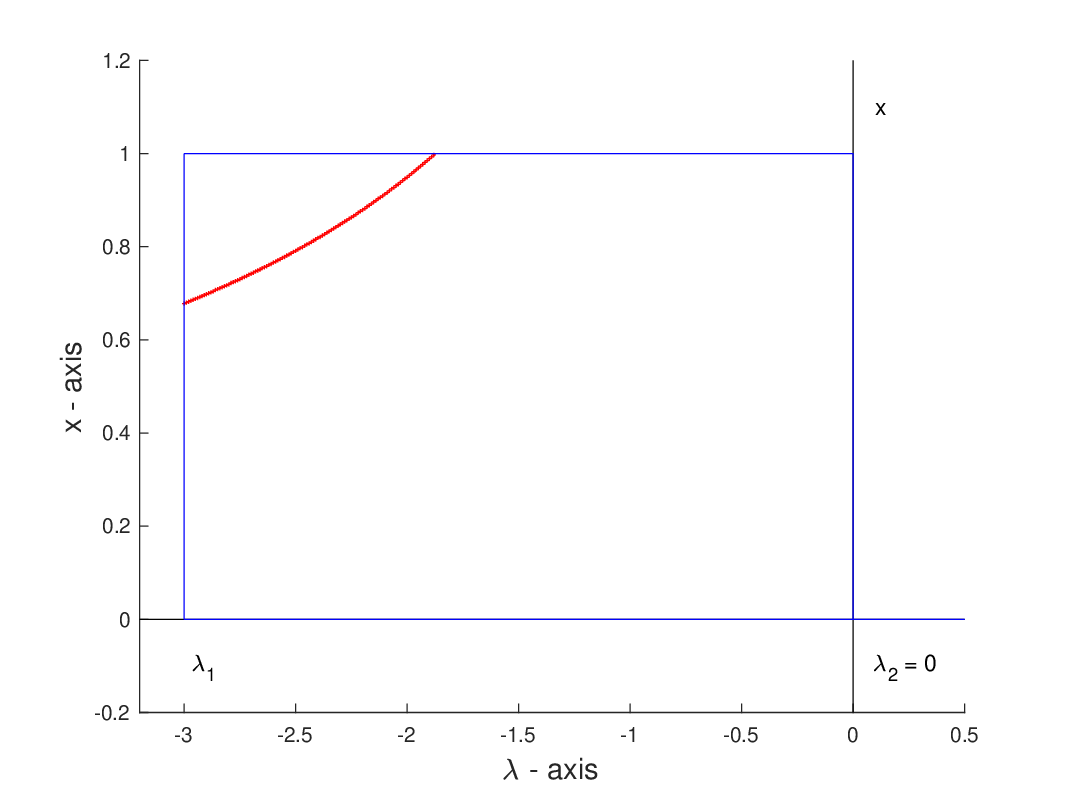}
\end{center}
\caption{Spectral curves for the Sturm-Liouville system example: approach of \cite{HJK2018, HS2016} on left; 
renormalized approach on right. \label{fourth-figure}}
\end{figure}

As a specific example in this case, we consider (\ref{fourth-equation})
with $n = 1$ and 
\begin{equation*}
    \begin{aligned}
    V_0 (x) &= -2 + 10 \sin(12x) \\
    V_2 (x) &= 10 \cos(2\pi x).
    \end{aligned}
\end{equation*}
The boundary conditions have the form {\bf (BC2)}
with $\Theta = (I_4 \,\,\, -I_4)$. 
Spectral curves for this equation are depicted in Figure \ref{fourth-figure},
with the approach of \cite{HJK2018, HS2016} on the left and the renormalized 
approach on the right. For the approach of \cite{HJK2018, HS2016}, 
each point on the bottom shelf is a crossing point, and in addition 
a spectral curve emerges from the bottom shelf. 

\section*{Appendix}

In this short appendix, we briefly discuss the view of our 
operator 
\begin{equation} \label{operator-pencil}
\mathcal{L} (\lambda) := J \frac{d}{dx} - \mathbb{B} (\cdot; \lambda)
\end{equation}
as an operator pencil. In order to keep the discussion brief, we
focus on the case of boundary conditions {\bf (BC1)}, for which the 
domain of $\mathcal{L} (\lambda)$ can be taken to be 
\begin{equation*}
\begin{aligned}
\mathcal{D} (\mathcal{L} (\lambda)) &:= \{y \in L^2 ((0,1), \mathbb{C}^{2n}): 
y \in \AC ([0,1], \mathbb{C}^{2n}), \\ 
& \quad \quad \mathcal{L} y \in L^2 ((0,1), \mathbb{C}^{2n}),
\, \alpha y(0) = 0, \, \beta y(1) = 0 \}.
\end{aligned}
\end{equation*}
We will confine the discussion in this appendix to the case in which 
$\mathbb{B} (\cdot; \lambda) \in L^2 ((0,1), \mathbb{C}^{2n \times 2n})$
for all $\lambda \in I$. Under this additional assumption,  
$\mathcal{D} (\mathcal{L} (\lambda))$ is independent of $\lambda$,
and in order to emphasize this independence we will express 
$\mathcal{D} (\mathcal{L} (\lambda))$ as $\mathcal{D}$. 
Here, $I \subset \mathbb{R}$ continues to be 
the interval specified in the introduction containing all values $\lambda$ 
under consideration. 

Following the development of \cite{BLR-M2014}, we specify the resolvent
set of $\mathcal{L}$ as 
\begin{equation} \label{pencil-resolvent}
\rho (\mathcal{L}) := \{\lambda \in I: \mathcal{L} (\lambda)^{-1} \in \mathcal{B} (L^2 ((0,1), \mathbb{C}^{2n}))\},
\end{equation}
where $\mathcal{B} (L^2 ((0,1), \mathbb{C}^{2n}))$ denotes the linear space of all
bounded linear operators mapping $L^2 ((0,1), \mathbb{C}^{2n})$ to itself, 
and we specify the spectrum of $\mathcal{L}$ as 
$\sigma (\mathcal{L}) = I \backslash \rho (\mathcal{L})$. 
More generally, operator pencils are often defined on open sets 
of the complex plane $\Omega \subset \mathbb{C}$, but such 
a specification is not necessary for this brief discussion. 
In order to be precise about terminology, 
we define what we mean by the essential spectrum and
the point spectrum (adapted from \cite{KP2013}). For this, 
we assume, as in the current setting, that 
$\mathcal{D} := \dom (\mathcal{L} (\lambda))$
is independent of $\lambda$, and we denote by 
$\mathbb{L} (L^2 ((0,1),\mathbb{C}^{2n}))$ the 
space of all closed linear operators 
mapping $\mathcal{D} \subset L^2 ((0,1),\mathbb{C}^{2n})$ 
to $L^2 ((0,1),\mathbb{C}^{2n})$.

\begin{definition}
We define the essential spectrum $\sigma_{\ess} (\mathcal{L})$
of an operator pencil $\mathcal{L}: I \to \mathbb{L} (L^2 ((0,1),\mathbb{C}^{2n}))$
as the set of $\lambda \in I$ for which either $\mathcal{L} (\lambda)$
is not Fredholm or $\mathcal{L} (\lambda)$ is Fredholm with 
Fredholm index $\ind (\mathcal{L} (\lambda)) \ne 0$. We define 
the point spectrum $\sigma_{\pt} (\mathcal{L})$ as the set of 
$\lambda \in I$ so that $\ind (\mathcal{L} (\lambda)) = 0$, 
but $\mathcal{L} (\lambda)$ is not invertible. 
\end{definition}

With these definitions, we see that the sets
$\rho (\mathcal{L})$, $\sigma_{\ess} (\mathcal{L})$, 
and $\sigma_{\pt} (\mathcal{L})$ are mutually 
exclusive, and 
\begin{equation*}
I = \rho (\mathcal{L}) \cup \sigma_{\ess} (\mathcal{L}) \cup \sigma_{\pt} (\mathcal{L});
\quad \sigma (\mathcal{L}) = \sigma_{\ess} (\mathcal{L}) \cup \sigma_{\pt} (\mathcal{L}).
\end{equation*}
Another way to view the definitions is as follows. A value $\lambda_0 \in I$
is categorized as an element of $\rho (\mathcal{L})$, $\sigma_{\ess} (\mathcal{L})$, 
or $\sigma_{\pt} (\mathcal{L})$ according precisely to the categorization
of $0$ relative to the operator $\mathcal{L} (\lambda_0)$.

Returning to our particular operator pencil from (\ref{operator-pencil}),
it's a straightforward application of the methods of \cite{Weidmann1987} to verify 
that under our Assumptions {\bf (A)}, we have the following: 
for each $\lambda \in I$, $\mathcal{L} (\lambda)$ is Fredholm with
index zero, and indeed is self-adjoint. We can conclude that 
$\sigma (\mathcal{L})$ is comprised entirely of point spectrum, 
and in particular that for each $\lambda \in \sigma (\mathcal{L})$ there exist 
a finite number of linearly independent eigenfunctions 
$\{y_i (x; \lambda)\}_{i=1}^m \subset \mathcal{D}$ so that 
$\mathcal{L} (\lambda) y_i (\cdot; \lambda) = 0$ for all 
$i \in \{1, 2, \dots, m\}$. In addition, our Assumption {\bf (B1)}
ensures that the eigenvalues are all discrete (i.e., isolated).
Our Theorems \ref{bc1_theorem}
and \ref{bc2_theorem} count the number of such discrete
eigenvalues, including geometric multiplicity, and it's 
natural to consider how this relates to the same count 
using algebraic multiplicity. First, proceeding as in 
\cite{KM2014}, we can define the algebraic multiplicity 
of an eigenvalue $\lambda_0$ of $\mathcal{L}$ in terms
of the nature of the Jordan chains associated with it. 
Readers interested in a complete definition along these
lines can find it in Definition 6 of \cite{KM2014}, but for
our purposes, we only require the following.

\begin{definition} \label{algebraic-multiplicity}
Let $\lambda_0 \in I$ be an 
eigenvalue of an operator pencil $\mathcal{L}: I \to \mathbb{L} (L^2 ((0,1),\mathbb{C}^{2n}))$
with geometric multiplicity $m$, and assume 
$\mathcal{L}'(\lambda_0) \in \mathbb{L} (L^2 ((0,1),\mathbb{C}^{2n}))$ exists, 
with additionally $\dom (\mathcal{L}'(\lambda_0)) = \mathcal{D}$. 
Suppose that for any pair $(y, \zeta)$ with $y \in \ker \mathcal{L} (\lambda_0)$,
and $\zeta \in \dom (\mathcal{L} (\lambda_0))$ satisfying 
\begin{equation}
\mathcal{L} (\lambda_0) \zeta = \mathcal{L}' (\lambda_0) y,
\end{equation}
we must have $y \equiv 0$. Then $\lambda_0$ has algebraic 
multiplicity $m$.
\end{definition}  

We are now in a position to verify that under slightly 
stronger conditions on $\mathbb{B} (x; \lambda)$ than 
assumed for Theorems \ref{bc1_theorem} and \ref{bc2_theorem},
the geometric and algebraic multiplicies of eigenvalues of 
the operator pencil $\mathcal{L}: I \to \mathbb{L} (L^2 ((0,1),\mathbb{C}^{2n}))$
coincide. 

\begin{claim} \label{pencil-claim} 
Let Assumptions {\bf (A)} and {\bf (B1)} hold, and 
additionally assume that for all $\lambda \in I$,
we have $\mathbb{B} (\cdot; \lambda), \mathbb{B}_{\lambda} (\cdot; \lambda) \in L^2 ((0,1),\mathbb{C}^{2n \times 2n})$.
Then for any eigenvalue $\lambda_0$ of the operator 
pencil $\mathcal{L}$, geometric and algebraic multiplicities
agree. 
\end{claim}

\begin{proof}
In our setting, $\mathcal{L}' (\lambda) = \mathbb{B}_{\lambda} (x; \lambda)$. 
Suppose $\lambda_0$ is an eigenvalue of $\mathcal{L}$, and that for some 
$y (\cdot; \lambda_0) \in \ker \mathcal{L} (\lambda_0)$,
there is a corresponding $\zeta (\cdot; \lambda_0) \in \mathcal{D} (\mathcal{L})$
so that 
\begin{equation*}
\mathcal{L} (\lambda_0) \zeta (x; \lambda_0) 
= \mathbb{B}_{\lambda} (x; \lambda_0) y (x; \lambda_0), \quad \text{a.e. } x \in (0,1). 
\end{equation*}
(Our additional assumption $\mathbb{B}_{\lambda} (\cdot; \lambda) \in L^2 ((0,1),\mathbb{C}^{2n \times 2n})$
ensures that $\mathcal{L}' (\lambda_0)$ maps 
$\mathcal{D}$ to $L^2 ((0,1),\mathbb{C}^{2n})$, and in particular that 
$\mathbb{B}_{\lambda} (x; \lambda_0) y (x; \lambda_0)$ is in the range of $\mathcal{L} (\lambda_0)$.)
If we take an $L^2$ inner product of this equation with $y$, we obtain the 
relation 
\begin{equation*}
\langle \mathcal{L} (\lambda_0) \zeta, y \rangle
= \langle \mathbb{B}_{\lambda} (x; \lambda_0) y, y\rangle. 
\end{equation*}
Since $\mathcal{L} (\lambda_0)$ is self-adjoint, the left-hand side 
can be computed as 
\begin{equation*}
\langle \mathcal{L} (\lambda_0) \zeta, y \rangle
 = \langle \zeta, \mathcal{L} (\lambda_0) y \rangle = 0.
\end{equation*}
We see that the right-hand side satisfies 
$\langle \mathbb{B}_{\lambda} (\cdot; \lambda_0) y, y \rangle = 0$,
and by our positivity condition {\bf (B1)} this means $y = 0$
for a.e. $x \in (0,1)$. According to Definition 
\ref{algebraic-multiplicity}, we can conclude that the 
algebraic multiplicity of $\lambda_0$ agrees with the 
geometric multiplicity of $\lambda_0$. 
\end{proof}

\bigskip
{\it Acknowledgments.} The authors are grateful to Yuri Latushkin for bringing 
\cite{GZ2017} to their attention, for suggesting that the problem could be 
approached with the Maslov index, and for several helpful conversations along
the way. A.S. acknowledges support from the National Science Foundation under
grant DMS-1910820.

\end{document}